\documentclass[11pt]{amsart}
\usepackage{amsmath}
\usepackage{amsxtra}
\usepackage{amscd}
\usepackage{amsthm}
\usepackage{amsfonts}
\usepackage{amssymb}
\usepackage{eucal}

\newtheorem{thm}{Theorem}[section]
\newtheorem{conj}{Conjecture}
\newtheorem{cor}[thm]{Corollary}
\newtheorem{lem}[thm]{Lemma}
\newtheorem{prop}[thm]{Proposition}

\theoremstyle{remark}
\newtheorem{remark}[thm]{Remark}

\theoremstyle{definition}

\numberwithin{equation}{section}

\newcommand{\bean}{\begin{eqnarray}}
\newcommand{\eean}{\end{eqnarray}}
\newcommand{\be}{\begin{displaymath}}
\newcommand{\ee}{\end{displaymath}}
\newcommand{\bea}{\begin{eqnarray*}}   
\newcommand{\eea}{\end{eqnarray*}}

\newcommand{\thmref}[1]{Theorem~\ref{#1}}
\newcommand{\secref}[1]{Section~\ref{#1}}
\newcommand{\lemref}[1]{Lemma~\ref{#1}}
\newcommand{\propref}[1]{Proposition~\ref{#1}}
\newcommand{\corref}[1]{Corollary~\ref{#1}}
\newcommand{\remref}[1]{Remark~\ref{#1}}
\newcommand{\conjref}[1]{Conjecture~\ref{#1}}
\newcommand{\nc}{\newcommand}
\nc{\on}{\operatorname}
\nc{\ch}{\mbox{ch}}
\nc{\Z}{{\mathbb Z}}
\nc{\C}{{\mathbb C}}
\nc{\pone}{{\mathbb P}^1}
\nc{\pa}{\partial}
\nc{\F}{{\mathcal F}}
\nc{\arr}{\rightarrow}
\nc{\larr}{\longrightarrow}
\nc{\al}{\alpha}
\nc{\ri}{\rangle}
\nc{\lef}{\langle}
\nc{\W}{{\mathcal W}}
\nc{\la}{\lambda}
\nc{\ep}{\epsilon}
\nc{\su}{\widehat{{\mathfrak s}{\mathfrak l}}_2}
\nc{\sw}{{\mathfrak s}{\mathfrak l}}
\nc{\g}{{\mathfrak g}}
\nc{\h}{{\mathfrak h}}
\nc{\n}{{\mathfrak n}}
\nc{\N}{\widehat{\n}}
\nc{\G}{\widehat{\g}}
\nc{\De}{\Delta}
\nc{\gt}{\widetilde{\g}}
\nc{\Ga}{\Gamma}
\nc{\one}{{\mathbf 1}}
\nc{\z}{{\mathfrak Z}}
\nc{\La}{\Lambda}
\nc{\wt}{\widetilde}
\nc{\wh}{\widehat}
\nc{\cri}{_{\kappa_c}}
\nc{\kk}{h^\vee}
\nc{\sun}{\widehat{\sw}_N}
\nc{\si}{\sigma}
\nc{\el}{\ell}
\nc{\bi}{\bibitem}
\nc{\om}{\omega}
\nc{\ol}{\overline}
\nc{\ds}{\displaystyle}
\nc{\dzz}{\frac{dz}{z}}
\nc{\Res}{\on{Res}}
\nc{\mc}{\mathcal}
\nc{\Cal}{\mathcal}
\nc{\bb}{{\mathfrak b}}
\nc{\ot}{\otimes}
\nc{\R}{{\mc R}}
\nc{\yy}{{\mc Y}}
\nc{\ga}{\gamma}

\nc{\us}{\underset}
\nc{\opl}{\oplus}
\nc{\beq}{\begin{equation}}
\nc{\Fq}{{\mathcal F}}
\nc{\Mq}{{\mathcal M}}
\nc{\Rep}{\on{Rep}}
\nc{\sssec}{\subsubsection}
\nc{\ssec}{\subsection}
\nc{\lan}{\langle}
\nc{\ran}{\rangle}

\nc{\D}{\mathcal D}
\nc{\Vect}{\on{Vect}}
\nc{\ghat}{\G}
\nc{\T}{\mc T}
\nc{\Tloc}{\T^\g_{\on{loc}}}
\nc{\vac}{|0\ran}
\nc{\Wick}{{\mb :}}
\nc{\mb}{\mathbf}
\nc{\delz}{\partial_z}
\nc{\K}{{\cali K}}
\nc{\cali}{\mathcal}
\nc{\li}{\mathfrak l}
\nc{\lt}{\widetilde{\li}}
\nc{\astar}{a^*}
\nc{\cA}{{\mc A}}
\nc{\ka}{\kappa}

\nc{\OO}{{\mc O}}
\nc{\AutO}{\on{Aut}\OO}
\nc{\DerO}{\on{Der}\OO}
\nc{\DerpO}{\on{Der}_+\OO}
\nc{\Au}{{\mc A}ut}
\nc{\mf}{\mathfrak}
\nc{\V}{{\mathbb V}}
\nc{\hh}{\wh{\h}}

\nc{\pp}{{\mathfrak p}}
\nc{\mm}{{\mathfrak m}}
\nc{\rr}{{\mathfrak r}}
\nc{\ket}{\rangle}
\nc{\zz}{{\mathfrak z}}
\nc{\gr}{\on{gr}}
\nc{\Spe}{\on{Spec}}
\nc{\rv}{\crho}
\nc{\can}{\on{can}}
\nc{\CC}{\on{Op}_G(D))}
\nc{\Op}{\on{Op}_G(D)}
\nc{\MOp}{\on{MOp}_G(D)}
\nc{\Db}{{\mathbb D}}
\nc{\ww}{w}

\nc{\af}{{\mathbb A}^1}
\nc{\bs}{\backslash}
\nc{\laa}{(\la_i)}
\nc{\zn}{(z_i)}

\nc{\cla}{\check{\la}}
\nc{\cmu}{\check{\mu}}
\nc{\crho}{\check{\rho}}
\nc{\chal}{\check{\al}}
\nc{\cc}{{\mathfrak c}}

\begin{document}

\title[Opers on the projective line, flag manifolds and Bethe
Ansatz]{Opers on the projective line, flag manifolds and Bethe Ansatz}

\author{Edward Frenkel}\thanks{Partially supported by grants from the
  Packard Foundation and the NSF}

\address{Department of Mathematics, University of California,
  Berkeley, CA 94720, USA}

\dedicatory{To Boris Feigin on his 50th birthday}

\date{August 2003; Revised March 2005}

\maketitle


\section{Introduction}

\subsection{}

Our starting point is the Gaudin model associated to a simple
finite-dimensional Lie algebra $\g$. Let us introduce some
notation. For any integral dominant weight $\la$, denote by $V_\la$
the irreducible finite-dimensional representation of $\g$ of highest
weight $\la$. Let $z_1,\ldots,z_N$ be a set of distinct complex
numbers and $\la_1,\ldots,\la_N$ a set of dominant integral weights of
$\g$. Set
$$V_{\laa} = V_{\la_1} \otimes \ldots \otimes V_{\la_N}.$$ Let $\{ J_a
\}, a=1,\ldots,d$, be a basis of $\g$ and $\{ J^a \}$ the dual basis
with respect to a non-degenerate invariant bilinear form on $\g$.

The {\em Gaudin hamiltonians} are linear operators on $V_{\laa}$:
\begin{equation}    \label{the gaudin ham}
\Xi_i = \sum_{j\neq i} \sum_{a=1}^d \frac{J_a^{(i)}
  J^{a(j)}}{z_i-z_j}, \qquad i=1,\ldots,N,
\end{equation}
They commute with the diagonal action of $\g$ on $V_{\laa}$ and hence
their action is well-defined on the subspace of highest weight vectors
in $V_{\laa}$ of an arbitrary dominant integral weight $\mu$ with
respect to the diagonal $\g$--action. We may decompose $V_{\laa}$ with
respect to the diagonal action of $\g$ as
$$
V_{\laa} = \bigoplus_\mu V_\mu \otimes \on{Hom}_\g(V_\mu,V_{\laa}).
$$
Then the space of highest weight vectors of weight $\mu$ is identified
with $\on{Hom}_\g(V_\mu,V_{\laa})$, or, equivalently, with
$$V_{\laa,\la_\infty}^G = (V_{\laa} \otimes V_{\la_\infty})^G,$$
if we write $\mu = -w_0(\la_\infty)$, where $w_0$ is the longest
element of the Weyl group of $\g$.

Consider the problem of simultaneous diagonalization of the Gaudin
hamiltonians in $V_{\laa}$ (or equivalently, in all spaces
$V_{\laa,\la_\infty}^G$). Set
$$
|0\ri = v_{\la_1} \otimes \ldots \otimes v_{\la_N} \in V_{\laa}.
$$ It is an eigenvector of the $\Xi_i$'s. Other eigenvectors are
constructed by a procedure known as the {\em Bethe Ansatz}. We explain
it for $\g=\sw_2$. Let $\{ e,h,f \}$ be the standard basis of $\sw_2$
and set
$$f(w) = \sum_{i=1}^N \frac{f^{(i)}}{w-z_i}.$$
Define the Bethe vector
$$
|w_{1},\ldots,w_{m}\ri = f(w_1) f(w_2) \ldots f(w_m)|0\ri.
$$
It is easy to show that it is an eigenvector of the Gaudin
hamiltonians if and only if the following equations are satisfied:
$$
\sum_{i=1}^N \frac{\la_i}{w_j-z_i} - \sum_{s\neq j} \frac{2}{w_j-w_s}
= 0,
\qquad j=1,\ldots,m.
$$

These are the {\em Bethe Ansatz equations} for $\g=\sw_2$.

One can write analogous systems of equations for a general simple Lie
algebra (or, more generally, a Kac-Moody algebra) $\g$. They are
equations on the set of points $w_1,\ldots,w_m$ colored by simple
roots of $\g$, which we denote by $\al_{i_1},\ldots,\al_{i_m}$:
\begin{equation}    \label{BAE}
\sum_{i=1}^N \frac{\langle \la_i,\chal_{i_j} \rangle}{w_j-z_i} -
\sum_{s \neq j} \frac{\langle \al_{i_s},\chal_{i_j}
\rangle}{w_j-w_s} = 0, \quad j=1,\ldots,m.
\end{equation}

If they are satisfied, then one can construct the corresponding Bethe
vector (see formula \eqref{genbv}) which is an eigenvector of the
Gaudin hamiltonians (see \cite{bf,FFR,RV}). This is a highest weight
vector of weight
$$
\sum_{i=1}^N \la_i - \sum_{j=1}^m \al_{i_j},
$$
so it can only be non-zero if
\begin{equation}    \label{dom wt cond}
\sum_{i=1}^N \la_i - \sum_{j=1}^m \al_{i_j} = \mu,
\end{equation}
where $\mu$ is a dominant integral weight which we write again as $\mu
= -w_0(\la_\infty)$. Then it belongs to $V_{\laa,\la_\infty}^G$,
considered as the subspace of highest weight vectors of weight
$-w_0(\la_\infty)$ in $V_{\laa}$.

For $\g=\sw_2$ it has been proved by I. Scherbak and A. Varchenko
\cite{SV} that for generic $z_i$'s the Bethe vectors form an
eigenbasis in the space of highest weight vectors in $V_{\laa}$ (some
important results in this direction have been obtained earlier by
E. Sklyanin \cite{Skl} using the so-called functional Bethe Ansatz).

\subsection{}

In \cite{FFR}, B. Feigin, N. Reshetikhin and myself have given an
interpretation of the Bethe Ansatz procedure using the spaces of
conformal blocks for representations of the affine Kac-Moody algebra
$\ghat$ associated to $\g$. We showed that the Gaudin hamiltonians
naturally arise from central elements of the completed universal
enveloping algebra of $\ghat$ at the critical level. This center has
been identified by Feigin and myself (see \cite{FF:gd,F:wak}) with the
algebra of functions on the space of $^L G$--{\em opers} on the
punctured disc. Here $^L G$ is the group (of adjoint type) which is
Langlands dual to the (simply-connected) Lie group of $G$. Recall that
passing from $G$ to $^L G$ means switching the sets of weights and
coweights, roots and coroots of $G$ (with respect to a maximal torus),
and at the level of Lie algebras it corresponds to taking the
transpose of the Cartan matrix.

An $^L G$--oper on a smooth curve (or the formal disc) $X$ are triples
$(\F,\nabla,\F_{^L B})$, where $\F$ is a $^L G$--bundle on $X$
equipped with a connection $\nabla$ and a reduction $\F_{^L B}$ to a
Borel subgroup $^L B$ of $^L G$. In more concrete terms, opers may be
described as gauge equivalence classes of first order differential
operators of a certain form. They were defined in this way first by
V. Drinfeld and V. Sokolov in their study \cite{DS} of generalized KdV
hierarchies, and later this definition was made more geometric and
coordinate-independent by A. Beilinson and V. Drinfeld
\cite{BD:opers}.

For example, in the case when $\g = \sw_n$ an oper on a smooth affine
curve (or on the disc) is an equivalence class of operators of the
form
$$
\pa_t + \left( \begin{array}{ccccc}
*&*&*&\cdots&*\\
-1&*&*&\cdots&*\\
0&-1&*&\cdots&*\\
\vdots&\ddots&\ddots&\ddots&\vdots\\ 
0&0&\cdots&-1&*
\end{array} \right),
$$
with respect to the gauge action of the group $N$ of the upper
triangular matrices with $1$'s on the diagonal. It is easy to see that
each gauge class contains a unique operator of the form
$$
\partial_t + \left( \begin{array}{ccccc}
0&v_1&v_2&\cdots&v_{n-1}\\
-1&0&0&\cdots&0\\
0&-1&0&\cdots&0\\
\vdots&\ddots&\ddots&\cdots&\vdots\\
0&0&\cdots&-1&0
\end{array}\right).
$$
But giving such an operator is the same as giving a scalar $n$th
order differential operator
\begin{equation}    \label{first time opers}
L = \partial_t^n+v_1(t) \partial_t^{n-2}+\ldots+v_{n-1}(t)
\end{equation}
(taking into account its transformation properties under changes of
variables, we obtain that it must act from $\Omega^{-(n-1)/2}$ to
$\Omega^{(n+1)/2}$). So the space of $PGL_n$--opers is the space of
operators of the form \eqref{first time opers}, which is incidentally
the phase space of the $n$th KdV hierarchy introduced by Adler and
Gelfand--Dickey when $X$ is the disc.

The interpretation of the Gaudin model in terms of the affine
Kac-Moody algebra of critical level allows us to construct a large
commutative algebra of hamiltonians acting on $V^G_{\laa,\la_\infty}$,
which includes the Gaudin hamiltonians \eqref{the gaudin ham}, and to
view their eigenvalues as $^L G$--opers. The first main result of this
paper (see \thmref{belongs}) is a precise statement as to what kind of
opers may appear as the eigenvalues of the generalized Gaudin
hamiltonians (in the case when $\g=\sw_2$ this was proved in my paper
\cite{F:icmp}).

\bigskip

\noindent{\bf Theorem 1.} {\em There is an injective map from the
spectrum of the generalized Gaudin hamiltonians acting on
$V^G_{\laa,\la_\infty}$ to the set of $^L G$--opers on $\pone$ with
regular singularities at $z_1,\ldots,z_N,\infty$ that have residues
$\la_1,\ldots,\la_N,\la_\infty$ and trivial monodromy representation.}

\bigskip

Thus, eigenvalues of the generalized Gaudin hamiltonians are encoded
by $^L G$--opers on $\pone$ with prescribed singularities and trivial
monodromy. We remark that if we remove the ``no monodromy'' condition,
then we obtain a description of the $^L G$--opers corresponding to the
eigenvalues of the generalized Gaudin hamiltonians acting on the
tensor product of the Verma modules $M_{\la_1} \otimes \ldots \otimes
M_{\la_N}$ (this description in fact holds for arbitrary weights
$\la_1,\ldots,\la_N$).

Though this subject is beyond the scope of this paper, it is worth
noting here that the correspondence between the eigenvectors of the
generalized Gaudin hamiltonians and $^L G$--opers on $\pone$ is an
example of the {\em geometric Langlands correspondence}. This is a
correspondence between $^L G$--local systems on a smooth projective
curve $X$ over $\C$ (possibly, with ramifications at marked points)
and certain sheaves (${\mc D}$--modules) on the moduli spaces of
$G$--bundles on $X$ (possibly, with additional structures at the
marked points). In the case at hand, $X=\pone$ and the local system is
represented by a $^L G$--oper on $\pone$. The corresponding ${\mc
D}$--module on the moduli space of $G$--bundles on $\pone$ with
parabolic structures at $z_1,\ldots,z_N$ and $\infty$ is represented
by the Gaudin system, according to a general construction of Beilinson
and Drinfeld \cite{BD} (see \cite{F:icmp} for more details on this
connection).

\subsection{}

Recall that only special solutions of the Bethe Ansatz equations
\eqref{BAE} may give rise to non-zero eigenvectors of the generalized
Gaudin hamiltonians, namely, the ones which satisfy the condition
\eqref{dom wt cond}. The corresponding eigenvalues are then encoded by
a $^L G$--oper on $\pone$.

Now we want to describe {\em all} solutions of the Bethe Ansatz
equations in geometric terms. It turns out that general solutions are
parameterized by {\em Miura opers}.

While a $^L G$--oper is a triple $(\F,\nabla,\F_{^L B})$, a {\em Miura
$^L G$--oper} is by definition a quadruple $(\F,\nabla,\F_{^L
B},\F'_{^L B})$ where $\F'_{^L B}$ is another $^L B$--reduction of
$\F$, which is preserved by $\nabla$. The space of Miura opers on a
curve $X$ (or on the disc) whose underlying oper has a regular
singularities and trivial monodromy representation (so that $\F$ is
isomorphic to the trivial bundle) is isomorphic to the {\em flag
manifold} $^L G/{}^L B$ of $^L G$. Indeed, in order to define the $^L
B$--reduction $\F'_{^L B}$ of such $\F$ everywhere, it is sufficient
to define it at one point $x \in X$ and then use the connection to
``spread'' it around. But choosing a $^L B$--reduction at one point
means choosing an element of the twist of $^L G/{}^L B$ by $\F_x$, and
so we see that the space of all reductions is isomorphic to the flag
manifold of $^L G$.

The flag manifold is the union of {\em Schubert cells} which are the
$^L B$--orbits. They are parameterized by the Weyl group $W$ of
$\g$. The cell attached to $1 \in W$ is open and dense. It corresponds
to those reductions $\F'_{^L B}$ which are in {\em generic position}
with $\F_{^L B}$. Then there are cells of codimension one labeled by
the simple reflections $s_i$, etc.

Suppose now that we are given a Miura oper corresponding to a $^L
G$--oper on $\pone$ with regular singularities at
$z_1,\ldots,z_N,\infty$. Let $^L H = {}^LB/[{}^L B,{}^L B]$ and $^L
\h$ be its Lie algebra. We construct an $^L H$--bundle on $\pone$
equipped with a connection with regular singularities, i.e., an
operator of the form $\pa_t + {\mb u}(t)$, where ${\mb u}(t)$ is an
$^L \h$--valued function which has poles of order at most one.
Namely, we intersect $\F_{^L B}$ with $\F'_{^L B} w_0$ -- this will be
an $^L H$--bundle, and it inherits a connection from $\F'_{^L
B}$. This map gives us a bijection between Miura $^L G$--opers and $^L
H$--connections. The corresponding map from $^L H$--connections to $^L
G$--opers is called the {\em Miura transformation}.

For example, for $\g=\sw_n$ we have ${\mb u}(t) =
(u_1(t),\ldots,u_n(t))$, and the Miura transformation is given by
the formula
$$
L = (\pa_t+u_1(t)) \ldots (\pa_t+u_n(t)),
$$
where $L$ is the operator \eqref{first time opers}. Hence for
$\g=\sw_2$ we have
$$
\pa^2_t - v(t) =(\pa_t - u(t))(\pa_t + u(t)),
$$
i.e.,
$$
v(t) = u(t)^2 - u'(t)
$$
(note that this is the Poisson map intertwining the KdV and mKdV
hierarchies of soliton equations discovered by R. Miura).

The reductions $\F_{^L B}$ and $\F'_{^L B}$ are going to be in generic
position everywhere on $\pone$ except at finitely many points (see
\lemref{is generic} below). Denote these points by $w_1,\ldots,w_m \in
\pone$. At these points the $^L H$--connection will develop a regular
singularity.  For a generic Miura oper the relative positions at
$w_j$'s will correspond to simple reflections from $W$.  An explicit
computation then shows that the corresponding $^L H$--connection will
have regular singularity with residue $\al_{i_j}$. In addition, our
$^L H$--connection will have regular singularity at $z_i$ with residue
$- \la_i$. So the connection will look like this:
\begin{equation}    \label{cartan}
\pa_t - \sum_{i=1}^N \frac{\la_i}{t-z_i} + \sum_{j=1}^m
\frac{\al_{i_j}}{t-w_j}.
\end{equation}
But the $^L G$--oper underlying our $^L H$--connection has
singularities only at the points $z_1,\ldots,z_N$ and no singularity
at $w_1,\ldots,w_m$. Therefore these singularities must be somehow
erased by the Miura transformation.

We have shown in \cite{FFR} that the oper obtained by applying the
Miura transformation to \eqref{cartan} has no singularity if and only
if the Bethe Ansatz equations are satisfied.  So we obtain an
interpretation of the Bethe Ansatz equations as the conditions that
the singularities of our $^L H$--connection at $w_1,\ldots,w_m$ be
erased by the Miura transformation.

Our connection also has a regular singularity at $\infty$, which is
determined by the relative position of $\F_{^L B}$ and $\F'_{^L B}$ at
$\infty$. If the relative position is $y \in W$, then the residue is
$-y(-w_0(\la_\infty)+\rho)+\rho$. The transformation properties of the
connection determine the residue at $\infty$, so we obtain the
``charge conservation law''
\begin{equation}    \label{special relation first}
\sum_{i=1}^N \la_i - \sum_{j=1}^m \al_{i_j} = y(-w_0(\la_\infty)+\rho)
- \rho.
\end{equation}
This leads to the following statement (implicit already in
\cite{F:icmp}), which is the second main result of this paper (see
\corref{fixed oper}):

\bigskip

\noindent {\bf Theorem 2.} {\em The set of those solutions of the
Bethe Ansatz equations which correspond to a fixed $^L G$--oper is in
bijection with an open and dense subset of the flag manifold $^LG/{}^L
B$.

Further, every solution must satisfy the equation \eqref{special
relation first} for some $y \in W$, and the solutions which satisfy
this equation with fixed $y \in W$ are in bijection with an open
subset of the Schubert cell $^L B w_0 y w_0 {}^L B \subset {}^L G/{}^L
B$.}

\bigskip

In particular, a solution for which we have $y=1$ corresponds to the
one-point Schubert cell in the flag manifold. If this point is
contained in the open dense subset of the flag manifold from Theorem
2, then this solution gives rise to a Bethe eigenvector. It was shown
in \cite{FFR} that the eigenvalues of the Gaudin hamiltonians on this
vector are encoded precisely by the $^L G$--oper obtained by applying
the Miura transformation to the $^L H$--connection
\eqref{cartan}. This follows immediately from the construction of the
Bethe eigenvectors using conformal blocks of {\em Wakimoto modules}
presented in \cite{FFR}.

\subsection{}

Let us summarize the results: the eigenvalues of the hamiltonians of
the Gaudin model associated to a simple Lie algebra $\g$ are encoded
by $^L G$--opers on $\pone$, where $^L G$ is the Langlands dual group
of $G$, which have regular singularities at the marked points and
trivial monodromy. We attach to each solution of the Bethe Ansatz
equations \eqref{BAE} an $^L H$--connection on $\pone$ with regular
singularities. There is a special map from $^L H$--con\-nections to
$^L G$--opers which is called the Miura transformation. The Bethe
Ansatz equations naturally arise as the conditions that the Miura
transformation erases the singularities of the corresponding $^L
H$--connection. The set of all solutions of the Bethe Ansatz equations
\eqref{BAE} is the union of certain open dense subsets of the flag
manifold of the Langlands dual group, one for each oper of the above
type. If the open subset corresponding to an oper $\tau$ contains the
one-point Schubert cell, then the corresponding solution gives rise to
a Bethe eigenvector of the Gaudin hamiltonians whose eigenvalues are
encoded by $\tau$.

One can easily write down the Bethe Ansatz equations for an arbitrary
Kac-Moody algebra $\g$, and it is natural to ask whether the set of
solutions is again the union of open dense subsets of the flag
manifold associated to the Langlands dual group. We show that this is
indeed the case. The subtle point is which flag manifold appears here,
because for infinite-dimensional Kac-Moody algebras there are
non-isomorphic flag manifolds: the ``thick'' flag variety, which is a
proalgebraic variety, and the ``thin'' one, which is an ind-scheme. It
turns out that the relevant flag variety is the thin one, $^L G/{}^L
B_-$. Here $^L G$ is an ind-group corresponding to the Lie algebra $^L
\g$ whose Cartan matrix is the transpose of the Cartan matrix of $\g$,
and $^L B_-$ is the Borel subgroup of $^L G$ that is a proalgebraic
group (see, e.g., \cite{flags}, Ch. VII, for the precise definition).

To establish the connection between solutions of the Bethe Ansatz
equations and points of the ind-flag variety, we proceed in the same
way as in the finite-dimensional case. First, we introduce suitable
notions of opers and Miura opers for Kac-Moody algebras. (Note that in
the case when $\g$ is an untwisted affine Kac-Moody algebra, the
notions of opers and Miura opers have been introduced earlier by
D. Ben-Zvi and myself \cite{FB1}; however, those notions are different
from the ones we introduce here, see \remref{different}.) Using them,
we show that the set of solutions of the Bethe Ansatz equations is an
open and dense subset in the set of Miura opers on the projective line
with prescribed residues at marked points. We then show that, as in
the finite-dimensional case, this set is in bijection with a disjoint
union of the sets of points of certain open dense subsets of the
ind-flag variety $^L G/{}^L B_-$. Thus, we generalize Theorem 2 to the
case of an arbitrary Kac-Moody algebra (see \thmref{final km}).

In the case of a symmetrizable Kac-Moody algebra $\g$ it is also easy
to write down analogues of the Gaudin hamiltonians acting on the
tensor product of integrable representations of $\g$. Then to any
solution of the Bethe Ansatz equations one associates a Bethe
eigenvector of these hamiltonians in the same way as in the
finite-dimensional case. But the connection between the eigenvalues of
the Gaudin hamiltonians on these vectors and opers on the projective
line is not obvious in the infinite-dimensional case. Recall that in
the finite-dimensional case it was based on the concept of conformal
blocks for modules over the {\em affinization of} $\g$ {\em at the
critical level}. It is not immediately clear what should be the
analogue of this Lie algebra when $\g$ is an infinite-dimensional
Kac-Moody algebra. This as a very interesting open problem, which in
fact served as one of our motivations. The fact that our results on
the solutions of the Bethe Ansatz equations apply to general Kac-Moody
algebras indicates that this question may have a good answer. We hope
to return to it in a future publication.

\subsection{}

As a corollary of the above description of the solutions of the Bethe
Ansatz equations we obtain a (rational) action of the group $^L G$ on
the set of solutions of the Bethe Ansatz equations. It is easy to
write down explicitly the action of the one-parameter subgroups
corresponding to the generators $e_i$ of the nilpotent Lie algebra $^L
{\mathfrak n}$. Taking the closure of an orbit of such a subgroup, we
obtain for each simple root a procedure for producing a projective
line (minus finitely many points) worth of new solutions of the Bethe
Ansatz equations from a given one (these projective lines are
precisely the ones appearing in the Bott-Samelson resolutions of the
closures of the Schubert cells).

These procedures were introduced independently and in a different way
by E. Mukhin and A. Varchenko \cite{MV}. They defined what they called
a ``population'' of solutions of the Bethe Ansatz equations as the
closure of the set of all solutions obtained from a given one by
iterating these procedures. In the case when $\g$ is of types $A_n,
B_n$ or $C_n$ they proved (by a method different from ours) that a
population of solutions is isomorphic to the flag manifold of $^L G$
(for $\g=\sw_2$ this had been proved earlier by Scherbak and Varchenko
\cite{SV}, and for $\g=G_2$ this was subsequently proved by
Borisov and Mukhin \cite{BM}). In \secref{comparison} it is shown
that the reproduction procedures of \cite{MV} are equivalent to the
action of one-parameter subgroups of $^L N$ on Miura opers.

\subsection{Plan of the paper.}

The paper is organized as follows. We start in \secref{first} by
defining opers and Miura opers following \cite{BD:opers} and
\cite{F:wak}. We describe opers with regular singularities, explain
the connection between Miura opers and Cartan connections and
establish the correspondence between relative positions of Borel
reductions in a Miura oper and residues of the corresponding Cartan
connection (this last result is borrowed from a forthcoming joint work
with D. Gaitsgory \cite{FG}). Next, we explain in \secref{second} the
connection between solutions of Bethe Ansatz equations and
non-degenerate Miura $G$--opers on $\pone$ with prescribed
singularities at marked points. Here, in order to simplify our
notation, we consider the Bethe Ansatz equations which correspond to
the Gaudin model of $^L \g$. Then the solutions of these equations
correspond to Miura $G$--opers, rather than $^L G$--opers (as
discussed in this Introduction) and hence to points of the flag
manifold of $G$ rather than $^L G$. Our main result is that the set of
solutions of the Bethe Ansatz equations is isomorphic to a union of
open dense subsets of the flag manifold $G/B$ (one for each $G$--oper
on $\pone$ with prescribed singularities).

\secref{third} is devoted to the Gaudin model. We recall how the Bethe
Ansatz equations arise naturally in the problem of diagonalization of
the Gaudin hamiltonians. We explain the construction of the
generalized Gaudin hamiltonians as central elements of the vertex
algebra corresponding to $\ghat$ acting on a suitable space of
coinvariants on $\pone$. We then show that the eigenvalues of the
generalized Gaudin hamiltonians are encoded by $^L G$--opers on
$\pone$ with fixed singularities and trivial monodromy
representation. We discuss an application of this result to the
problem of completeness of Bethe Ansatz. Finally, in \secref{fourth}
we generalize our results on solutions of the Bethe Ansatz equations
to the case of an arbitrary Kac-Moody algebra. We introduce suitable
notions of opers and Miura opers and exhibit the connection between
the solutions of the Bethe Ansatz equations and Miura opers on
$\pone$, much like in the finite-dimensional case.

\subsection{Acknowledgments.} It is a great pleasure to dedicate
this paper, with gratitude and admiration, to Boris Feigin, my
teacher, friend and collaborator of many years. Especially so,
since the results of this paper are based on or motivated by the
results of our previous joint works.

I thank E. Mukhin and A. Varchenko for stimulating discussions, which
encouraged me to revisit my earlier work concerning the Bethe Ansatz
equations and led me to consider these equations when the underlying
Lie algebra is infinite-dimensional.

I also thank D. Gaitsgory and I. Scherbak for valuable discussions and
the referee for useful comments.

\section{Opers and Miura opers}    \label{first}

In this section we first recall the notions of opers and Miura opers
(see \cite{DS,BD,F:wak}). We will then discuss the spaces of opers and
Miura opers with singularities.

\subsection{Opers}    \label{opers}

Let $G$ be a simple algebraic group of adjoint type, $B$ a Borel
subgroup and $N = [B,B]$ its unipotent radical, with the corresponding
Lie algebras $\n \subset \bb\subset \g$. There is an open $B$--orbit
${\bf O}\subset [\n,\n]^\perp/\bb \subset \g/\bb$, consisting of
vectors which are stabilized by the radical $N\subset B$, and such
that all of their negative simple root components, with respect to the
adjoint action of $H = B/N$, are non-zero. This orbit may also be
described as the $B$--orbit of the sum of the projections of simple
root generators $f_i$ of any nilpotent subalgebra $\n_-$, which is in
generic position with $\bb$, onto $\g/\bb$. The torus $H = B/N$ acts
simply transitively on ${\bf O}$, so ${\bf O}$ is an $H$--torsor.

We will often choose a splitting $H \to B$ of the homomorphism $B \to
H$ and the corresponding splitting $\h \to \bb$ at the level of Lie
algebras. Then we will have a Cartan decomposition $\g = \n_- \oplus
\h \oplus \n$. We will choose generators $\{ e_i \}, i=1,\ldots,\ell$,
of $\n$ and generators $\{ f_i \}, i=1,\ldots,\ell$ of $\n_-$
corresponding to simple roots, and denote by $\crho \in \h$ the sum of
the fundamental coweights of $\g$. Then we will have the following
relations: $[\crho,e_i] = 1, [\crho,f_i] = -1$.

Suppose we are given a principal $G$--bundle $\F$ on $X$, which is a
smooth curve, or a disc $D \simeq \on{Spec} \C[[t]]$, or a punctured
disc $D^\times \simeq \on{Spec} \C((t))$, together with a connection
$\nabla$ (automatically flat) and a reduction $\F_B$ to the Borel
subgroup $B$ of $G$. Then we define the relative position of $\nabla$
and $\F_B$ (i.e., the failure of $\nabla$ to preserve $\F_B$) as
follows. Locally, choose any flat connection $\nabla'$ on $\F$
preserving $\F_B$, and take the difference $\nabla - \nabla'$.  It is
easy to show that the resulting local sections of $(\g/\bb)_{\F_B}
\otimes \Omega$, where $\Omega$ is the canonical line bundle of $X$,
are independent of $\nabla'$, and define a global
$(\g/\bb)_{\F_B}$--valued one-form on $X$, denoted by $\nabla/\F_B$.

\medskip

Let $X$ be as above. A $G$--{\em oper} on $X$ is by definition a
triple $(\F,\nabla,\F_B)$, where $\F$ is a principal $G$--bundle $\F$
on $X$, $\nabla$ is a connection on $\F$ and $\F_B$ is a
$B$--reduction of $\F$, such that the one--form $\nabla/\F_B$ takes
values in ${\bf O}_{\F_B} \subset(\g/\bb)_{\F_B}$.

\medskip

This definition is due to A. Beilinson and V. Drinfeld \cite{BD} (in
the case when $X$ is the punctured disc opers were first introduced in
\cite{DS}). Note that ${\bf O}$ is $\C^\times$--invariant, so that
${\bf O} \otimes \Omega$ is a well-defined subset of $(\g/\bb)_{\F_B}
\otimes \Omega$.

Equivalently, the above condition may be reformulated as follows. Let
$U$ be an open subset of a smooth curve $X$ (in the analytic or
Zariski topology) which admits a coordinate $t: U \to {\mathbb A}^1$
(analytic or \'etale, respectively) and a trivialization of $\F_B$,
then with respect to this coordinate and this trivialization the
connection will have the form
\begin{equation}    \label{form of nabla}
\nabla = \pa_t + \sum_{i=1}^\ell \psi_i(t) f_i + {\mb v}(t),
\end{equation}
where each $\psi_i(t)$ is a nowhere vanishing function, and ${\mb
v}(t)$ is a $\bb$--valued function. If we change the trivialization of
$\F_B$, then this operator will get transformed by the corresponding
gauge transformation from the group $B(R)$, where $R$ is the ring of
functions (analytic or algebraic, respectively) on $U$. This
observation allows us to describe opers on $U$ in more concrete terms.

Namely, we obtain from the above description that the space
$\on{Op}_G(U)$ of $G$--opers on $U$ is the quotient of the space of
all operators of the form \eqref{form of nabla}, where $\psi_i(t) \in
R$ is nowhere vanishing, and ${\mb v}(t) \in \bb(R)$, by the action
of the group $B(R)$ by gauge transformations:
$$
g \cdot (\pa_t + A(t)) = \pa_t + g A(t) g^{-1} - \pa_t g \cdot g^{-1}.
$$
The same description applies if $U = D$ or $U = D^\times$, with $t$
being the topological generator of $R$, which is equal to $\C[[t]]$ or
$\C((t))$, respectively.

Since the $B$--orbit ${\bf O}$ is an $H$--torsor, we can use the
$H$--action to make all functions $\psi_i(t)$ equal to $1$ (or any
other non-zero constant). Thus, we obtain that $\on{Op}_G(U)$ is
equal to the quotient of the space $\wt{\on{Op}}_G(U)$ of operators of
the form
\begin{equation}    \label{another form of nabla}
\nabla = \pa_t + \sum_{i=1}^\ell f_i + {\mb v}(t), \qquad {\mb v}(t)
\in \bb(R),
\end{equation}
by the action of the group $N(R)$.

The operator $\on{ad} \crho$ defines the principal gradation on $\bb$,
with respect to which we have a direct sum decomposition $\bb =
\bigoplus_{i\geq 0} \bb_i$. Set
$$
p_{-1} = \sum_{i=1}^\ell f_i.
$$
Let $p_1$ be the unique element of degree 1 in $\n$, such that $\{
p_{-1},2\rv,p_1 \}$ is an $\sw_2$--triple. Let $V_{\can} = \oplus_{i
\in E} V_{\can,i}$ be the space of $\on{ad} p_1$--invariants in
$\n$. Then $p_1$ spans $V_{\on{can},1}$. Choose a linear
generator $p_j$ of $V_{\can,d_j}$ (if the multiplicity of $d_j$ is
greater than one, which happens only in the case $\g=D^{(1)}_{2n},
d_j=2n$, then we choose linearly independent vectors in
$V_{\on{can},d_j}$).

\begin{lem}[\cite{DS}]    \label{free}
The gauge action of $N(R)$ on $\wt{\on{Op}}_G(\on{Spec} R)$ is free,
and each gauge equivalence class contains a unique operator of the
form $\nabla = \pa_t + p_{-1} + {\mathbf v}(t)$, where ${\mathbf v}(t)
\in V_{\can}(R)$, so that we can write
$${\mathbf v}(t) = \sum_{j=1}^\ell v_j(t) \cdot p_j.$$
\end{lem}

\begin{proof}
The operator $\on{ad} p_{-1}$ acts from $\bb_{i+1}$ to $\bb_{i}$
injectively for all $i\geq 0$ and we have $\bb_i = [p_{-1},\bb_{i+1}]
\oplus V_{\on{can},i}$. In particular, $V_0=0$. We claim that each
element of $\pa_t + p_{-1} + {\mb v}(t) \in \wt{\on{Op}}_G(\on{Spec}
R)$ can be uniquely represented in the form
\begin{equation}    \label{gauge}
\pa_t + p_{-1} + {\mb v}(t) = \exp \left( \on{ad} M \right) \cdot
\left( \pa_t + p_{-1} + {\mb c}(t) \right),
\end{equation}
where $M \in \n \otimes R$ and ${\mb c}(t) \in V_{\on{can}} \otimes
R$. To see that, we decompose with respect to the principal gradation:
$M=\sum_{j\geq 0} M_j$, ${\mb v}(t)=\sum_{j\geq 0} {\mb v}_j(t)$,
${\mb c}(t) = \sum_{j\in E} {\mb c}_j(t)$. Equating the homogeneous
components of degree $j$ on both sides of \eqref{gauge}, we obtain
that ${\mb c}_i + [M_{i+1},p_{-1}]$ is expressed in terms of ${\mb
v}_i,{\mb c}_j, j<i$, and $M_j, j\leq i$. The injectivity of $\on{ad}
p_{-1}$ then allows us to determine uniquely ${\mb c}_i$ and
$M_{i+1}$. Hence $M$ and ${\mb c}$ satisfying equation \eqref{gauge}
may be found uniquely by induction, and the lemma follows.
\end{proof}

If we choose another coordinate $s$ such that $t = \varphi(s)$, then
the operator \eqref{another form of nabla} will become
$$
\nabla = \pa_s + \varphi'(s) \sum_{i=1}^\ell f_i + \varphi'(s) \cdot
       {\mb v}(\varphi(s)).
$$
In order to bring it back to the form \eqref{another form of nabla} we
need to apply the gauge transformation by $\crho(\varphi'(s))$, where
we choose a splitting $H \to B$ of the homomorphism $B \to H$ and view
$\crho$ as a homomorphism $\C^\times \to H$. We have
$$
\crho(\varphi'(s)) \cdot \left( \pa_s + \varphi'(s)
\sum_{i=1}^\ell f_i + \varphi'(s) \cdot {\mb v}(\varphi(s))\right)
$$
\begin{equation}    \label{change of var}
= \pa_s + \sum_{i=1}^\ell f_i + \varphi'(s) \crho(\varphi'(s))
  \cdot {\mb v}(\varphi(s)) \cdot \crho(\varphi'(s))^{-1}  -
\crho \cdot \frac{\varphi''(s)}{\varphi'(s)}.
\end{equation}
This formula allows us to glue together opers defined on various open
subsets of a general curve $X$ and thus describe the space
$\on{Op}_G(X)$ in terms of first order differential operators. It also
allows us to describe the space of opers on the disc $D_x = \on{Spec}
\OO_x$, where $\OO_x$ is the completion of the local ring of $U$ at
$x$, or on the punctured disc $D_x^\times = \on{Spec} \K_x$, where
$\K_x$ is the field of fractions of $\OO_x$.

In particular, we obtain the following result. Consider the
$H$--bundle $\Omega^{\crho}$ on $D$. It is uniquely determined by the
following property: for any character $\la: H \to \C^\times$, the line
bundle $\Omega^{\crho} \us{H}\times \la$ associated to the
corresponding one-dimensional representation of $H$ is
$\Omega^{\langle \la,\crho \rangle}$.

\begin{lem}[\cite{F:wak}, Lemma 10.1]    \label{FH}
The $H$--bundle $\F_H = \F_B \underset{B}\times H = \F_B/N$ is
isomorphic to $\Omega^{\crho}$.
\end{lem}

Moreover, it is easy to find transformation formulas for the canonical
representatives of opers. Indeed, by \lemref{free}, there exists a
unique operator $\pa_s + p_{-1} + \ol{{\mathbf v}}(s)$ with
$\ol{{\mathbf v}}(s) \in V_{\can}(R)$ and $g \in B(R)$, such that
\begin{equation}    \label{cano}
\pa_s + p_{-1} + \ol{{\mathbf v}}(s) = g \cdot \left( \pa_s +
\varphi'(s) \sum_{i=1}^\ell f_i + \varphi'(s) \cdot {\mb
v}(\varphi(s)) \right).
\end{equation}
It is straightforward to find that (see \cite{F:wak})
\begin{align}
g &= \exp \left(\frac{1}{2} \frac{\varphi''}{\varphi'} \cdot p_1
\right) \crho(\varphi'), \notag \\ \label{Schwarzian}
\ol{v}_1(s) &= v_1(\varphi(s)) \left( \varphi'
\right)^2 - \frac{1}{2} \{ \varphi,s \}, \\ \ol{v}_j(s) &=
v_j(\varphi(s)) \left( \varphi' \right)^{d_j+1},
\quad \quad j>1, \notag
\end{align}
where $$\{ \varphi,s \} = \frac{\varphi'''}{\varphi'} - \frac{3}{2}
\left( \frac{\varphi''}{\varphi'} \right)^2$$ is the Schwarzian
derivative.

The above formulas describe the transition functions of the
bundle $\F_B$ and hence of $\F$. Namely, they are equal to
$$
\exp \left(\frac{1}{2} \frac{\varphi''}{\varphi'} \cdot p_1
\right) \crho(\varphi'),
$$
where $\varphi(s)$ is the change of coordinate function. Thus, we find
that the bundles $\F$ and $\F_B$ are the same for all opers.

These formulas also imply that under changes of variables, $v_1$
transforms as a projective connection, and $v_j, j>1$, transforms as a
$(d_j+1)$--differential on $U$. Thus, we obtain an isomorphism
\begin{equation}    \label{repr}
\on{Op}_G(X) \simeq  {\mc P}roj(X) \times \bigoplus_{j=2}^\el
\Gamma(X,\Omega^{\otimes(d_j+1)}),
\end{equation}
where ${\mc P}roj(X)$ is the $\Gamma(X,\Omega^{\otimes 2})$--torsor of
projective connections on $X$ (see, e.g., \cite{FB}, Sect. 8.2).

\subsection{Opers for classical Lie groups}    \label{opers class}

For Lie groups of classical types opers may be described in terms of
scalar differential operators. Consider first the case of
$\g=\sw_n$. Then the space of opers on $U = \on{Spec} R$ is the
quotient of the space of operators of the form
\begin{equation}    \label{sln-oper2}
\pa_t + \left( \begin{array}{ccccc}
*&*&*&\cdots&*\\
-1&*&*&\cdots&*\\
0&-1&*&\cdots&*\\
\vdots&\ddots&\ddots&\ddots&\vdots\\ 
0&0&\cdots&-1&*
\end{array} \right),
\end{equation}
where the stars stand for elements of $R$, by the gauge action of the
group $N(R)$ of upper triangular matrices over $R$ with $1$'s on the
diagonal.  It is easy to see that each gauge orbit contains a unique
operator of the form
\begin{equation}    \label{sln-oper1}
\partial_t + \left( \begin{array}{ccccc}
0&v_1&v_2&\cdots&v_{n-1}\\
-1&0&0&\cdots&0\\
0&-1&0&\cdots&0\\
\vdots&\ddots&\ddots&\cdots&\vdots\\
0&0&\cdots&-1&0
\end{array}\right).
\end{equation}

But giving such an operator is the same as giving a scalar $n$th
order differential operator
\begin{equation}    \label{sln-oper}
L = \partial_t^n+v_1(t) \partial_t^{n-2}+\ldots+v_{n-1}(t).
\end{equation}
Thus, we obtain representatives of the gauge equivalence classes that
is different from those described by \lemref{free}. If we look at how
these operators transform under changes of variables, we find that
they transform as operators acting from the $-(n-1)/2$--densities,
i.e., section of the $-(n-1)/2$th power of $\Omega$, to the
$(n+1)/2$--densities. This completely describes the transformation
formulas of the $v_i(t)$'s. Note that if $n$ is even, we need to
choose a square root of $\Omega$, but the resulting space of
differential operators will not depend on this choice. For example, if
$n=2$, we obtain the space of projective connections, i.e., operators
of the form $\pa_t^2 + v(t)$ acting from $\Omega^{-1/2}$ to
$\Omega^{3/2}$. Under changes of coordinates $v(t)$ transforms
according to formula \eqref{Schwarzian}.

For the classical Lie algebra ${\mathfrak sp}_{2n}$ and ${\mathfrak
o}_{2n+1}$ opers may also be realized as scalar differential
operators, as explained by Drinfeld and Sokolov \cite{DS}, Sect. 8
(see also \cite{BD:opers}, Sect. 2). Observe that using the residue
pairing we can identify the dual space of the space of sections of the
line bundle $\Omega^m$ on the punctured disc with that of
$\Omega^{1-m}$. Then the adjoint of a differential operator from
$\Omega^m$ to $\Omega^k$ acts from $\Omega^{1-k}$ to
$\Omega^{1-m}$. Now the space $\on{Op}_{{\mathfrak
sp}_{2n}}(D^\times)$ (resp., $\on{Op}_{{\mathfrak
so}_{2n+1}}(D^\times)$) is realized as the space of self-adjoint
differential operators $L: \Omega^{-n+1/2} \to \Omega^{n+1/2}$ of
order $2n$ (resp., anti-self adjoint operators $L: \Omega^{-n} \to
\Omega^{n+1}$ of order $2n+1$) with the principal symbol $1$.

In the case of $\g={\mathfrak so}_{2n}$ opers may be realized as
scalar pseudo-differential operators (see \cite{DS,BD:opers}).

\subsection{Opers with regular singularities}    \label{reg sing}

Let $x$ be a point of a smooth curve $X$ and $D_x = \on{Spec} \OO_x,
D^\times_x = \on{Spec} \K_x$, where $\OO_x$ is the completion of the
local ring of $x$ and $\K_x$ is the field of fractions of $\OO_x$.
Choose a formal coordinate $t$ at $x$, so that $\OO_x \simeq \C[[t]]$
and $\K_x = \C((t))$. Recall that the space $\on{Op}_G(D_x)$ (resp.,
$\on{Op}_G(D_x^\times)$) of $G$--opers on $D_x$ (resp., $D_x^\times$)
is the quotient of the space of operators of the form \eqref{form of
nabla} where $\psi_i(t)$ and ${\mb v}(t)$ take values in $\OO_x$
(resp., in $\K_x$) by the action of $B(\OO_x)$ (resp., $B(\K_x)$).

A $G$--oper on $D_x$ with regular singularity at $x$ is by definition
(see \cite{BD}, Sect. 3.8.8) a $B(\OO_x)$--conjugacy class of
operators of the form
\begin{equation}    \label{oper with RS1}
\nabla = \pa_t + t^{-1} \left( \sum_{i=1}^\ell \psi_i(t) f_i + {\mb
  v}(t) \right),
\end{equation}
where $\psi_i(t) \in \OO_x, \psi_i(0) \neq 0$, and ${\mb v}(t) \in
\bb(\OO_x)$. Equivalently, it is an $N(\OO_x)$--equivalence class of
operators
\begin{equation}    \label{oper with RS}
\nabla = \pa_t + \frac{1}{t} \left( p_{-1} + {\mb v}(t) \right),
\qquad {\mb v}(t) \in \bb(\OO_x).
\end{equation}
Denote by $\on{Op}_G^{\on{RS}}(D_x)$ the space of opers on $D_x$ with
regular singularity.  It is easy to see (\cite{BD} or \propref{can
real RS}) that the natural map $\on{Op}_G^{\on{RS}}(D_x) \to
\on{Op}_G(D_x^\times)$ is injective. Therefore an oper with regular
singularity may be viewed as an oper on the punctured disc. But to an
oper with regular singularity one can unambiguously attach a point in
$$\g/G := \on{Spec} \C[\g]^G \simeq \C[\h]^W =: \h/W,$$ its residue,
which in our case is equal to $p_{-1} + {\mb v}(0)$.

In particular, the residue of a regular oper $\pa_t + p_{-1} + {\mb
v}(t)$, where ${\mb v}(t) \in \bb(\OO_x)$, is equal to $-\crho$ (see
\cite{BD}). Indeed, a regular oper may be brought to the form
\eqref{oper with RS} by using the gauge transformation with
$\crho(t) \in B(\K_x)$, after which it takes the form
$$
\pa_t + \frac{1}{t} \left( p_{-1} - \crho + t \cdot
\crho(t) ({\mb v}(t)) \crho(t)^{-1} \right).
$$
If ${\mb v}(t)$ is regular, then so is $\crho(t) ({\mb v}(t))
\crho(t)^{-1}$. Therefore the residue of this oper in $\h/W$ is
equal to $-\crho$.

Given $\cla \in \h$, we denote by
$\on{Op}_G^{\on{RS}}(D_x)_{\cla}$ the subvariety of
$\on{Op}_G^{\on{RS}}(D_x)$ which consists of those opers that have
residue $- \cla - \crho \in \h/W$ (in particular,
$\on{Op}_G(D_x) = \on{Op}_G^{\on{RS}}(D_x)_{0}$).

Denote by $\g_{\on{can}}$ the affine subspace of $\g$ consisting of
all elements of the form $$p_{-1} + \sum_{j \in E} y_j p_j.$$ Recall
from \cite{Ko} that the adjoint orbit of any regular element in the
Lie algebra $\g$ contains a unique element that belongs to
$\g_{\on{can}}$, and the corresponding morphism $\g_{\on{can}} \to
\h/W$ is an isomorphism.

\begin{prop}[\cite{BD}, Prop. 3.8.9]    \label{can real RS}
The canonical representatives of opers with regular singularities have
the form
\begin{equation}    \label{can form2}
\pa_t + p_{-1} + \sum_{j \in E} t^{-j-1} c_j(t) p_j,
\qquad c_j(t) \in \C[[t]].
\end{equation}
Moreover, the residue of this oper is realized in $\g_{\on{can}}$ as
\begin{equation}    \label{res2}
p_{-1} +  \left( c_1(0) + \frac{1}{4} \right) p_1 + \sum_{j \in E,
  j>1} c_j(0) p_j.
\end{equation}
\end{prop}

Let $(\F,\nabla,\F_B) \in \on{Op}_G^{\on{RS}}(D_x)$. For each
finite-dimensional representation $V$ of $G$, consider the system of
differential equations with regular singularities $\nabla \cdot
\phi_V(t) = 0$, where $\phi_V(t)$ takes values in $V$. For varying $V$
the solutions of these equations give rise to a well-defined solution
with values in $G$, whose monodromy around $x$ is a well-defined
conjugacy class in $G$.

Now let $\cla$ be a dominant integral coweight of $\g$. Following
Drinfeld, introduce the variety $\on{Op}_G(D_x)_{\cla}$ as the
quotient of the space of operators of the form
\begin{equation}    \label{psi la}
\nabla = \pa_t + \sum_{i=1}^\ell \psi_i(t) f_i + {\mb v}(t),
\end{equation}
where $$\psi_i(t) = t^{\langle \al_i,\cla \rangle}(\kappa_i +
t(\ldots)) \in \OO_x, \qquad \kappa_i \neq 0$$ and ${\mb v}(t) \in
\bb(\OO_x)$, by the gauge action of $B(\OO_x)$. Equivalently,
$\on{Op}_G(D_x)_{\cla}$ is the quotient of the space of operators of
the form
\begin{equation}    \label{psi la1}
\nabla = \pa_t + \sum_{i=1}^\ell t^{\langle \al_i,\cla \rangle} f_i +
{\mb v}(t),
\end{equation}
where ${\mb v}(t) \in \bb(\OO_x)$, by the gauge action of
$N(\OO_x)$. Considering the $N(\K_x)$--class of such an operator, we
obtain an oper on $D_x^\times$. Thus, we have a map
$\on{Op}_G(D_x)_{\cla} \to \on{Op}_G(D_x^\times)$.

\begin{lem}    \label{no mon}
The map $\on{Op}_G(D_x)_{\cla} \to \on{Op}_G(D_x^\times)$ is injective
and its image is contained in the subvariety
$\on{Op}_G^{\on{RS}}(D_x)_{\cla}$. Moreover, the points of
$\on{Op}_G(D_x)_{\cla}$ are precisely those $G$--opers with regular
singularity and residue $\cla$ which have no monodromy around $x$.
\end{lem}

\begin{proof}
By using the gauge transformation with $(\cla+\crho)(t)$, we
bring the operator \eqref{psi la} to the form \eqref{oper with RS1},
with
\begin{multline}    \label{after conj}
{\mb v}(t) = -(\cla+\crho) + {\mb v}_0(t) + \sum_{\al \in \De_+}
{\mb v}_\al(t), \\ {\mb v}_0 \in \h \otimes t\C[[t]], \quad
{\mb v}_\al(t) \in \n_\al \otimes t^{\langle \al,\cla+\crho \rangle}
\C[[t]],
\end{multline}
and the $N[[t]]$--equivalence class of \eqref{psi la} is mapped to the
$(\cla+\crho)(t) N[[t]] (\cla+\crho)(t)^{-1}$ class of the conjugate
operator. It is then easy to see that the subgroup of $N[[t]]$ which
preserves the operators with ${\mb v}(t)$ of the form \eqref{after
conj} is precisely $(\cla+\crho)(t) N[[t]] (\cla+\crho)(t)^{-1}$. This
proves the first statement.

To prove the second statement, observe that the monodromy of $\nabla$
is trivial if and only if $\nabla$ is gauge equivalent, under the
gauge action of the entire loop group $G((t))$, to a regular
connection (not necessarily an oper). Therefore the second statement
is equivalent to the statement that an oper $\tau \in
\on{Op}_G^{\on{RS}}(D_x)_{\cla}$ is gauge equivalent to a regular
connection if and only if it belongs to $\on{Op}_G(D_x)_{\cla}$. But
$G((t)) = G[[t]] B((t))$, and the gauge action of $G[[t]]$ preserves
the space of regular connections. Therefore if an oper is gauge
equivalent to a regular connection, then its $B((t))$ gauge class
already must contain a regular connection. The oper condition then
implies that this gauge class contains a connection operator of the
form \eqref{psi la}, where $\psi_i(t) = t^{\langle \al_i,\cmu
\rangle}(\kappa_i + t(\ldots)) \in \OO_x, \kappa_i \neq 0$ and ${\mb
v}(t) \in \bb(\OO_x)$ for some integral dominant coweight $\cmu$ of
$\g$. But according to the above calculation, the residue of such an
oper is equal to $-\cmu-\crho$. This gives us the second statement of
the lemma.
\end{proof}

\subsection{Miura opers}    \label{miura opers}

By definition (see \cite{F:wak}, Sect. 10.3), a {\em Miura $G$--oper}
on $X$ (which is a smooth curve or a disc) is a quadruple
$(\F,\nabla,\F_B,\F'_B)$, where $(\F,\nabla,\F_B)$ is a $G$--oper on
$X$ and $\F'_B$ is another $B$--reduction of $\F$ which is preserved
by $\nabla$.

We denote the space of Miura $G$--opers on $X$ by $\on{MOp}_G(X)$.

\medskip

A $B$--reduction of $\F$ which is preserved by the connection $\nabla$
is uniquely determined by a $B$--reduction of the fiber $\F_x$ of $\F$
at any point $x \in X$ (in the case when $U=D$, $x$ has to be the origin
$0 \in D$). The set of such reductions is the $\F_x$--twist
\begin{equation}    \label{twist of flag}
(G/B)_{\F_x} = \F_x \us{G}\times G/B = \F'_{B,x} \us{B}\times G/B =
  (G/B)_{\F'_{B,x}}
\end{equation}
of the flag manifold $G/B$. If $X$ is a curve or a disc and the oper
connection has a regular singularity and trivial monodromy
representation, then this connection gives us a global (algebraic)
trivialization of the bundle $\F$. Then any $B$--reduction of the
fiber $\F_x$ gives rise to a global (algebraic) $B$--reduction of
$\F$. Thus, we obtain:

\begin{lem}    \label{isom with flags}
Suppose that we are given an oper $\tau$ on a curve $X$ (or on the
disc) such that the oper connection has a regular singularity and
trivial monodromy. Then for each $x \in X$ there is a canonical
isomorphism between the space of Miura opers with the underlying oper
$\tau$ and the twist $(G/B)_{\F'_{B,x}}$.
\end{lem}

Recall that the $B$--orbits in $G/B$, known as the Schubert cells, are
parameterized by the Weyl group $W$ of $G$. Let $w_0$ be the longest
element of the Weyl group of $G$. Denote the orbit $B w_0 w B \subset
G/B$ by $S_w$ (so that $S_1$ is the open orbit). We obtain from the
second description of $(G/B)_{\F'_x}$ given in formula \eqref{twist of
flag} that $(G/B)_{\F'_x}$ decomposes into a union of locally closed
subvarieties $S_{w,\F'_{B,x}}$, which are the $\F'_{B,x}$--twists of
the Schubert cells $S_w$. The $B$--reduction $\F_{B,x}$ defines a
point in $(G/B)_{\F'_{B,x}}$. We will say that the $B$--reductions
$\F_{B,x}$ and $\F'_{B,x}$ are in {\em relative position} $w$ if
$\F_{B,x}$ belongs to $S_{w,\F'_{B,x}}$. In particular, if it belongs
to the open orbit $S_{1,\F'_{B,x}}$, we will say that $\F_{B,x}$ and
$\F'_{B,x}$ are in generic position.

A Miura $G$--oper is called {\em generic} at the point $x \in X$ if
the $B$--reductions $\F_{B,x}$ and $\F'_{B,x}$ of $\F_x$ are in
generic position. In other words, $\F_{B,x}$ belongs to the stratum
$\on{Op}_G(X) \times S_{1,\F'_{B,x}} \subset \on{MOp}_G(X)$. Being
generic is an open condition. Therefore if a Miura oper is generic at
$x \in X$, then there exists an open neighborhood $U$ of $x$ such that
it is also generic at all other points of $U$. We denote the space of
generic Miura opers on $U$ by $\on{MOp}_G(U)_{\on{gen}}$.

\begin{lem}    \label{is generic}
Suppose we are given a Miura oper on the disc $D_x$ around a point $x
\in X$. Then its restriction to the punctured disc $D_x^\times$ is
generic.
\end{lem}

\begin{proof}
Since being generic is an open condition, we obtain that if a Miura
oper is generic at $x$, it is also generic on the entire $D_x$. Hence
we only need to consider the situation where the Miura oper is not
generic at $x$, i.e., the two reductions $\F_{B,x}$ and $\F'_{B,x}$
are in relative position $w \neq 1$. Let us trivialize the $B$--bundle
$\F_B$, and hence the $G$--bundle $\F_G$ over $D_x$. Then $\nabla$
gives us a connection on the trivial $G$--bundle which we can bring to
the canonical form
$$
\nabla = \pa_t + p_{-1} + \sum_{j=1}^\ell v_j(t) \cdot p_j
$$
(see \lemref{free}). It induces a connection on the trivial
$G/B$--bundle. We are given a point $gB$ in the fiber of the latter
bundle which lies in the orbit $S_w = B w_0 w B$, where $w \neq
1$. Consider the horizontal section whose value at $x$ is $gB$, viewed
as a map $D_x \to G/B$. We need to show that the image of this map
lies in the open $B$--orbit $S_1 = B w_0 B$ over $D_x^\times$, i.e.,
it does not lie in the orbit $S_y$ for any $y \neq 1$.

Suppose that this is not so, and the image of the horizontal section
actually lies in the orbit $S_y$ for some $y \neq 1$. Since all
$B$--orbits are $H$--invariant, we obtain that the same would be true
for the horizontal section with respect to the connection $\nabla' = h
\nabla h^{-1}$ for any constant element of $H$. Choosing $h =
\crho(a)$ for $a \in \C^\times$, we can bring the connection to the
form
$$
\pa_t + a^{-1} p_{-1} + \sum_{j=1}^\ell a^{d_j} v_j(t) \cdot p_j.
$$
Changing the variable $t$ to $s = a^{-1} t$, we obtain the connection
$$
\pa_s + p_{-1} + \sum_{j=1}^\ell a^{d_j+1} v_j(t),
$$
so choosing small $a$ we can make the functions $v_j(t)$ arbitrarily
small. Therefore without loss of generality we can consider the case
when our connection operator is $\nabla = \pa_t + p_{-1}$.

In this case our assumption that the horizontal section lies in $S_y,
y \neq 1$, means that the vector field $\xi_{p_{-1}}$ corresponding to
the infinitesimal action of $p_{-1}$ on $G/B$ is tangent to an orbit
$S_y, y \neq 1$, in the neighborhood of some point $gB$ of $S_w
\subset G/B, w \neq 1$. But then, again because of the $H$--invariance
of the $B$--orbits, the vector field $\xi_{h p_{-1} h^{-1}}$ is also
tangent to this orbit for any $h \in H$. For any $i=1\ldots,\ell$,
there exists a one-parameter subgroup $h_{\ep}^{(i)}, \ep \in
\C^\times$ in $H$, such that $\underset{\ep \to 0}\lim \; \ep p_{-1}
\ep^{-1} = f_i$. Hence we obtain that each of the vector fields
$\xi_{f_i}, i=1\ldots,\ell$, is tangent to the orbit $S_y, y \neq 1$,
in the neighborhood of $gB \in S_w, w \neq 1$. But then all
commutators of these vectors fields are also tangent to this
orbit. Hence we obtain that all vector fields of the form $\xi_p, p
\in \n_-$, are tangent to $S_y$ in the neighborhood of $gB \in
S_w$.

Consider any point of $G/B$ that does not belong to the open dense
orbit $S_1$. Then the quotient of the tangent space to this point by
the tangent space to the $B$--orbit passing through this point is
non-zero and the vector fields from the Lie algebra $\n_-$ map
surjectively onto this quotient. Therefore they cannot be tangent to
the orbit $S_y, y \neq 1$, in a neighborhood of $gB$. Therefore our
Miura oper is generic on $D_x^\times$.
\end{proof}

This lemma shows that any Miura oper on any smooth curve $X$ is
generic over an open dense subset.

Consider the $H$--bundles $\F_H = \F_B/N$ and $\F'_H = \F'_B/N$
corresponding to a generic Miura oper $(\F,\nabla,\F_B,\F'_B)$ on
$X$. If ${\mc P}$ is an $H$--bundle, then applying to it the
automorphism $w_0$ of $H$, we obtain a new $H$--bundle which we denote
by $w_0^*(\F_H)$.

\begin{lem}[\cite{F:wak},Lemma 10.3]    \label{H bundles isom}
For a generic Miura oper $(\F,\nabla,\F_B,\F'_B)$ the $H$--bundle
$\F'_H$ is isomorphic to $w_0^*(\F_H)$.
\end{lem}

\begin{proof}
Consider the vector bundles $\g_\F = \F \underset{G}\times \g$,
$\bb_{\F_B} = \F_B \underset{B}\times \bb$ and $\bb_{\F'_B} = \F'_B
\underset{G}\times \bb$. We have the inclusions $\bb_{\F_B},
\bb_{\F'_B} \subset \g_\F$ which are in generic position. Therefore
the intersection $b_{\F_B} \cap b_{\F'_B}$ is isomorphic to
$\bb_{\F_B}/[\bb_{\F_B},\bb_{\F_B}]$, which is the trivial vector
bundle with the fiber $\h$. It naturally acts on the bundle $\g_\F$
and under this action $\g_\F$ decomposes into a direct sum of $\h$ and
the line subbundles $\g_{F,\al}, \al \in \De$. Furthermore,
$\bb_{\F_B} = \bigoplus_{\al \in \De_+} \g_{F,\al}, \bb_{\F'_B} =
\bigoplus_{\al \in \De_+} \g_{F,w_0(\al)}$. Since the action of $B$ on
$\n/[\n,\n]$ factors through $H = B/N$, we find that
$$
\F_H \underset{H}\times \bigoplus_{i=1}^\ell \C_{\al_i} \simeq
\bigoplus_{i=1}^\ell \g_{\F,\al_i}, \qquad \F'_H \underset{H}\times
\bigoplus_{i=1}^\ell \C_{\al_i} \simeq \bigoplus_{i=1}^\ell
\g_{\F,w_0(\al_i)}.
$$
Therefore we obtain that
$$
\F_H \underset{H}\times \C_{\al_i} \simeq \F'_H
\underset{H}\times \C_{w_0(\al_i)}, \qquad i=1,\ldots,\ell.
$$
Since $G$ is of adjoint type by our assumption, the above associated
line bundles completely determine $\F_H$ and $\F'_H$, and the above
isomorphisms imply that $\F'_H \simeq w_0^*(\F_H)$.
\end{proof}

Since the $B$--bundle $\F'_B$ is preserved by the oper connection
$\nabla$, we obtain a connection $\ol{\nabla}$ on $\F'_H$ and hence on
$\F_H \simeq \Omega^{\crho}$. Therefore we obtain a morphism ${\mb
a}$ from the variety $\on{MOp}_G(U)_{\on{gen}}$ of generic Miura opers
on $U$ to the variety of connections $\on{Conn}_U$ on the $H$--bundle
$\Omega^{\crho}$ on $U$.

Explicitly, connections on $\Omega^{\crho}$ may be described as
follows. If we choose a local coordinate $t$ on $U$, then we
trivialize $\Omega^{\crho}$ and represent the connection as an
operator $\pa_t + {\mb u}(t)$, where ${\mb u}(t)$ is an $\h$--valued
function on $U$. If $s$ is another coordinate such that
$t=\varphi(s)$, then this connection will be represented by the
operator
\begin{equation}    \label{trans for conn}
\pa_s + \varphi'(s) {\mb u}(\varphi(s)) - \crho \cdot
\frac{\varphi''(s)}{\varphi'(s)}.
\end{equation}

\begin{prop}[\cite{F:wak},Prop. 10.4]    \label{map beta}
The morphism ${\mb a}: \on{MOp}_G(U)_{\on{gen}} \to
\on{Conn}_U$ is an isomorphism of algebraic varieties.
\end{prop}

\begin{proof}
We define a morphism ${\mb b}$ in the opposite direction. Suppose we
are given a connection $\ol\nabla$ on the $H$--bundle
$\Omega^{\crho}$ on $D$. We associate to it a generic Miura oper
as follows. Let us choose a splitting $H \to B$ of the homomorphism $B
\to H$ and set $\F = \Omega^{\crho} \underset{H}\times G, \F_B =
\Omega^{\crho} \underset{H}\times B$, where we consider the
adjoint action of $H$ on $G$ and on $B$ obtained through the above
splitting. The choice of the splitting also gives us the opposite
Borel subgroup $B_-$, which is the unique Borel subgroup in generic
position with $B$ containing $H$. Let again $w_0$ be the longest
element of the Weyl group of $\g$. Then $w_0 B$ is a $B$--torsor
equipped with a left action of $H$, so we define the $B$--subbundle
$\F'_B$ of $\F$ as $\Omega^{\crho} \underset{H}\times w_0 B$.

Observe that the space of connections on $\F$ is isomorphic to the
direct product
$$
\on{Conn}_U \times \bigoplus_{\al \in \De}
\Gamma(U,\Omega^{\al(\crho) + 1}).
$$ Its subspace corresponding to negative simple roots is isomorphic
to $\left( \bigoplus_{i=1}^\ell \g_{-\al_i} \right) \otimes R$. Having
chosen a basis element $f_i$ of $\g_{-\al_i}$ for each
$i=1,\ldots,\ell$, we now construct an element $p_{-1} =
\sum_{i=1}^\ell f_i \otimes 1$ of this space. Now we set $\nabla =
\ol\nabla + p_{-1}$. By construction, $\nabla$ has the correct
relative position with the $B$--reduction $\F_B$ and preserves the
$B$--reduction $\F'_B$. Therefore the quadruple
$(\F,\nabla,\F_B,\F'_B)$ is a generic Miura oper on $U$. We define the
morphism ${\mb b}$ by setting ${\mb b}(\ol{\nabla}) =
(\F,\nabla,\F_B,\F'_B)$.

This map is independent of the choice of a splitting $H \to B$ and of
the generators $f_i, i=1,\ldots,\ell$. Indeed, changing the splitting
$H \to B$ amounts to conjugating of the old splitting by an element
of $N$. This is equivalent to applying to $\nabla$ the gauge
transformation by this element. Therefore it will not change the
underlying Miura oper structure. Likewise, rescaling of the generators
$f_i$ may be achieved by a gauge transformation by a constant element
of $H$, and this again does not change the Miura oper structure. Thus,
the morphism ${\mb b}$ is well-defined. It is clear from the
construction that ${\mb a}$ and ${\mb b}$ are mutually inverse
isomorphisms.
\end{proof}

More generally, we define, for any dominant integral coweight $\cla
\in \h$, {\em Miura $G$--opers of coweight} $\cla$ on $D_x$ as
quadruples $(\F,\nabla,\F_B,\F'_B)$, where $(\F,\nabla,\F_B)$ is a
$G$--oper on $D_x^\times$ with regular singularity which belongs to
$\on{Op}_G(D_x)_{\cla}$ and $\F'_B$ is another $B$--reduction of $\F$
which is preserved by $\nabla$.

We denote the space of Miura $G$--opers of coweight $\cla$ on $D_x$ by
$\on{MOp}_G(D_x)_{\cla}$. In particular, if $\cla=0$ we obtain the old
definition of Miura opers on $D_x$. It is clear that we have an
isomorphism
$$
\on{MOp}_G(D_x)_{\cla} \simeq \on{Op}_G(D_x)_{\cla} \times
(G/B)_{\F'_{B,x}}.
$$
We define the relative positions of $\F_B$ and $\F'_B$ in the same way
as for $\cla=0$ and denote by $\on{MOp}_G(D_x)_{\cla,\on{gen}}$ the
variety of generic Miura opers of coweight $\cla$.

Let $\on{Conn}^{\on{RS}}_{D_x,\cla}$ be the variety of connections on
the $H$--bundle $\Omega^{\crho}$ over $D_x$ with regular singularity
at $x$ and residue $-\cla$. With respect to a coordinate $t$ at $x$,
the corresponding connection operator has the form $$\ol\nabla = \pa_t
+ \cla t^{-1} + {\mb u}(t), \qquad {\mb u}(t) \in \h[[t]].$$ Denote by
$\Omega^{\crho}(- \cla \cdot x)$ the $H$--bundle on $D_x$ determined
by the associated line bundles
$$
\Omega^{\crho}(- \cla \cdot x) \underset{H}\times \C_{\al_i} =
\Omega(-\langle \al_i,\cla \rangle x).
$$ We have a morphism $${\mb b}_{\cla}: \on{Conn}^{\on{RS}}_{D_x,\cla}
\to \on{MOp}_G(D_x)_{\cla,\on{gen}}$$ sending such a connection
$\ol\nabla$ to the triple $(\F,\nabla,\F_B,\F'_B)$, where $$\F =
\Omega^{\crho} \underset{H}\times G, \qquad \F_B = \Omega^{\crho}
\underset{H}\times B, \qquad \F_B = \Omega^{\crho} \underset{H}\times
w_0 B,$$ and $\nabla = \ol{\nabla} + p_{-1}$, or equivalently,
$$\F = \Omega^{\crho}(- \cla \cdot x) \underset{H}\times G, \qquad
\F_B = \Omega^{\crho}(- \cla \cdot x) \underset{H}\times B, \qquad
\F_B = \Omega^{\crho}(- \cla \cdot x) \underset{H}\times w_0 B,$$
$$
\nabla = \cla(t)^{-1}(\ol{\nabla} + p_{-1})\cla(t) = \pa_t +
\sum_{i=1}^\ell t^{\langle \al_i,\cla \rangle} f_i + {\mb u}(t)
$$
(the corresponding oper does not depend on the choice of $t$). We
prove in the same way as in \propref{map beta} that ${\mb b}_{\cla}$
is an isomorphism.

\subsection{Miura transformation}    \label{miura trans}

Under the isomorphism of \propref{map beta}, the natural forgetful
morphism $\on{MOp}_G(U)_{\on{gen}} \to \on{Op}_G(U)$ becomes a map
$\on{Conn}_U \to \on{Op}_G(U)$. We call this map the {\em Miura
transformation}. The origin of this terminology is as follows.  In
\secref{opers class} we described a realization of opers for Lie
algebras of classical types in terms of scalar differential
operators. These realizations may be used to describe explicitly the
Miura transformation as well.

In the case of $\sw_n$ the space $\on{Op}_{\sw_n}(D^\times)$ consists
of differential operators of the form \eqref{sln-oper}. The space
$\on{Conn}_{D^\times}$ consists of the operators $\pa_t + {\mb u}(t)$,
where ${\mb u}(t) \in \h((t))$ may be viewed as an $n$--tuple
$(u_1(t),\ldots,u_n(t))$ such that $\sum_{i=1}^n u_i(t) = 0$. The
Miura transformation sends $\pa_t + {\mb u}(t)$ to the operator
\begin{equation}    \label{miura trans for sln}
L = (\pa_t+u_1(t)) \ldots (\pa_t+u_n(t)).
\end{equation}
In particular, for $\g=\sw_2$ we obtain a map sending a connection
$\pa_t + u(t)$ to the projective connection $\pa^2_t - v(t)$ where
$$
\pa^2_t - v(t) = (\pa_t - u(t))(\pa_t + u(t)),
$$
i.e.,
$$
u(t) \mapsto v(t) = u(t)^2 - u'(t).
$$
This map was first introduced by R. Miura as the Poisson map from the
phase space of the mKdV hierarchy (the space $\on{Conn}_{D^\times}$ in
our notation) to the phase space of the KdV hierarchy (the space
$\on{Op}_G(D^\times)$ in our notation). This is the reason why we call
this map (for an arbitrary $\g$) the Miura transformation.

One can also write down explicit formulas for the Miura transformation
for other simple Lie algebras of classical types. As we have seen in
\secref{opers class} (following \cite{DS}), in the case of the Lie
algebras ${\mathfrak sp}_{2n}$ and ${\mathfrak so}_{2n+1}$ the spaces
of opers consist of self-adjoint differential operators $L:
\Omega^{-n+1/2} \to \Omega^{n+1/2}$ of order $2n$ (resp., anti-self
adjoint operators $L: \Omega^{-n} \to \Omega^{n+1}$ of order $2n+1$)
with principal symbol $1$. Identifying the Cartan subalgebras of these
Lie algebras with $\C^n$, we obtain an identification of the
corresponding space of connections with the space of $n$--tuples
$(u_1(t),\ldots,u_n(t))$, where $u_i(t) \in \C((t))$. Then the Miura
transformation takes the form
$$
L = (\pa_t+u_1(t)) \ldots (\pa_t+u_n(t)) (\pa_t-u_n(t)) \ldots
(\pa_t-u_1(t))
$$
for $\g={\mathfrak sp}_{2n}$ and
$$
L = (\pa_t+u_1(t)) \ldots (\pa_t+u_n(t)) \pa_t (\pa_t-u_n(t)) \ldots
(\pa_t-u_1(t))
$$
for $\g={\mathfrak so}_{2n+1}$.

Finally, in the case of $\g={\mathfrak so}_{2n}$ the Miura
transformation is realized by the formula
$$
L = (\pa_t+u_1(t)) \ldots (\pa_t+u_n(t)) \pa_t^{-1} (\pa_t-u_n(t))
\ldots (\pa_t-u_1(t))
$$
(see \cite{DS}, Sect. 8).

\subsection{Singularities of Miura opers}

Now suppose that we are given a Miura oper of coweight $\cla$ on the
disc $D_x$ such that the reduction $\F'_{B,x}$ has relative position
$w$ with $\F_{B,x}$ at $x$. The restriction of this Miura oper to the
punctured disc $D_x^\times$ is generic by \lemref{is generic} (which
is easily generalized to the case of an arbitrary $\cla$), and hence
it corresponds, by \propref{map beta}, to a connection $\ol\nabla$ on
the $H$--bundle $\Omega^{\crho}$ over $D_x^\times$. We would like to
describe the singularity of this connection at $x$.

First of all, we claim that $\ol{\nabla}$ has a regular singularity at
$x$. Indeed, the corresponding $G$--oper may be represented by the
connection operator $\nabla = \ol{\nabla} + p_{-1}$. Let us choose a
coordinate $t$ at $x$ and write $\ol\nabla = \pa_t + {\mathbf u}(t)$,
where ${\mb u}(t) \in \h((t))$. Then
\begin{equation}    \label{Miura conn}
\nabla = \pa_t + p_{-1} + {\mb u}(t), \qquad {\mb u}(t) \in \h((t)).
\end{equation}
The corresponding $G$--oper should be regular, i.e., there should
exist an element $g \in N((t))$ such that $g \nabla g^{-1}$ has no
singularity at $t=0$, so that the equation $g \nabla g^{-1} \cdot
\phi(t) = 0$ has solutions in $G[[t]]$ for arbitrary initial
conditions in $G$. But then the equation $\nabla \phi(t) = 0$ would
have solutions in $G((t))$ (for arbitrary initial conditions in
$G$). This implies that $\ol{\nabla}$ has at most regular singularity.
Suppose that this is not so. Then there would exist a dominant
integral weight $\chi$ such that $\langle \chi,{\mb u}(t) \rangle$ has
a pole of order higher than $1$. But then consider the equation
$\ol{\nabla} \phi(t) = 0$, where $\phi(t)$ takes values in
$V_{-w_0(\chi)}$. Clearly, the component of the solution lying in the
subspace of lowest weight $-\chi$ would not belong to $\C((t))$, which
is a contradiction.

Thus, $${\mb u}(t) = \cmu t^{-1} + \on{reg}.$$ for some integral
coweight $\cmu$. Using the gauge transformation with $\crho(t) \in
B((t))$, we obtain that the operator $\nabla$ given by \eqref{Miura
conn} is gauge equivalent to the operator
$$
\pa_t + \frac{1}{t} \left( p_{-1} - \crho + \cmu + t(\ldots)
\right).
$$
Therefore the $G$--oper corresponding to $\nabla$ is an oper with
regular singularity (see \secref{reg sing}), whose residue in $\h/W$
is equal to the image of $- \crho + \cmu \in \h$. But by our
assumption this oper belongs to $\on{Op}_G(D_x)_{\cla}$, hence its
residue is the image of $-\cla-\crho$ in $\h/W$. Therefore we obtain
that there exists $y \in W$ such that $- \crho + \cmu = -
y(\cla+\crho)$, i.e., $\cmu = \crho - y(\cla+\crho)$.

We wish to show that $y=w$, where $w$ is the relative position of
$\F'_{B,x}$ with $\F_{B,x}$. Let us make a more precise statement.

Denote by $\on{Conn}^{\on{RS}}_{D_x,\cla,w}$ the variety of all
connections on the $H$--bundle $\Omega^{\crho}$ with regular
singularity at $x$ and residue $-w(\cla+\crho)+\crho$. We have a
morphism
$$
{\mb b}^{\on{RS}}_{\cla,w}: \on{Conn}^{\on{RS}}_{D_x,\cla,w} \to
\on{Op}_G^{\on{RS}}(D_x)
$$
defined as in \secref{miura opers}. Namely, we send a connection
$\ol\nabla \in \on{Conn}^{\on{RS}}_{D_x,\cla,w}$ to the oper
$(\F,\nabla,\F_B)$ where we set $\F = \Omega^{\crho}
\underset{H}\times G, \F_B = \Omega^{\crho} \underset{H}\times B$ and
$\nabla = \ol{\nabla} + p_{-1}$.

Explicitly, after choosing a coordinate $t$ on $D$, we can
write $\ol\nabla$ as $\pa_t + t^{-1} {\mb u}(t)$, where ${\mb u}(t)
\in \h[[t]]$. Then the corresponding oper with regular singularity is
the $N((t))$--equivalence class of the operator
$$
\nabla = \pa_t + p_{-1} + t^{-1} {\mb u}(t),
$$
which is the same as the $N[[t]]$--equivalence class of the operator
$$
\crho(t) \nabla \crho(t)^{-1} = \pa_t + t^{-1} (p_{-1} - \crho + {\mb
  u}(t))
$$
(so it is indeed an oper with regular singularity).

Denote by $\on{Conn}^{\on{reg}}_{D_x,w,\cla}$ the reduced part of the
preimage of $\on{Op}_G(D_x)_{\cla} \subset
\on{Op}_G^{\on{RS}}(D_x)$ under this morphism. Then we have a morphism
$$
{\mb b}_{\cla,w}: \on{Conn}^{\on{reg}}_{D_x,\cla,w} \to
\on{MOp}_G(D_x)_{\cla},
$$
which sends $\ol\nabla \in \on{Conn}^{\on{reg}}_{D_x,\cla,w}$ to the
Miura oper $(\F,\nabla,\F_B,\F'_B)$, where $(\F,\nabla,\F_B)$ are as
above and $\F'_B = \Omega^{\crho} \underset{H}\times w_0 B$,
where $w_0$ is the longest element of the Weyl group.

Recall that in \secref{miura opers} we have established an isomorphism
between the variety $\on{MOp}_G(D_x)_{\cla}$ of Miura opers of
coweight $\cla$ on $D_x$ and the product $\on{Op}_G(D_x)_{\cla} \times
(G/B)_{\F'_{B,x}}$. Denote by $\on{MOp}_G(D_x)_{\cla,w} \subset
\on{MOp}_G(D_x)_{\cla}$ the subvariety of those Miura opers of coweight
$\cla$ which have relative position $w$ at $x$. Then
$\on{MOp}_G(D_x)_{\cla,w} \simeq \on{Op}_G(D_x)_{\cla} \times
S_{w,\F'_{B,x}}$.

We wish to show that each map ${\mb b}_{\cla,w}$ is an isomorphism
between $\on{Conn}^{\on{reg}}_{D_x,\cla,w}$ and
$\on{MOp}_G(D_x)_{\cla,w}$. The following result is proved my joint
work \cite{FG} with D. Gaitsgory in a more general setting.

\begin{prop}    \label{isom w}
For each $w \in W$ the morphism ${\mb b}_{\cla,w}$ is an isomorphism
between the varieties $\on{Conn}^{\on{reg}}_{D_x,\cla,w}$ and
$\on{MOp}_G(D_x)_{\cla,w}$.
\end{prop}

\begin{proof}
First we observe that at the level of points the map defined by ${\mb
b}_{\cla,w}, w \in W$, from the union of
$\on{Conn}^{\on{reg}}_{D_x,\cla,w}, w \in W$, to
$\on{MOp}_G(D_x)_{\cla}$, is a bijection. Indeed, by \propref{map
beta} we have a map taking a Miura oper from $\on{MOp}_G(D_x)_{\cla}$,
considered as a Miura oper on the punctured disc $D_x^\times$, to a
connection $\ol\nabla$ on the $H$--bundle $\Omega^{\crho}$ over
$D_x^\times$. We have shown above that $\ol\nabla$ has regular
singularity at $x$ and that its residue is of the form
$-w(\cla+\crho)+\crho, w \in W$. Thus, we obtain a map from the set of
points of $\on{MOp}_G(D_x)_{\cla}$ to the union of
$\on{Conn}^{\on{reg}}_{D_x,\cla,w}, w \in W$, and by \propref{map
beta} it is a bijection.

It remains to show that if the Miura oper belongs to
$\on{MOp}_G(D_x)_{\cla,w}$, then the corresponding connection has
residue precisely $-w(\cla+\crho)+\crho$.

Thus, we are given a $G$--oper $(\F,\nabla,\F_B,\F'_B)$ of coweight
$\cla$. Let us choose a trivialization of the $B$--bundle $\F_B$. Then
the connection operator reads
\begin{equation}    \label{conn op oper}
\nabla = \pa_t + \sum_{i=1}^\ell t^{\langle \al_i,\cla \rangle} f_i +
       {\mb v}(t), \qquad {\mb v}(t) \in \bb[[t]].
\end{equation}
Suppose that the horizontal $B$--reduction $\F'_B$ of our Miura oper
has relative position $w$ with $\F_B$ at $x$ (see \secref{miura opers}
for the definition of relative position). We need to show that the
corresponding connection on $\F'_H \simeq \Omega^{\crho}$ has residue
$-w(\cla+\crho)+\crho$.

This is equivalent to the following statement. Let $\Phi(t)$ be the
$G$--valued solution of the equation
\begin{equation}    \label{again equation}
\left( \pa_t + \sum_{i=1}^\ell t^{\langle \al_i,\cla \rangle} f_i +
     {\mb v}(t) \right) \Phi(t) = 0,
\end{equation}
such that $\Phi(0) = 1$. Since the connection operator is regular at
$t=0$, this solution exists and is unique. Then $\Phi(t) w^{-1} w_0$
is the unique solution of the equation \eqref{again equation} whose
value at $t=0$ is equal to $w^{-1} w_0$.

By \lemref{is generic}, we have
$$
\Phi(t) w^{-1} w_0 = X_w(t) Y_w(t) Z_w(t) w_0,
$$
where
$$
X_w(t) \in N((t)), \qquad Y_w(t) \in H((t)), \qquad Z_w(t) \in
N_-((t)).
$$
We can write $Y_w(t) = \cmu_w(t) \wt{Y}_w(t)$, where $\cmu_w$ is a
coweight and $\wt{Y}_w(t) \in H[[t]]$.

Since the connection $\nabla$
preserves $$\Phi(t) w_0 \bb_+ w_0 \Phi(t)^{-1} = \Phi(t)
\bb_- \Phi(t)^{-1},$$ the connection $X(t)_w^{-1} \nabla
X_w(t)$ preserves $$Y_w(t) Z_w(t) \bb_- Z_w(t)^{-1} Y_w(t)^{-1} =
\bb_-,$$ and therefore has the form
$$
\pa_t + \sum_{i=1}^\ell t^{\langle \al_i,\cla \rangle} f_i -
\frac{\cmu_w}{t} + {\mb u}(t), \qquad {\mb u}(t) \in \h[[t]].
$$
By conjugating it with $\cla(t)$ we obtain a connection
$$
\pa_t + p_{-1} - \frac{\cla+\cmu_w}{t} + {\mb u}(t), \qquad {\mb
  u}(t) \in \h[[t]].
$$
Therefore we need to show that
\begin{equation}    \label{desired}
\cmu_w = w(\cla+\crho) - (\cla+\crho).
\end{equation}

To see that, let us apply the identity $\Phi(t) w^{-1} = X_w(t) Y_w(t)
Z_w(t)$ to a non-zero vector $v_{w_0(\nu)}$ of weight $w_0(\nu)$ in a
finite-dimensional irreducible $\g$--module $V_\nu$ of highest weight
$\nu$ (so that $v_{w_0(\nu)}$ is a lowest weight vector and hence is
unique up to scalar). The right hand side will then be equal to a
$P(t) v_{w_0(\nu)}$ plus the sum of terms of weights greater than
$w_0(\nu)$, where $P(t) = c t^{\langle w_0(\nu),\cmu_w \rangle}, c
\neq 0$, plus the sum of terms of higher degree in $t$. Applying the
left hand side to $v_{w_0(\nu)}$, we obtain $\Phi(t) v_{w^{-1}
w_0(\nu)}$, where $v_{w^{-1} w_0(\nu)} \in V_\nu$ is a non-zero vector
of weight $w^{-1} w_0(\nu)$ which is also unique up to a scalar.

Thus, we need to show that the coefficient with which $v_{w_0(\nu)}$
enters $\Phi(t) v_{w^{-1} w_0(\nu)}$ is a polynomial in $t$ whose
lowest degree is equal to
$$
\langle w_0(\nu),w(\cla+\crho) - (\cla+\crho) \rangle,
$$
because if this is so for all dominant integral weights $\nu$, then we
obtain the desired equality \eqref{desired}. But this formula is easy
to establish. Indeed, from the form \eqref{conn op oper} of the oper
connection $\nabla$ it follows that we can obtain a vector
proportional to $v_{w_0}$ by applying the operators $t^{\langle
\al_i,\cla \rangle + 1} f_i, i=1,\ldots,\ell$, to $v_{w^{-1}
w_0(\nu)}$ in some order. The linear combination of these monomials
appearing in the solution is the term of the lowest degree in $t$ with
which $v_{w_0(\nu)}$ enters $\Phi(t) v_{w^{-1} w_0(\nu)}$. It follows
from \lemref{is generic} that it is non-zero. The corresponding power
of $t$ is nothing but the difference between the
$(\cla+\crho)$--degrees of the vectors $v_{w^{-1}w_0}$ and $v_{w_0}$,
i.e.,
$$
\langle w^{-1} w_0(\nu),\cla+\crho \rangle - \langle
w_0(\nu),\cla+\crho \rangle = \langle w_0(\nu),w(\cla+\crho) -
(\cla+\crho) \rangle,
$$
as desired. This completes the proof.
\end{proof}

Suppose we are given a Miura oper on $D_x$ with $\cla=0$ that has
relative position $s_i$ at $x$. Then the corresponding connection on
$\Omega^{\crho}$ has residue $- s_i(\crho) +\crho = \chal_i$. Choosing
a coordinate $t$ at $x$, we write this connection as
\begin{equation}    \label{alpha i}
\ol{\nabla} = \pa_t + \frac{\chal_i}{t} + {\mb u}(t), \qquad {\mb
u}(t) \in \h[[t]].
\end{equation}

\begin{lem}    \label{si}
A connection of the form \eqref{alpha i} belongs to
$\on{Conn}^{\on{reg}}_{D_x,s_i}$ (i.e., the
corresponding $G$--oper is regular at $x$) if and only if $\langle
\al_i,{\mb u}(0) \rangle = 0$.
\end{lem}

\begin{proof}
Let $V_{\omega_i}$ be the $i$th fundamental representation of $\g$. It
contains a one-dimensional subspace $L_{\omega_i}$ of
$B$--invariants. There is a canonical two-dimensional subspace
$W_{\omega_i}$ of $V_{\omega_i}$ stable under $B$, containing
$L_{\omega_i}$, and on which the $SL_2$ subgroup corresponding to the
$i$th simple root acts irreducibly. Moreover, the generators $f_j, j
\neq i$, act on $W_{\omega_i}$ by $0$. Consider the vector bundle
$$
V_{\omega_i,\F} = \F \underset{G}\times V_{\omega_i} = \F_B
\underset{B}\times V_{\omega_i}
$$
and the corresponding rank two subbundle $W_{\omega_i,\F_B}$. The
connection $\nabla$ preserves $W_{\omega_i,\F_B}$, and its restriction
to $W_{\omega_i,\F_B}$ is equal to
$$
\pa_t + \begin{pmatrix} \frac{1}{t} + \frac{1}{2} u_i(t) & 0 \\
1 & - \frac{1}{t} - \frac{1}{2} u_i(t) \end{pmatrix},
$$
where $u_i(t) = \langle \al_i,{\mb u}(t) \rangle$.

The corresponding equation $\nabla \Phi = 0$ has two linearly
independent solutions:
$$
\Phi_1 = \begin{pmatrix} 0 \\ t e^{\int u_i(t) dt} \end{pmatrix},
\qquad \Phi_2 = \begin{pmatrix} - t^{-1} e^{- \int u_i(t) dt} \\ t
  e^{\int u_i(t) dt} \int t^{-2} e^{-2 \int u_i(t) dt} dt
  \end{pmatrix}.
$$
Hence the monodromy of these solutions is equal to $$\begin{pmatrix} 1 &
  -4 \pi i u_i(0) \\ 0 & 1 \end{pmatrix}.$$ This implies that this
oper is regular only if $u_i(0) = 0$.

Conversely, if $u_i(0) = 0$, then applying to the connection $\nabla$
the gauge transformation with $\exp(-e_i/t)$ we obtain a regular
connection. Hence the corresponding oper is regular. This completes
the proof.
\end{proof}

The above calculation also implies that the scheme-theoretic preimage
of $\on{Op}_G(D_x)$ under the morphism $\on{Conn}^{\on{RS}}_{D_x,s_i}
\to \on{Op}_G^{\on{RS}}(D_x)$ is in fact reduced, and therefore it is
equal to $\on{Conn}^{\on{reg}}_{D_x,s_i}$.

\section{Bethe Ansatz equations and Miura opers on $\pone$}
\label{second}

In this section we consider Miura opers and the corresponding
$H$--connections on $\pone$. We show that the Miura opers having the
simplest possible degenerations are described by the solutions of the
so-called Bethe Ansatz equations. This will allow us eventually to
describe the set of solutions of the Bethe Ansatz equations as an open
dense subset of the flag variety.

\subsection{Miura opers on $\pone$}

Let us fix a set of distinct points $z_1,\ldots,z_N$ on $\pone$ such
that $z_i \neq \infty$ for all $i=1,\ldots,N$, and a set of dominant
coweights $\cla_1,\ldots,\cla_N$ of $\g$. Let
$\on{Op}^{\on{RS}}_{G}(\pone)_{(z_i),\infty;(\cla_i),\cla_\infty}$ be
the set of $G$--opers on $\pone$ which are regular at all points other
than $z_1,\ldots,z_N,\infty$ and have regular singularities at
$z_1\ldots,z_N,\infty$ with the residues
$\cla_1,\ldots,\cla_N,\cla_\infty$. More precisely, this is the subset
of the set $\on{Op}_{G}(\pone \bs \{ z_1,\ldots,z_N,\infty \})$
consisting of those opers whose restriction to the punctured disc
$D_{z_i}^\times$ at the point $z_i$ (resp., $D^\times_\infty$) belongs
to $\on{Op}^{\on{RS}}_{G}(D_{z_i})_{\cla_i}$ for all $i=1,\ldots,N$
(resp., to $\on{Op}^{\on{RS}}_{G}(D_\infty)_{\cla_\infty}$).

Let $\cla_\infty$ be another dominant coweight of $\g$. Introduce a
subset $$\on{Op}_{G}(\pone)_{(z_i),\infty;(\cla_i),\cla_\infty}
\subset
\on{Op}^{\on{RS}}_{G}(\pone)_{(z_i),\infty;(\cla_i),\cla_\infty}$$ of
those $G$--opers whose restriction to $D_{z_i}^\times$ belongs to
$\on{Op}_{G}(D_{z_i})_{\cla_i} \subset
\on{Op}^{\on{RS}}_{G}(D_{z_i})_{\cla_i}$ for all $i=1,\ldots,N$ and
whose restriction to $D^\times_\infty$ belongs to
$\on{Op}_{G}(D_{\infty})_{\cla_\infty} \subset
\on{Op}^{\on{RS}}_G(D_\infty)_{\cla_\infty}$. Denote by
$\on{MOp}_{G}(\pone)_{(z_i),\infty;(\cla_i),\cla_\infty}$ the space of
Miura opers on $\pone \bs \{ z_1,\ldots,z_N,\infty \}$ whose
underlying opers belong to
$\on{Op}_{G}(\pone)_{(z_i),\infty;(\cla_i),\cla_\infty}$.

Let $\tau = (\F,\nabla,\F_{B})$ be an oper from
$\on{Op}_{G}(\pone)_{(z_i),\infty;(\cla_i),\cla_\infty}$. The above
conditions mean that the oper bundle $\F$, which is a priori defined
on $\pone \bs \{ z_1,\ldots,z_N,\infty \}$, has a canonical extension
to the entire $\pone$. By \lemref{no mon}, the monodromy around each
of the points $z_1,\ldots,z_N,\infty$ is trivial. Therefore the flat
connection $\nabla$ has the trivial monodromy representation and
therefore defines a global trivialization of the oper bundle
$\F$. Hence, by \lemref{isom with flags}, the space
$\on{MOp}_{G}(\pone)_\tau$ of Miura $G$--opers on $\pone$ whose
underlying oper is $\tau$ is isomorphic to the flag variety $G/{}B$ of
$G$. Indeed, any $B$--reduction of the fiber of this bundle at an
arbitrary point $x$ of $\pone$ uniquely extends to a horizontal
$B$--reduction of $\F$ on the entire $\pone$. But a $B$--reduction of
the fiber of the trivial bundle at $x$ is the same as a point of
$(G/{}B)_{\F_x}$ which is isomorphic to $G/B$.

On the other hand, let
$\on{Conn}(\pone)_{(z_i),\infty;(\cla_i),\cla_\infty}^{\on{RS}}$ be
the space of connections on the $H$--bundle $\Omega^{\crho}$ on
$\pone$ with regular singularities at the points
$z_1,\ldots,z_N,\infty$ and a finite number of other points
$w_1,\ldots,w_m$ such that the residue at $z_i$ (resp., $\infty$,
$w_j$) is equal to $-y_i(\cla_i+\crho)+\crho$ (resp.,
$-y_\infty(\cla_\infty+\crho)+\crho, -y'_j(\crho)+\crho$) for some
elements $y_i,y_\infty,y'_j \in W$. Such a connection then has the
form
\begin{equation}    \label{conn RS}
\pa_t - \sum_{i=1}^N \frac{y_i(\cla_i+\crho)-\crho}{t-z_i}
- \sum_{j=1}^m \frac{y'_j(\crho)-\crho}{t-w_j}
\end{equation}
on $\af = \pone \bs \infty$. According to formula \eqref{trans for
conn}, connection $\pa_t + f(t)$ on $\Omega^{\crho}$ over $\af$ has
the following expansion on the disc around $\infty \in \pone$ with
respect to the coordinate $u=t^{-1}$:
$$
\pa_u - u^{-2} f(u^{-1}) + 2\crho u^{-1}.
$$
Therefore the residue of the connection \eqref{conn RS} at $\infty$
is equal to
$$
2\crho + \sum_{i=1}^N (y_i(\cla_i+\crho)-\crho) +
\sum_{j=1}^m (y'_j(\crho)-\crho).
$$
On the other hand, by our assumption, it should be equal to
$-y_\infty(\cla_\infty+\crho)+\crho$ for some $y_\infty \in
W$. Denoting $y_\infty$ by $y'_\infty w_0$, we obtain the following
equation relating the residues of our connection:
\begin{equation}    \label{relation}
\sum_{i=1}^N (y_i(\cla_i+\crho)-\crho) + \sum_{j=1}^m
(y'_j(\crho)-\crho) = y'_\infty(-w_0(\cla_\infty) +
\crho) - \crho.
\end{equation}

To a connection of this form we associate a $G$--oper on $\pone$ with
regular singularities at $(z_i),\infty,(w_j)$ in the same way as
above. Namely, we set $$\F = \Omega^{\crho} \underset{H}\times G,
\qquad \F_{B} = \Omega^{\crho} \underset{H}\times {}B, \qquad \nabla =
\ol{\nabla} + p_{-1}.$$ This oper belongs to the set
$$\on{Op}^{\on{RS}}_G(\pone)_{(z_i),(w_j),\infty;(\cla_i),(0),\cla_\infty}$$
defined in the same way as at the beginning of this section. Note that
by construction the residue of this oper at $z_i$ is in the $W$--orbit
of $\cla_i$, whereas at $w_j$ it is in the $W$--orbit of $0$.

Thus, we have a map
$$
\on{Conn}(\pone)^{\on{RS}}_{(z_i),\infty;(\cla_i),\cla_\infty} \to
\on{Op}^{\on{RS}}_G(\pone)_{(z_i),(w_j),\infty;(\cla_i),(0),\cla_\infty}.
$$
Let $\on{Conn}(\pone)_{(z_i),\infty;(\cla_i),\cla_\infty}$ be the
subset of
$\on{Conn}(\pone)^{\on{RS}}_{(z_i),\infty;(\cla_i),\cla_\infty}$
consisting of those connections for which the resulting oper $\tau$ on
$\pone$ belongs to
$\on{Op}_{G}(\pone)_{(z_i),\infty;(\cla_i),\cla_\infty}$. The
resulting map
$$
\ol{\mb b}_{(z_i),\infty;(\cla_i),\cla_\infty}:
\on{Conn}(\pone)_{(z_i),\infty;(\cla_i),\cla_\infty} \to
\on{Op}_{G}(\pone)_{(z_i),\infty;(\cla_i),\cla_\infty}
$$
may be lifted to a map
$${\mb b}_{(z_i),\infty;(\cla_i),\cla_\infty}:
\on{Conn}(\pone)_{(z_i),\infty;(\cla_i),\cla_\infty} \to
\on{MOp}_{G}(\pone)_{(z_i),\infty;(\cla_i),\cla_\infty}.$$ Here
$\on{MOp}_{G}(\pone)_{(z_i),\infty;(\cla_i),\cla_\infty}$ is the space
of Miura opers on $\pone \bs \{ z_1,\ldots,z_N,\infty \}$ such that
the underlying oper belongs to
$\on{Op}_{G}(\pone)_{(z_i),\infty;(\cla_i),\cla_\infty}$. Namely, we
give an oper $\tau$ that is in the image of $\ol{\mb
  b}_{(z_i),\infty;(\cla_i),\cla_\infty}$ the structure of a Miura
oper by defining a horizontal $B$--reduction $\F'_{B}$ of $\F$ by the
formula $$\F'_{B} = \Omega^{\crho} \underset{H}\times {}w_0 B.$$

Next, we construct the map
$$
{\mb a}_{(z_i),\infty;(\cla_i),\cla_\infty}:
\on{MOp}_{G}(\pone)_{(z_i),\infty;(\cla_i),\cla_\infty} \to
\on{Conn}(\pone)_{(z_i),\infty;(\cla_i),\cla_\infty}.
$$
Any Miura oper on $\pone \bs \{ z_1,\ldots,z_N,\infty \}$ from
$\on{MOp}_{G}(\pone)_{(z_i),\infty;(\cla_i),\cla_\infty}$ becomes
generic after removing finitely many points $w_j$. Hence it gives rise
to a connection $\ol\nabla$ on the bundle $\F_H \simeq \Omega^{\crho}$
over $\pone \bs \{ (z_i),(w_j),\infty \}$. The restrictions of this
connection to the discs around the points $z_i$, $w_j$ and $\infty$
must have regular singularities with the residues being in the
$W$--orbits of $\cla_i$, $0$, and $\cla_\infty$,
respectively. Therefore this connection must be of the form
\eqref{conn RS}. This defines a map ${\mb
  a}_{(z_i),\infty;(\cla_i),\cla_\infty}$. In the same way as in
\propref{map beta} and \propref{isom w} we show that the maps ${\mb
  a}_{(z_i),\infty;(\cla_i),\cla_\infty}$ and $\ol{\mb
  b}_{(z_i),\infty;(\cla_i),\cla_\infty}$ are mutually inverse
bijections.

Let us fix an oper $$\tau \in
\on{Op}_{G}(\pone)_{(z_i),\infty;(\cla_i),\cla_\infty}$$ and
trivialize the underlying $G$--bundle $\F$ by identifying the fiber at
the point $\infty \in \pone$ with $G$. Then the connection trivializes
the bundle $\F$ by identifying all fibers with the fiber at $\infty$
and hence with $G$. Therefore we also obtain a trivialization of the
corresponding $G/{}B$--bundle, and so the reduction $\F_{B}$ gives us
a map $\phi_\tau: \pone \to G/{}B$.

Note that giving $\tau$ the structure of a Miura oper amounts to
picking a point in the flag variety $G/{}B$. If this point belongs to
the $B$--orbit $S_{y_\infty} = B y_\infty^{-1} w_0 B \subset G/B$
(i.e., if the corresponding $B$--reduction and the oper reduction of
the fiber at $\infty$ are in relative position $y_\infty$; see the
definition in \secref{miura opers}), then the corresponding connection
has the residue $-y_\infty(\cla_\infty+\crho)+\crho$ at
$\infty$. Furthermore, identifying the fiber of $\F$ at $\infty$ with
the fiber of $\F$ at $z_i$, we obtain a reduction of $\F_{z_i}$ to
$B$. According to \propref{isom w}, this reduction then has relative
position $y_i$ with the oper reduction $\F_{B,z_i}$ precisely when the
residue of our connection at $z_i$ is equal to
$-y_i(\cla_i+\crho)+\crho$.

Likewise, the points $w_j$'s are the points where our $B$--reduction
is not in generic position with the oper reduction, and it is then in
relative position $y'_j$ precisely when the residue of our connection
at $w_j$ is equal to $-y'_j(\crho)+\crho$, by \propref{isom w}. All of
these residues must satisfy the relation \eqref{relation}. Thus, we
obtain the following result.

\begin{thm}    \label{strongest}
The map ${\mb b}_{(z_i),\infty;(\cla_i),\cla_\infty}$ is a bijection
at the level of points. Thus, the set of all connections from
$\on{Conn}(\pone)_{(z_i),\infty;(\cla_i),\cla_\infty}$ which
correspond to a fixed $G$--oper $\tau \in
\on{Op}_{G}(\pone)_{(z_i),\infty;(\cla_i),\cla_\infty}$ is isomorphic
to the set of points of the flag variety $G/{}B$.

Moreover, the residues of these connections at the points $z_i$
(resp., the points $w_j$) are equal to $R_{z_i} = -
y_i(\cla_i+\crho) + \crho$ (resp., $R_{w_j} = - y'_j(\crho) + \crho$)
for some elements $y_i, y'_j \in W$ and they must satisfy the relation
$$
\sum_{i=1}^n R_{z_i} + \sum_{j=1}^m R_{w_j} = -
y'_\infty(-w_0(\cla_\infty)+\crho) + \crho
$$
for some $y'_\infty \in W$. The set of those connections which satisfy
this relation is in bijection with the Schubert cell $B w_0
y'_\infty w_0 B$ in $G/B$.
\end{thm}

\subsection{Bethe Ansatz equations}    \label{miura opers on pone}

Let $\tau$ be again an oper from
$\on{Op}_{G}(\pone)_{(z_i),\infty;(\cla_i),\cla_\infty}$ and suppose
we have a Miura oper projecting onto $\tau$ under the Miura
transformation. Let $\phi_\tau: \pone \to G/B$ be the map
corresponding to the reduction $\F_B$. Recall that its value at $x \in
\pone$ is $\F_{B,x}$ considered as a point of $(G/B)_{\F'_B,x} \simeq
(G/B)_{\F'_B,\infty} \simeq G/B$, where the first isomorphism is
obtained from the identification of the fibers of $\F'_B$ induced by
the oper connection and the second isomorphism corresponds to a choice
of trivialization of $(G/B)_{\F'_B,\infty}$. Consider the subvariety
$(G/{}B)_\tau$ of $G/{}B$ whose points $p$ satisfy the following
conditions:
\begin{itemize}
\item[(1)] $\phi_\tau(z_i)$ is in generic position with
$p$ for all $i=1,\ldots,N$;

\item[(2)] the relative position of $\phi_\tau(x)$ and $p$ is either
generic or corresponds to a simple reflection $s_i \in W$ for all $x
\in \pone \bs \{ z_1,\ldots,z_N,\infty \}$.
\end{itemize}

It is clear that $(G/{}B)_\tau$ is an open and dense subvariety of
$G/{}B$. Indeed, $(G/{}B)_\tau$ is contained in the intersection
$U_\tau$ of finitely many open and dense subsets, namely, the sets of
points of $G/B$ which are in generic relative position with
$\phi_\tau(z_i)$ (each is isomorphic to the big Schubert cell). The
complement of $(G/{}B)_\tau$ in $U_\tau$ is a subvariety of
codimension one. This
subvariety consists of all points in $G/B$ which are in relative
position $w$ with $\phi_\tau(x), x
\in \pone \bs \{ z_1,\ldots,z_N,\infty \}$, where $w$ runs over the
subset of $W$ of all elements of length $l(w) \geq 2$. The subvariety
of these points for fixed $x$ has codimension
two, and therefore their
union, as $x$ moves along the curve $\pone \bs \{
z_1,\ldots,z_N,\infty \}$, has codimension (at least) one.

Let
$$\on{MOp}_{G}(\pone)^{\on{gen}}_{(z_i),\infty;(\cla_i),\cla_\infty}
\subset \on{MOp}_{G}(\pone)_{(z_i),\infty;(\cla_i),\cla_\infty}$$ be
the open dense subvariety which is the union of $(G/{}B)_\tau, \tau
\in \on{Op}_{G}(\pone)_{(z_i),\infty;(\cla_i),\cla_\infty}$.

Note that equation \eqref{relation} now reads
\begin{equation}    \label{special relation}
\sum_{i=1}^N \cla_i - \sum_{j=1}^m \chal_{i_j} =
y(-w_0(\cla_\infty)+\crho) - \crho,
\end{equation}
where we write $y = y'_\infty$ to simplify notation.

Consider the image of
$\on{MOp}_{G}(\pone)^{\on{gen}}_{(z_i),\infty;(\cla_i),\cla_\infty}$ in
$\on{Conn}(\pone)_{(z_i),\infty;(\cla_i),\cla_\infty}$ under the
bijection ${\mb b}_{(z_i),\infty;(\cla_i),\cla_\infty}$. We denote it
by $\on{Conn}(\pone)^{\on{gen}}_{(z_i),\infty;(\cla_i),\cla_\infty}$.
Then according to \lemref{si},
$\on{Conn}(\pone)^{\on{gen}}_{(z_i),\infty;(\cla_i),\cla_\infty}$ is
precisely the set of all connections of the form
\begin{equation}    \label{main conn}
\ol{\nabla} = \pa_t - \sum_{i=1}^N \frac{\cla_i}{t-z_i} + \sum_{j=1}^m
\frac{\chal_{i_j}}{t-w_j},
\end{equation}
where $w_1,\ldots,w_m$ are points of $\pone \bs \{
z_1,\ldots,z_N,\infty \}$ and $i_j \in I$ for all $j=1,\ldots,m$, such
that if
$$
\pa_t + \chal_{i_j}/(t-w_j) + {\mb u}_j(t-w_j), \qquad {\mb u}_j(u)
\in \h[[u]],
$$
is the expansion of the connection \eqref{main conn} at the point
$w_j$, then $\langle \al_{i_j},{\mb u}_j(0) \rangle = 0$ for all
$j=1,\ldots,m$. Explicitly, these equations read
\begin{equation}    \label{bethe}
\sum_{i=1}^N \frac{\langle \al_{i_j},\cla_i \rangle}{w_j-z_i} -
\sum_{s \neq j} \frac{\langle \al_{i_j},\chal_{i_s} \rangle}{w_j-w_s}
= 0, \qquad j=1,\ldots,m.
\end{equation}
They are called the {\em Bethe Ansatz equations}. We have an obvious
action of a product of symmetric groups permuting the points $w_j$
corresponding to simple roots of the same kind. In what follows, by a
{\em solution} of the Bethe Ansatz equations we will understand a
solution defined up to these permutations. We will adjoin to the
set of all solutions associated to all possible collections $\{
\al_{i_j} \}$ of simple roots of $\g$, the unique ``empty'' solution,
corresponding to the empty set of simple roots.

Now we obtain the following

\begin{thm}    \label{main}
The set of solutions of the Bethe Ansatz equations \eqref{bethe} is in
bijection with the set of points of
$\on{MOp}_{G}(\pone)^{\on{gen}}_{(z_i),\infty;(\cla_i),\cla_\infty}$.
\end{thm}

Let us fix an oper $\tau$. Then considering the value of the map
$\phi_\tau$ at $\infty \in \pone$, we obtain an identification of the
space
$\on{MOp}_{G}(\pone)^{\on{gen}}_{(z_i),\infty;(\cla_i),\cla_\infty}$
with an open dense subset of $G/{}B$. Furthermore, it follows from
\thmref{strongest} that those elements of
$\on{MOp}_{G}(\pone)^{\on{gen}}_{(z_i),\infty;(\cla_i),\cla_\infty}$
which satisfy formula \eqref{special relation} correspond to points
that lie in the Schubert cell $B w_0 y w_0 B$ in $G/B$.

Note that except for the big cell $B w_0 B$, the intersection between
the Schubert cell $B w_0 y w_0 B$ and the open dense subset
$(G/B)_\tau \subset G/B$ could be either an open dense subset of $B
w_0 y w_0 B$ or empty.\footnote{For example, it follows from the
results of Mukhin and Varchenko in \cite{MV:new} that sometimes this
open set may not contain the one point Schubert cell $B \subset G/B$
even if we allow $z_1,\ldots,z_N$ to be generic.} Therefore we obtain

\begin{cor}    \label{fixed oper}
The set of those solutions of the Bethe Ansatz equations which
correspond to a fixed $G$--oper $\tau \in
\on{Op}_{G}(\pone)_{(z_i),\infty;(\cla_i),\cla_\infty}$ is in
bijection with the set of points of an open and dense subset
$(G/B)_\tau$ of the flag variety $G/{}B$. Further, every solution must
satisfy the equation \eqref{special relation} for some $y \in W$, and
the solutions which satisfy this equation with fixed $y \in W$ are in
bijection with an open subset of the Schubert cell $B w_0 y w_0 B \in
G/B$.
\end{cor}

\subsection{The action of $N$ on solutions of the Bethe Ansatz
equations}    \label{action of N}

The group $N$ naturally acts on $G/{}B$, and thus we obtain an action
of $N$ on the set of those solutions of the Bethe Ansatz equations
which correspond to a fixed $G$--oper. This action is however
rational, because solutions of the Bethe Ansatz equations correspond
to points of an open dense subset of $G/B$, not the entire $G/B$.

Let us identify the set of solutions of the Bethe Ansatz equations
with an open dense subset of the flag variety by using the fiber at $0
\in \pone$, instead of $\infty \in \pone$. Then the action of $g \in
N$ is given by gauge transformations on a connection of the form
$$
\pa_t + p_{-1} + {\mb u}(t), \qquad {\mb u}(t) \in {}\h,
$$
by a rational $N$--valued function $g(t)$ such that $g(0) = g$ and
\begin{equation}    \label{g conj}
g(t)(\pa_t + p_{-1} + {\mb u}(t))g(t)^{-1} = \pa_t + p_{-1} + \wt{\mb
u}(t),
\end{equation}
where $\wt{\mb u}(t)$ is again in $\h$. Clearly, $g(t)$ is uniquely
determined by these conditions.

By our assumptions, the connection $\pa_t + p_{-1} + {\mb
u}(t)$ has trivial monodromy representation. Therefore there
exists a (unique) polynomial $G$--valued solution $\Phi(t)$ of the
equation $(\pa_t + p_{-1} + {\mb u}(t)) \Phi(t) = 0$ with the initial
condition $\Phi(0) = 1$. Because of the form of the connection, we
find that $\Phi(t)$ actually takes values in $B_-$. Further, for any
constant element $M$ of $G$, the solution of the equation $(\pa_t +
p_{-1} + {\mb u}(t)) \Psi(t) = 0$ with the initial condition $\Psi(0)
= M$ is $\Psi(t) = \Phi(t) M$.

Now if $\wt{\Phi}(t)$ is the solution of the equation $(\pa_t + p_{-1}
+ \wt{\mb u}(t)) \wt{\Phi}(t) = 0$ with the initial condition
$\wt{\Phi}(0) = 1$ (like $\Phi(t)$, it takes values in $B_-$), we
obtain the following equation:
\begin{equation}    \label{factorization}
\Phi(t) g^{-1} = g(t)^{-1} \wt{\Phi}(t).
\end{equation}
Thus, to find $g(t)$, we need to find $\Phi(t)$ and to project the
function $\Phi(t) g^{-1}$ onto $N$ considered as an open dense subset
of $G/B_-$ (in general, this may only be done for generic values of $t$).

Let us consider more explicitly the case when $g = \exp(a e_i), i \in
I, a \in \C$ (these one-parameter subgroups generate the action of the
group $N$). We claim that $g(t)$ is then necessarily of the form $g(t)
= \exp(f(t)e_i)$, where $f(t)$ is a rational function in $t$ such that
$f(0) = a$.

Indeed, since $\wt{\Phi}(t) \in {}B_-$, we obtain that the left hand
side of \eqref{factorization} belongs to the $i$th minimal parabolic
subgroup of $G$ generated by $B_-$ and the $SL_2$ subgroup
corresponding to the $i$th simple root. Hence the left hand side must
also belong to this parabolic subgroup, and therefore $g(t)$
necessarily has the form $\exp(f(t) e_i)$ for some rational function
$f(t)$ satisfying $f(0) = a$.

Let us compute $f(t)$. We have
\begin{multline*}
\exp(f(t)e_i)(\pa_t + p_{-1} + {\mb u}(t))\exp(-f(t)e_i) = \\ \pa_t +
({\mb u}(t) + f(t) \chal_i) - (f'(t) + f(t)^2 + f(t) u_i(t)) e_i,
\end{multline*}
where $u_i(t) = \langle \al_i,{\mb u}(t) \rangle$. Therefore
$f(t)$ has to be a rational solution of the differential equation
\begin{equation}    \label{diff eq}
f'(t) + f(t)^2 + f(t) u_i(t) = 0
\end{equation}
with the initial condition $f(0) = a$. Note that from the previous
discussion we already know that such a solution exists and is
unique. Then in the new connection we will have $\wt{\mb u}(t) = {\mb
u}(t) + f(t) \chal_i$. For generic values of $a$ the function $\wt{\mb
u}(t)$ will have the same form as ${\mb u}(t)$. In particular, only
the positions of the poles $w_j$ such that $i_j=i$ will be changed in
${\mb u}(t)$, and hence its poles will again give us a solution of the
Bethe Ansatz equations \eqref{bethe}.

Thus we obtain a rational action of the elements of the form $\exp(a
e_i)$ on the set of solutions of the Bethe Ansatz equations. These
actions generate a rational action of the group $N$ on the set of
solutions of the Bethe Ansatz equations corresponding to a fixed
$G$--oper $\tau$. By our construction, this action becomes the natural
action of $N$ on the flag variety $G/{}B$ under the embedding of the
set of solutions of \eqref{bethe} into the flag variety as an open
dense subset.

\subsection{Comparison with results of Mukhin and Varchenko}
\label{comparison}

In \cite{MV}, Mukhin and Varchenko associated to each solution $\{
w_1,\ldots,w_m \}$ of the Bethe Ansatz equations \eqref{bethe} an
$I$--tuple of polynomials $(y_i(x)), i \in I = \{ 1,\ldots,\ell \}$,
each defined up to a scalar, such that the roots of $y_i(x)$ are
precisely those $w_j$'s for which $i_j=i$. Thus, they obtained an
embedding of the set of solutions of \eqref{bethe} into the product of
$\ell$ copies of ${\mathbb P}(\C[x])$. Then they defined the ``$i$th
reproduction procedure'' of solutions as follows. Set
$$
T_i(x) = \prod_{j=1}^N (x-z_j)^{\langle \cla_j,\al_i \rangle}
$$
and let $\wt{y}_i$ be a new polynomial that has the form
\begin{equation}    \label{mv}
\wt{y}_i(x) = y_i(x) \int^x T_i(t) \prod_{j \in I} y_j(t)^{-\langle
\chal_j,\al_i \rangle} dt.
\end{equation}
The closure of the set of polynomials of this form in ${\mathbb
P}(\C[x])^\ell$ is isomorphic to a projective line. It is proved in
\cite{MV} that for all but finitely many points of this line the
$I$--tuple $(y_1(x),\ldots,\wt{y}_i(x),\ldots,y_\ell(x))$ will again
correspond to a solution of the Bethe Ansatz equations \eqref{bethe},
and thus they obtain a rational map from $\pone$ to the set of
solutions of \eqref{bethe}.

Let us show that the image of this map coincides with the closure of
the orbit of the group $\{ \exp(a e_i) \}$ acting on the set of
solutions of \eqref{bethe} as explained above (this observation is
due to Mukhin, Varchenko and myself). Indeed, given an $I$--tuple
$(y_i(x))_{i \in I}$ encoding a solution of the equations
\eqref{bethe}, the corresponding connection is given by the formula
$\pa_t + p_{-1} + {\mb u}(t)$, where
$$
{\mb u}(t) = - \sum_{j=1}^N \frac{\cla_j}{t-z_j} + \sum_{i \in I}
\chal_i \frac{d}{dt} \log y_i(t),
$$
so that
$$
u_i(t) = \langle \al_i,{\mb u}(t) \rangle = - \frac{d}{dt} \log
\left( T_i(t) \prod_{j \in I} y_j(t)^{-\langle \chal_j,\al_i
\rangle} \right).
$$
But \eqref{mv} implies that $\wt{y}_i(t)$ satisfies the equation
$$
\frac{d}{dt} \log \wt{y}_i(t) = \frac{d}{dt} \log y_i(t) +
\frac{d}{dt} \log \int^t T_i(x) \prod_{j \in I} y_j(x)^{-\langle
\chal_j,\al_i \rangle} dx.
$$
Therefore the connection corresponding to the new $I$--tuple
$(y_1(x),\ldots,\wt{y}_i(x),\ldots,y_\ell(x))$ has the form $\pa_t +
p_{-1} + \wt{\mb u}(t)$, where $\wt{\mb u}(t) = {\mb u}(t) + f(t)
\chal_i$, and $f(t)$ satisfies the differential equation \eqref{diff
eq}.

Therefore the $i$th reproduction procedure of \cite{MV} on the set of
solutions of \eqref{bethe} coincides with the action of the group $\{
\exp(a e_i) \}$ on the image of this set in $G/{}B$. Mukhin and
Varchenko use the reproduction procedures to construct ``populations''
of solutions of \eqref{bethe}. A population is by definition the
closure in $({\mathbb P}(\C[x]))^\ell$ of the set of all polynomials
obtained by applying consecutively all possible reproduction
procedures to an $I$--tuple of polynomials corresponding to a
particular solution of \eqref{bethe}. Conjecture 3.10 in \cite{MV}
then asserts that each population is isomorphic to the flag variety
$G/B$ and that the subset of the population corresponding to the
$I$--tuples of polynomials of fixed degrees is isomorphic to a Schubert
cell in $G/{}B$. This has been proved in \cite{MV} for $\g$ of types
$A_n, B_n$ and $C_n$ and in \cite{BM} for $G_2$ (note that in the
convention of \cite{MV} the Bethe Ansatz equations \eqref{bethe}
correspond to the Langlands dual group $^L G$, and so the relevant
flag manifold is $^L G/{}^L B$ rather than $G/B$).

Now this assertion immediately follows for an arbitrary simple Lie
algebra $\g$ from \corref{fixed oper} and the above
discussion. Indeed, we have found in \corref{fixed oper} that the set
of solutions of the Bethe Ansatz equations is identified with an open
dense subset of $G/{}B$ and that those solutions which correspond to
$I$--tuples of polynomials of fixed degrees correspond to points of
Schubert cells in $G/{}B$.  Moreover, we have identified the
reproduction procedures with the action of the one-parameter subgroups
$\{ \exp(a e_i) \}$ on this set inside $G/{}B$. But the closure of the
union of the consecutive orbits of these subgroups is equal to the
entire flag manifold $G/B$. Hence we obtain that any population of
solutions (in the terminology of \cite{MV}) is indeed isomorphic to
$G/B$.

\section{The Gaudin model and the Bethe Ansatz}    \label{third}

In this section we consider a simple Lie algebra $\g$ and its
Langlands dual Lie algebra $^L \g$ (whose Cartan matrix is the
transpose of that of $\g$). We will identify the set of roots of $\g$
with the set of coroots of $^L \g$ and the set of weights of $\g$ with
the set of coweights of $^L \g$. The results on opers and Miura opers
from the previous sections will be applied here to the Lie algebra $^L
\g$.

\subsection{The definition of the Gaudin model}

Here we recall the definition of the Gaudin model and the realization
of the Gaudin hamiltonians in terms of the spaces of conformal blocks
for affine Kac-Moody algebras of critical level. We follow closely the
paper \cite{FFR}.

For a dominant integral weight $\la$ denote by $V_\la$ the irreducible
representation of $\g$ of highest weight $\la$. Choose a
non-degenerate invariant inner product $\ka_0$ on $\g$. Let $\{ J_a
\}, a=1,\ldots,d$, be a basis of $\g$ and $\{ J^a \}$ the dual basis
with respect to $\ka_0$. Denote by $\Delta$ the quadratic Casimir
operator from the center of $U(\g)$:
$$\Delta = \frac{1}{2} \sum_{a=1}^d J_a J^a.$$

Let $\laa$ be a set of of dominant highest weights of $\g$. Denote by
$V_{\laa}$ the tensor product $V_{\la_1} \otimes \ldots \otimes
V_{\la_N}$. Let $z_1,\ldots,z_N$ be a set of distinct complex
numbers. The {\em Gaudin hamiltonians} are the linear operators
\begin{equation}    \label{hi}
\Xi_i = \sum_{j\neq i} \sum_{a=1}^d \frac{J_a^{(i)}
  J^{a(j)}}{z_i-z_j}, \qquad i=1,\ldots,N,
\end{equation}
acting on $V_{\laa}$. Note that $$\sum_{i=1}^N \Xi_i = 0.$$ These
operators commute with the diagonal action of $\g$ on $V_{\laa}$ and
hence their action is well-defined on the subspace of highest weight
vectors in $V_{\laa}$ of an arbitrary dominant integral weight $\mu$
with respect to the diagonal $\g$--action. Writing $\mu =
-w_0(\la_\infty)$ where $\la_\infty$ is another dominant integral
weight, we identify this subspace with $(V_{\laa} \otimes
V_{\la_\infty})^G$.

Consider the problem of simultaneous diagonalization of the Gaudin
hamiltonians. Set $$|0\ri := v_{\la_1} \otimes
\ldots \otimes v_{\la_N} \in V_{\laa}.$$ Clearly, it is an eigenvector
of the $\Xi$'s. Other eigenvectors are constructed by a procedure
known as the {\em Bethe Ansatz}.

Let
$$F_j(w) = \sum_{i=1}^N \frac{F_j^{(i)}}{w-z_i}, \quad \quad
j=1,\ldots,\ell.$$
For a set of distinct complex numbers $w_1,\ldots,w_m$ and a
collection of labels $i_1,\ldots,i_m \in I$ we introduce the Bethe
vector
\begin{equation}    \label{genbv}
|w_1^{i_1},\ldots,w_m^{i_m}\ri = \sum_{p=(I^1,\ldots,I^N)} \prod_{j=1}^N
\frac{F_{i^j_1}^{(j)} F_{i^j_2}^{(j)} \ldots
F_{i^j_{a_j}}^{(j)}}{(w_{i^j_1}-w_{i^j_2})(w_{i^j_2}-w_{i^j_3}) \ldots
(w_{i^j_{a_j}}-z_j)} |0\ri.
\end{equation}
Here the summation is taken over all {\em ordered} partitions $I^1 \cup
I^2 \cup \ldots \cup I^N$ of the set $\{1,\ldots,m\}$, where $I^j = \{
i^j_1,i^j_2,\ldots,i^j_{a_j} \}$. Note that one can consider vector
\eqref{genbv} as an element of the tensor product of Verma modules
$M_{\la_1} \otimes \ldots \otimes M_{\la_N}$ with arbitrary highest
weights $\la_1,\ldots,\la_N$.

The following result is proved in \cite{bf,FFR,RV}.

\begin{prop}    \label{bethe vectors}
The vector $|w_1^{i_1},\ldots,w_m^{i_m}\ri$ is an eigenvector of the
Gaudin hamiltonians $\Xi_i, i=1,\ldots,\ell$, if and only if the Bethe
Ansatz equations
\begin{equation}    \label{bethe dual}
\sum_{i=1}^N \frac{\langle \la_i,\chal_{i_j} \rangle}{w_j^{i_j}-z_i} -
\sum_{s \neq j} \frac{\langle \al_{i_s},\chal_{i_j}
\rangle}{w_j^{i_j}-w_s^{i_s}} = 0, \qquad j=1,\ldots,m.
\end{equation}
are satisfied.
\end{prop}

Note that the equations \eqref{bethe dual} are nothing but the
equations \eqref{bethe} for the Langlands dual Lie algebra $^L \g$.

One checks also that if this vector is an eigenvector, then it is
automatically a highest weight vector of weight $\ds \sum_{i=1}^N
\la_i - \sum_{j=1}^m \al_{i_j}$. Hence for this vector to be
non-zero, we must have
\begin{equation}    \label{imp cond}
\sum_{i=1}^N \la_i - \sum_{j=1}^m \al_{i_j} = -w_0(\la_\infty)
\end{equation}
for some dominant integral weight $\la_\infty$. Note that this
relation is nothing but a special case of equation \eqref{special
  relation} for $^L\g$ (when $y = 1$).

\subsection{Gaudin model and coinvariants}    \label{coinvariants}

In \cite{FFR} \propref{bethe vectors} is proved using the following
interpretation of the Gaudin hamiltonians.

Let $\G$ be the affine Kac-Moody algebra corresponding to $\g$. It is
the extension of the Lie algebra $\g \otimes \C((t))$ by the
one-dimensional center $\C K$. The commutation relations in $\G$ read
\begin{equation}    \label{comm rel}
[A \otimes f(t),B \otimes g(t)] = [A,B] \otimes fg - \ka_c(A,B)
\on{Res}_{t=0} f dg \cdot K,
\end{equation}
where $\ka_c(\cdot,\cdot)$ is the {\em critical} invariant inner
product on $\g$ defined by the formula $$\ka_c(A,B) = - \frac{1}{2}
\on{Tr}_\g \on{ad} A \on{ad} B.$$

Denote by $\G_+$ the Lie subalgebra $\g \otimes \C[[t]] \oplus \C K$
of $\G$. We extend the action of $\g$ on the finite-dimensional
representation $V_\la$ to $\G_+$ in such a way that $\g \otimes
t\C[[t]]$ acts trivially and $K$ acts as the identity. Denote by
$\V_{\la}$ the Weyl module which is the induced representation of $\G$
$$\V_{\la} = U(\G) \underset{U(\G_+)}\otimes V_{\la}.$$ These are the
representations of {\em critical level}. In the normalization of
\cite{Kac}, the central element $K$ acts as minus the dual Coxeter
number.

Consider the projective line $\pone$ with a global coordinate $t$ and
$N$ distinct finite points $z_1,\ldots,z_N \in \pone$. In the
neighborhood of each point $z_i$ we have the local coordinate $t-z_i$
and in the neighborhood of the point $\infty$ we have the local
coordinate $t^{-1}$. Set $\widetilde{\g}(z_i) = \g \otimes
\C((t-z_i))$ and $\wt\g(\infty) = \g \otimes \C((t^{-1}))$. Let $\G_N$
be the extension of the Lie algebra $\bigoplus_{i=1}^N
\widetilde{\g}(z_i) \oplus \wt\g(\infty)$ by a one-dimensional
center $\C K$ whose restriction to each summand $\widetilde{\g}(z_i)$
or $\wt\g(\infty)$ coincides with the above central extension.
The Lie algebra $\G_N$ naturally acts on the tensor product
$$
\V_{\laa,\la_\infty} = \V_{\la_1} \otimes \ldots \otimes \V_{\la_N}
\otimes \V_{\la_\infty};
$$
in particular, $K$ acts as the identity.

Let $\g_{\zn}=\g_{z_1,\ldots,z_N}$ be the Lie algebra of $\g$--valued
regular functions on $\pone\backslash\{ z_1, \ldots,$ $z_N,\infty \}$
(i.e. rational functions on $\pone$, which may have poles only at the
points $z_1,\ldots,z_N$ and $\infty$). Clearly, such a function can be
expanded into a Laurent power series in the corresponding local
coordinates at each point $z_i$ and at $\infty$. Thus, we obtain an
embedding $$\g_{\zn} \hookrightarrow \bigoplus_{i=1}^N
\widetilde{\g}(z_i) \oplus \wt\g(\infty).$$ It follows from the
residue theorem and formula \eqref{comm rel} that the restriction of
the central extension to the image of this embedding is trivial. Hence
this embedding lifts to an embedding $\g_{\zn} \to \G_N$.

Denote by $H_{\laa,\la_\infty}$ the space of coinvariants of
$\V_{\laa,\la_\infty}$ with respect to the action of the Lie algebra
$\g_{\zn}$. By construction, we have a canonical embedding of the
finite-dimensional representation $V_\la$ into the module $\V_{\la}$:
$$x \in V_\la \arr 1 \otimes x \in \V_\la,$$ which commutes with the
action of $\g$ on both spaces (where $\g$ is embedded into $\G$ as the
constant subalgebra). Thus we have an embedding $V_{\laa,\la_\infty} =
V_{\laa} \otimes V_{\la_\infty}$ into $\V_{\laa,\la_\infty}$. We will
use the same notation $V_{\laa,\la_\infty}$ for the image of this
embedding. Denote by $V^G_{\laa,\la_\infty}$ the subspace of
$G$--invariants (equivalently, $\g$--invariants) in
$V_{\laa,\la_\infty}$ with respect to the diagonal action.

\begin{lem}[\cite{FFR}, Lemma 1]    \label{iso}
The composition of the embedding $V^G_{\laa,\la_\infty} \hookrightarrow
\V_{\laa,\la_\infty}$ and the projection $\V_{\laa,\la_\infty}
\twoheadrightarrow H_{\laa,\la_\infty}$ is an isomorphism.
\end{lem}

Let $\V_0$ be the representation of $\G$, which corresponds to the
one-dimensional trivial $\g$--module $V_0$; it is called the {\em
vacuum module}. Denote by $v_0$ the generating vector of $\V_0$. We
assign the vacuum module to a point $u \in \pone$ which is different
from $z_1,\ldots,z_N,\infty$. Denote by $H_{(\laa,\la_\infty,0)}$ the
space of $\g_{\zn,u}$--invariant functionals on $\V_{\laa,\la_\infty}
\otimes \V_0$ with respect to the Lie algebra
$\g_{\zn,u}$. \lemref{iso} tells us that the composition of the
embedding $V^G_{\laa,\la_\infty} \otimes v_0 \hookrightarrow
\V_{\laa,\la_\infty} \otimes \V_0$ and the projection
$\V_{\laa,\la_\infty} \otimes \V_0 \to H_{(\laa,\la_\infty,0)}$ is an
isomorphism.

Let $v$ be an arbitrary vector in $\V_0$. For any $x \in
V_{\laa,\la_\infty}^G$ consider the vector $x \otimes v \in
\V_{\laa,\la_\infty} \otimes \V_0$. By \lemref{iso}, the projection of
this vector onto $H_{(\laa,\la_\infty,0)}$ is equal to the projection
of a vector of the form $(\Psi_v(u) \cdot x) \otimes v_0$, where
$\Psi_v(u) \cdot x \in V_{\laa,\la_\infty}$. Thus we obtain a
well-defined linear operator $\Psi_v(u)$ on $V_{\laa,\la_\infty}$
corresponding to any $v \in \V_0$ and any point $u \in \pone \bs \{
z_1,\ldots,z_N,\infty \}$.

For $A \in \g$ and $m \in \Z$, denote by $A_m$ the element $A \otimes
t^m \in \G$. Now introduce the following Segal-Sugawara vector in
$\V_0$:
\begin{equation}    \label{sugawara}
S = \frac{1}{2} \sum_{a=1}^d J_{a,-1} J^a_{-1} v_0.
\end{equation}
This vector defines a linear operator $\Psi_S(u)$ on
$V_{\laa,\la_\infty}$.

Denote by $\Delta(\la)$ the scalar by which the Casimir operator
$\Delta$ acts on $V_\la$.

\begin{prop}[\cite{FFR},Prop. 1]    \label{coincide}
We have
$$
\Psi_S(u) = \sum_{i=1}^N \frac{\Xi_i}{u-z_i} + \sum_{i=1}^N
\frac{\Delta(\la_i)}{(u-z_i)^2},
$$
where the $\Xi_i$'s are the Gaudin operators \eqref{hi}.
\end{prop}

Now consider the subspace $\zz(\G)$ of all $\G_+$--invariant vectors
in $\V_0$. One checks that $S \in \zz(\G)$.

\begin{prop}[\cite{FFR}, Prop. 2]
For any $Z_1, Z_2 \in \zz(\G)$ and any points $u_1, u_2$ $\in \pone
\bs \{ z_1,\ldots,z_N,\infty \}$ the linear operators
$\Psi_{Z_1}(u_1)$ and $\Psi_{Z_2}(u_2)$ commute.
\end{prop}

Taking the coefficients in the expansions of the operators of the form
$\Psi_Z(u)$ at $z_1,\ldots,z_N$ we obtain a family of commuting linear
operators on $V_{\laa,\la_\infty}^G$ which includes the Gaudin
hamiltonians. It is natural to call them the {\em generalized Gaudin
hamiltonians}.

\subsection{The center of $\V_0$ and $^L G$--opers}

In order to describe the algebra of generalized Gaudin hamiltonians
and its spectrum we need to recall the description of $\zz(\G)$ from
\cite{FF:gd,F:wak}.

First, observe that each element $v$ of $\zz(\G)$ gives rise to an
endomorphism of $\V_0$ commuting with the action of $\G$ which sends
the generating vector $v_0$ to $v$. Conversely, any $\G$--endomorphism
of $V_0$ is uniquely determined by the image of $v_0$ which
necessarily belongs to $\zz(\G)$. Thus, we obtain an isomorphism
$\zz(\G) \simeq \on{End}_{\G}(\V_0)$ which gives $\zz(\G)$ an algebra
structure. The opposite algebra structure on $\zz(\G)$ coincides with
the algebra structure induced by the identification of $\V_0$ with the
algebra $U(\g \otimes t^{-1} \C[t^{-1}])$.

The realization of $\zz(\G)$ as $\on{End}_{\G}(\V_0)$ allows us to
interpret the action of $\zz(\G)$ on $V_{\laa,\la_\infty}^G$ as
follows. We identify $H_{\laa,\la_\infty,0}$ with
$V_{\laa,\la_\infty}^G$. By functoriality, any endomorphism of $\V_0$
gives rise to an endomorphism of $H_{\laa,\la_\infty,0}$, and hence of
$V_{\laa,\la_\infty}^G$. In particular, we see immediately that the
map $\Psi: \zz(\G) \to \on{End} V_{\laa,\la_\infty}^G$ is an algebra
homomorphism with respect to the algebra structure on $\zz(\G)$ that
we introduced above.

Let $\DerO = \C[[t]] \pa_t$ be the Lie algebra of continuous
derivations of the topological algebra $\OO = \C[[t]]$. The action of
its Lie subalgebra $\on{Der}_0 \OO = t \C[[t]] \pa_t$ on $\OO$
exponentiates to an action of the group $\AutO$ of formal changes of
variables. Both $\DerO$ and $\AutO$ naturally act on $V_0$ in a
compatible way, and these actions preserve $\zz(\G)$. They also act on
the space $\on{Op}_{^L G}(D)$ of $^L G$--opers on the disc $D =
\on{Spec} \C[[t]]$.

Denote by $\on{Fun} \on{Op}_{^L G}(D)$ the algebra of regular
functions on $\on{Op}_{^L G}(D)$. In view of \lemref{free}, it is
isomorphic to the algebra of functions on the space of $\ell$--tuples
$(v_1(t),\ldots,v_\ell(t))$ of formal Taylor series, i.e., the space
$\C[[t]]^\ell$. If we write $v_i(t) = \sum_{n\geq 0} v_{i,n} t^n$,
then we obtain
\begin{equation}    \label{descr of opers}
\on{Fun} \on{Op}_{^L G}(D) \simeq \C[v_{i,n}]_{i \in I, n\geq 0}.
\end{equation}
Note that the vector field $-t \pa_t$ acts naturally on $\on{Op}_{^L
G}(D)$ and defines a $\Z$--grading on $\on{Fun} \on{Op}_{^L G}(D)$
such that $\deg v_{i,n} = d_i+n+1$. The vector field $-\pa_t$ acts as
a derivation such that $-\pa_t \cdot v_{i,n} = -(d_i+n+1) v_{i,n+1}$.

\begin{thm}[\cite{FF:gd,F:wak}]    \label{center}
There is a canonical isomorphism $$\zz(\G) \simeq \on{Fun} \on{Op}_{^L
G}(D)$$ of algebras which is compatible with the action of $\DerO$ and
$\AutO$.
\end{thm}

The module $\V_0$ has a natural $\Z$--grading defined by the formulas
$\deg v_0 = 0, \deg J^a_n = -n$, and it carries a translation operator
$T$ defined by the formulas $T v_0 = 0, [T,J^a_n] = -n
J^a_{n-1}$. \thmref{center} and the isomorphism \eqref{descr of opers}
imply that there exist non-zero vectors $S_i \in \V_0$ of degrees
$d_i+1, i \in I$, such that
$$
\zz(\G) = \C[T^n S_i]_{i \in I,n\geq 0} v_0.
$$
Then under the isomorphism of \thmref{center} we have $S_i \mapsto
v_{i,0}$, the $\Z$--gradings on both algebras get identified and the
action of $T$ on $\zz(\G)$ becomes the action of $-\pa_t$ on $\on{Fun}
\on{Op}_{^L G}(D)$. Note that the vector $S_1$ is nothing but the
vector \eqref{sugawara}, up to a non-zero scalar.

Recall from \cite{FB} that $\V_0$ is a vertex algebra, and $\zz(\G)$
is its commutative vertex subalgebra; in fact, it is the center of
$\V_0$. Consider the corresponding enveloping algebra $U(\zz(\G))$ as
defined in \cite{F:wak}. It is shown in \cite{F:wak} that $U(\zz(\G))$
is isomorphic to the algebra of functions on the space $\on{Op}_{^L
G}(D^\times)$ of $^L G$--opers on the punctured disc. Moreover,
$U(\zz(\G))$ is the center $Z(\G)$ of the completed universal
enveloping algebra of $\G$ at the critical level (see
\cite{BD}). For each integral dominant weight $\la$ we have a
homomorphism $\on{Fun} \on{Op}_{^L G}(D^\times) \simeq Z(\G) \to
\on{End}_{\G} \V_\la$. The following result is proved in \cite{F:wak}.

\begin{thm}    \label{End V la}
The homomorphism $\on{Fun} \on{Op}_{^L G}(D^\times) \to \on{End}_{\G}
\V_\la$ is surjective. Moreover, it identifies $\on{End}_{\G} \V_\la$
with the algebra $\on{Fun} \on{Op}_{^L G}(D)_{\la}$ so that this
homomorphism becomes the natural surjection $\on{Fun} \on{Op}_{^L
G}(D^\times) \to \on{Fun} \on{Op}_{^L G}(D)_{\la}$ induced by the
embedding $\on{Op}_{^L G}(D)_{\la} \hookrightarrow \on{Op}_{^L
G}(D^\times)$.
\end{thm}

In particular, if $\la=0$ we obtain the statement of \thmref{center},
because $\on{End}_{\G}(\V_0) = \zz(\G)$.

\subsection{Eigenvalues of the generalized Gaudin hamiltonians and $^L
  G$--opers}

Now suppose we have an eigenvector $A \in V_{\laa,\la_\infty}^G$ of
the generalized Gaudin hamiltonians. The action of $\zz(\G)$ on it
defines, for any $u \in \pone$, a homomorphism $\zz(\G) \simeq
\on{Fun} \on{Op}_{^L G}(D_u) \to \C$, i.e., a $^L G$--oper on
$D_u$. Let us denote this oper by $\eta_{A,u}$.

The following theorem asserts that these opers on the discs $D_u$ for
different values of $u$ are restrictions of one and the same regular
$^L G$--oper on $\pone \bs \{ z_1,\ldots,z_N,\infty \}$. Moreover,
this oper has regular singularities at $z_1,\ldots,z_N,\infty$ with
residues $-\la_1-\rho,\ldots,\la_N-\rho,-\la_\infty-\rho$ at those
points, and it has trivial monodromy.

\begin{thm}    \label{belongs}
The $^L G$--opers $\eta_{A,u}$ on $D_u$ corresponding to the
eigenvalues of the generalized Gaudin hamiltonians on $A \in
V_{\laa,\la_\infty}^G$ are restrictions to the respective discs of a
unique (regular) $^L G$--oper $\eta_A$ on $\pone \bs \{
z_1,\ldots,z_N,\infty \}$. Moreover, the oper $\eta_A$ belongs to the
space $\on{Op}_{^L G}(\pone)_{(z_i),\infty;\laa,\la_\infty}$. In
particular, it has trivial monodromy representation.
\end{thm}

\begin{proof}
We will give first an abridged version of the proof and then explain
the details.

In \cite{FB} we defined, for any quasi-conformal vertex algebra $V$, a
smooth projective curve $X$, a set of points $x_1,\ldots,x_N \in X$
and a collection of $V$--modules $M_1,\ldots,M_N$, the space of
coinvariants $H_V(X,(x_i),(M_i))$ and its dual space, the space of
conformal blocks $C_V(X,(x_i),(M_i))$. This construction (which is
recalled below) is functorial: if $W \to V$ is a homomorphism of
vertex algebras, then we have natural maps $H_W(X,(x_i),(M_i)) \to
H_V(X,(x_i),(M_i))$ and $C_V(X,(x_i),(M_i)) \to C_W(X,(x_i),(M_i))$.

In the case of the vertex algebra $\V_0$ associated to the affine
Kac-Moody algebra $\G$, the curve $\pone$, the points
$z_1,\ldots,z_N,\infty$ and the modules
$\V_{\la_1},\ldots,\V_{\la_N},\V_{\la_\infty}$, the space of
coinvariants is nothing but the space $H_{\laa,\la_\infty}$, which we
have identified with $V_{\laa,\la_\infty}^G$ in \lemref{iso}.

The subspace $\zz(\G)$ of $\g[[t]]$--invariant vectors in $\V_0$ is a
commutative vertex subalgebra of $\V_0$; in fact, it is the center of
$\V_0$ (see \cite{FB}). The embedding $\zz(\G) \to \V_0$ then gives
rise to a map
$$H_{\zz(\G)}(\pone;(z_i),\infty;(\V_{\la_i}),\V_{\la_\infty}) \to
H_{\laa,\la_\infty}.$$ Applying the results of \cite{FB}, Sect. 8.4,
(as explained below) we obtain that each eigenvector $A$ of the
generalized Gaudin hamiltonians in $V_{\laa,\la_\infty}^G \simeq
H_{\laa,\la_\infty}$ gives rise to a character (i.e., an algebra
homomorphism) $$\on{Fun} \on{Op}_{^L G}(\pone \bs \{
z_1,\ldots,z_N,\infty \}) \to \C,$$ i.e., to a $^L G$--oper on $\pone
\bs \{ z_1,\ldots,z_N,\infty \}$. This is precisely the $^L G$--oper
$\eta_A$ that we are looking for. Next, we apply \thmref{End V la} to
show that $\eta_A$ actually belongs to the space $\on{Op}_{^L
G}(\pone)_{(z_i),\infty;\laa,\la_\infty}$.

Let us now explain all of this in detail. First, we recall the
definition of the space of conformal blocks from \cite{FB}, Ch. 8. For
that we will need a coordinate-independent description of the
structure of a module over a vertex algebra given in \cite{FB}, Ch. 5
and Sect. 6.3.9, where we refer the reader for more details.

Let $X$ be a smooth algebraic curve and ${\mc A}ut_X$ be the principal
$\AutO$--bundle over $X$ whose fiber ${\mc A}ut_x$ at $x \in X$ is the
space of formal coordinates at $x$. Let $V$ be a quasi-conformal
vertex algebra (see \cite{FB}, Sect. 5.2.4). It then carries an action
of $\AutO$. We define a vector bundle ${\mc V} = {\mc V}_{X}$ on
$\pone$ as the twist ${\mc A}ut_X \underset{\AutO}\times V$. This
bundle carries a (flat) connection. If we choose a coordinate $t$ and
trivialize ${\mc A}ut_X$ and ${\mc V}$ using this coordinate, then the
connection operator reads $\nabla = \pa_t + T$.

Let $M$ be a $V$--module which carries an action of $\on{Der}_0 \OO =
t\C[[t]] \pa_t$ compatible with that of $V$ such that the action of
$-t\pa_t$ is semi-simple and the eigenvalues belong to the union of
the sets $\ka_i + \Z_+$, where $\{ \ka_i \}$ is a finite set of
complex numbers. The action of the Lie algebra $\on{Der}_+ \OO = t^2
\C[[t]] \pa_t$ on $M$ may be exponentiated to an action of the group
$\on{Aut}_+ \OO$ consisting of the formal coordinate changes of the
form $z \mapsto z + z^2(\ldots)$. Let us fix a non-zero tangent vector
$\tau$ at $x$ and consider the $\on{Aut}_+ \OO$--torsor ${\mc
A}_{x,\tau}$ consisting of all formal coordinates at $x$ whose one-jet
is equal to $\tau$. We define the twist ${\mc M}_x = {\mc
A}ut_{x,\tau} \underset{\on{Aut}_+ \OO}\times M$ of $M$ at $x \in X$.

Let us pick a formal coordinate $t_x$ at $x$ whose one-jet is equal to
$\tau$. We use this coordinate to trivialize ${\mc V}|_{D_x}$ and
${\mc M}_x$ and to define an $\on{End} {\mc M}_x$--valued section
${\mc Y}_x^M$ of ${\mc V}^*|_{D_x^\times}$ as follows. The value
$\langle \varphi,{\mc Y}_x^M \cdot v \rangle$ of this section on $v
\in {\mc M}_x \simeq M$, $\varphi \in {\mc M}_x^* \simeq M^*$ and the
constant section $s_A$ of ${\mc V}|_{D_x}$ corresponding to the vector
$A \in V$ with respect to our trivialization, is equal to $\langle
\varphi,Y^M(A,t_x) v \rangle$. It is proved in \cite{FB} that the
section ${\mc Y}_x^M$ is well-defined, i.e., independent of the choice
of the coordinate $t_x$. Moreover, this section is horizontal with
respect to the connection on $\V^*$ which is the transpose of the
connection $\nabla$ on $\V$ (see \cite{FB}, Theorem 5.5.3).

Let $x_1,\ldots,x_N$ be a collection of distinct points on $X$. We
will fix once and for all a non-zero tangent vector $\tau_i$ at $x_i$
for each $i=1,\ldots,N$. The space of conformal blocks
$C_V(X,(x_i),(M_i))$ is by definition the space of linear functionals
$\varphi$ on ${\mc M}_{1,x_1} \otimes \ldots \otimes {\mc M}_{N,x_N}$
satisfying the following condition: for any $A_i \in {\mc M}_{i,x_i},
i=1,\ldots,N$, there exists a regular section of $\V^*$ on $X \bs \{
x_1,\ldots,x_N \}$ such that for all $i=1,\ldots,N$ its restriction to
$D^\times_{x_i}$ is equal to
$$
\langle \varphi,A_1 \otimes \ldots \otimes {\mc Y}^{M_i}_{x_i} \cdot
A_i \otimes \ldots \otimes A_N \rangle.
$$
This regular section is then automatically horizontal. This section
may be constructed explicitly as follows. In \cite{FB}, Theorem 9.3.1,
we established, for all $y \in X, u \neq x_i$, an isomorphism
$$C_V(X,(x_i),(M_i)) \simeq C_V(X;(x_i),u;(M_i),V).$$ This
means that the space of conformal blocks does not change if we insert
the vacuum module $V$ at a point $u \in X$ different from all the
$x_i$'s; note that we have considered in \secref{coinvariants} a
special case of this isomorphism. Let $\wt\varphi$ be the functional
in $C_V(X;(x_i),y;(M_i),V)$ corresponding to $\varphi \in
C_V(X,(x_i),(M_i))$ under this isomorphism. Then the value of our
section at $y \in X$ on an element $A \in \V_y$ is precisely equal to
$\langle \wt\varphi,A_1 \otimes \ldots \otimes A_N \otimes A \rangle$.

Now we set $V$ to be the affine Kac-Moody vertex algebra $\V_0$,
$X=\pone$, with the marked points $z_1,\ldots,z_N,\infty$, and take as
the modules attached to these points the Weyl modules
$\V_{\la_1},\ldots,\V_{\la_N},\V_{\la_\infty}$. Our global coordinate
$t$ on $\pone$ gives rise to the coordinate $t-z_i$ at each point
$z_i$ and the coordinate $t^{-1}$ at $\infty$. Hence we obtain an
identification of ${\mc V}_{\la_i,z_i}$ with $\V_{\la_i}$. It is
proved in \cite{FB} (see Theorem 8.3.3 and Remark 8.3.10) that the
corresponding space of conformal blocks is the space of
$\g_{\zn}$--invariant functionals on $\V_{\laa,\la_\infty}$, i.e., the
dual space to $H_{\laa,\la_\infty}$.

Let $\varphi$ be a linear functional on $\V_{\laa,\la_\infty}$. For
each vector $A_1 \otimes \ldots \otimes A_N \otimes A_\infty \in
\V_{\laa,\la_\infty}$ we then obtain a section
\begin{equation}    \label{section}
\langle \varphi,A_1 \otimes \ldots \otimes {\mc Y}_{z_i}^{\V_{\la_i}}
\cdot A_i \otimes \ldots \otimes A_N \otimes A_\infty \rangle
\end{equation}
of ${\mc V}_{0,\pone}|_{D^\times_{z_i}}$ for all $i=1,\ldots,N$, and
likewise at the point $\infty$. According to the above discussion, the
functional $\varphi$ is $\g_{\zn}$--invariant if and only if the above
sections are restrictions to the respective punctured discs of a
single rational section of ${\mc V}_{0,\pone}^*$ with poles only at
the points $z_1,\ldots,z_N$ and $\infty$, which is horizontal with
respect to the connection $\nabla$. Moreover, the value of this
section at $u \in \pone \bs \{ z_1,\ldots,z_N,\infty \}$ may be
obtained as explained above, by inserting $\V_0$ at $u$.

For any eigenvector $A$ of the Gaudin hamiltonians in
$V_{\laa,\la_\infty}^G \simeq H_{\laa,\la_\infty}$, there is a linear
functional on $H_{\laa,\la_\infty}$, taking a non-zero value on $A$,
which is an eigenvector of the transposed Gaudin operators and has the
same eigenvalues. We view this functional as a conformal block. Then
it satisfies the above condition, namely, that the sections
\eqref{section} are restrictions to the respective punctured discs of
a single horizontal section of ${\mc V}_{0,\pone}^*$ that is regular
on $\pone \bs \{ z_1,\ldots,z_N,\infty \}$. We will denote this
section by $s_{\varphi}$.  Evaluating $s_{\varphi}$ on an arbitrary
section of ${\mc V}_{0,\pone}$, we obtain a rational function on
$\pone$ with poles at $z_1,\ldots,z_N,\infty$.

Consider the subbundle ${\mc Z}_{\pone}$ of ${\mc V}_{0,\pone}$
obtained by twisting $\zz(\G) \subset \V_0$. The fiber ${\mc Z}_u$ of
${\mc Z}_{\pone}$ at $u \in \pone$ is just the algebra of functions on
the space $\on{Op}_{^L G}(D_u)$ of $^L G$--opers on $D_u$. Therefore
${\mc Z}_{\pone}$ is nothing but the algebra of functions on the
scheme ${\mc O}p_{^L G}(X)$ of jets of $^L G$--opers on $\pone$, whose
fiber at $u \in \pone$ is the space $\on{Op}_{^L G}(D_u)$. This scheme
carries a natural connection and the corresponding connection on ${\mc
Z}_{\pone}$ coincides with the connection $\nabla$ described
above. Note that horizontal sections of ${\mc O}p_{^L G}(\pone)$ over
$U \subset \pone$ are the same as the regular $^L G$--opers on $U$.

We now evaluate $s_{\varphi}$ on sections of ${\mc
Z}_{\pone}$. According to our construction of \secref{coinvariants},
the value of $s_\varphi$ on $v \in {\mc Z}_u$ at the point $u \in
\pone \bs \{ z_1,\ldots,z_N,\infty \}$ is precisely the eigenvalue of
the generalized Gaudin hamiltonian $\Psi_v(u)$ on our
eigenvector. Moreover, these eigenvalues are multiplicative with
respect to the commutative algebra structure on the bundle ${\mc
Z}_{\pone}$, which is inherited from that on $\zz(\G)$. Therefore
these eigenvalues define an algebra homomorphism ${\mc Z}_u \to \C$
for all $u \in \pone \bs \{ z_1,\ldots,z_N,\infty \}$. This is the
same as an algebra homomorphism from the sheaf of algebras
$$
{\mc Z}_{\pone \bs  \{ z_1,\ldots,z_N,\infty \}} \simeq \on{Fun} {\mc
  O}p_{^L   G}(\pone \bs \{ z_1,\ldots,z_N,\infty \})
$$
to $\C$ (considered as the constant sheaf over $\pone \bs \{
z_1,\ldots,z_N,\infty \}$). Moreover, according to the general results
on conformal blocks (see above), this homomorphism must be
horizontal. But such a homomorphism is the same as a horizontal
section of the bundle ${\mc O}p_{^L G}(\pone)$ of jets of $^L
G$--opers on $\pone \bs \{ z_1,\ldots,z_N,\infty \}$, which is the
same as a regular $^L G$--oper on $\pone \bs \{
z_1,\ldots,z_N,\infty \}$. This is the desired oper $\eta_A$. By
construction, its restriction to $D_u$ for each $u \in \pone \bs \{
z_1,\ldots,z_N,\infty \}$ equals to the $^L G$--oper on $D_u$ which
records the eigenvalues of the generalized Gaudin hamiltonians
corresponding to the point $u$.

Let us now look at the restrictions of the oper $\eta_A$ to the
punctured discs around the points $z_1,\ldots,z_N$ and $\infty$. By
definition of the section $s_\varphi$, these restrictions are equal to
the sections \eqref{section}, where the vertex operations ${\mc
Y}^{V_{\la_i}}_{z_i}$ and ${\mc Y}^{V_{\la_\infty}}_\infty$ are
restricted to $\zz(\G)_{z_i}$ and $\zz(\G)_{\infty}$ (which are the
twists of $\zz(\G)$ by ${\mc A}ut_{z_i}$ and ${\mc A}ut_\infty$,
respectively). Each section gives rise to a homomorphism $Z(\G)_{z_i}
\to \C$, and hence to a point in $\on{Spec} Z(\G)_{z_i} = \on{Op}_{^L
G}(D_{z_i}^\times)$. But we know from \thmref{End V la} that the
action of the center $Z(\G)$ on $\V_\la$ factors through the algebra
of functions on $\on{Op}_{^L G}(D)_{\la}$. Therefore this point
belongs to $\on{Op}_{^L G}(D_{z_i})_{\la_i} \subset \on{Op}_{^L
G}(D_{z_i}^\times)$. Hence we find that the restriction of our $^L
G$--oper on $\pone \bs \{ z_1,\ldots,z_N,\infty \}$ to the disc
$D_{z_i}^\times$ (resp., $D_\infty^\times$) belongs to $\on{Op}_{^L
G}(D_{z_i})_{\la_i}$ (resp., $\on{Op}_{^L
G}(D_\infty)_{\la_\infty}$). Therefore this oper belongs to
$\on{Op}_{^L G}(\pone)_{(z_i),\infty;\laa,\la_\infty}$, which is what
we wanted to prove.
\end{proof}

In more concrete terms, the oper $\eta_A$ may be described as
follows. From the description of $\zz(\G)$ we know that all
eigenvalues are encoded in the rational functions $v^A_i(u)$ which are
the eigenvalues of the operators $\Psi_{S_i}(u), i=1,\ldots,\ell$, on
$A$. The corresponding $^L G$--oper connection then reads (with
respect to our trivialization of $\F$ and the global coordinate $t$ on
$\pone$)
\begin{equation}    \label{conn eigenvalue}
\nabla = \pa_t + p_{-1} + \sum_{i \in I} v^A_i(t) p_i.
\end{equation}

\subsection{Completeness of the Bethe Ansatz}

According to \thmref{belongs}, each point in the spectrum of the
generalized Gaudin hamiltonians occurring in $V_{\laa,\la_\infty}^G$
(i.e., a collection of joint eigenvalues of these operators on
$V_{\laa,\la_\infty}^G$) is encoded by a $^L G$--oper on $\pone$ with
regular singularities at $z_1,\ldots,z_N,\infty$ which has trivial
monodromy. Moreover, two different points of the spectrum give rise to
different opers. Thus, we obtain the following

\begin{cor}    \label{inj map}
There is an injective map from the spectrum of the generalized Gaudin
hamiltonians on $V^G_{\laa,\la_\infty}$ (not counting multiplicities)
to the set $\on{Op}_{^L G}(\pone)_{(z_i),\infty;\laa,\la_\infty}$ of
$^L G$--opers on $\pone$ with regular singularities at
$z_1,\ldots,z_N,\infty$ which have trivial monodromy.
\end{cor}

On the other hand, suppose that we are given a $^L G$--oper $\tau$ in
$\on{Op}_{^L G}(\pone)_{(z_i),\infty;\laa,\la_\infty}$. Consider the
(unique) Miura oper structure on it for which the horizontal Borel
reduction coincides with the oper reduction at the point
$\infty$. Suppose that this Miura oper satisfies the conditions (1)
and (2) from \secref{miura opers on pone}, i.e., it belongs to the
space
$\on{MOp}_{G}(\pone)^{\on{gen}}_{(z_i),\infty;(\cla_i),\cla_\infty}$.
Then we will call $\tau$ a {\em non-degenerate} oper.

According to \thmref{main}, there is a bijection between the points of
the space
$\on{MOp}_{G}(\pone)^{\on{gen}}_{(z_i),\infty;(\cla_i),\cla_\infty}$
and the set of solutions of the Bethe Ansatz equations. Note that we
have switched to the Langlands dual group $^L G$, and so these
equations are given by formula \eqref{bethe dual}. Our Miura $^L
G$--oper, for which the horizontal Borel reduction coincides with the
oper reduction at the point $\infty$, gives rise to a unique solution
of equations \eqref{bethe dual} which satisfies the condition
\eqref{imp cond}. We will refer to it as the {\em special solution}
corresponding to the non-degenerate $^L G$--oper $\tau$ from
$\on{Op}_{^L G}(\pone)_{(z_i),\infty;\laa,\la_\infty}$.

According to \propref{bethe vectors}, we associate to this special
solution of the Bethe Ansatz equations an eigenvector of the Gaudin
hamiltonians by formula \eqref{genbv}. Let us denote this eigenvector
by $v_\tau$.

A natural question is what are the eigenvalues of the generalized
Gaudin hamiltonians on $v_\tau$. By \thmref{belongs}, these
eigenvalues are encoded by a $^L G$--oper in $\on{Op}_{^L
G}(\pone)_{(z_i),\infty;\laa,\la_\infty}$. Not surprisingly, the
answer is that this oper is $\tau$ itself (see \cite{FFR}, Theorem 3):

\begin{prop}
The eigenvalues of the generalized Gaudin hamiltonians acting on the Bethe
eigenvector $v_\tau$ constructed from the special solution of the
Bethe Ansatz equations corresponding to $\tau \in \on{Op}_{^L
G}(\pone)_{(z_i),\infty;\laa,\la_\infty}$ are encoded precisely by the
$^L G$--oper $\tau$.
\end{prop}

Thus, the $^L G$--oper on $\pone$ corresponding to the eigenvalues of
the Gaudin hamiltonians on a given Bethe vector \eqref{genbv} may be
found by applying the Miura transformation (see \secref{miura trans})
to the $^L H$--connection
$$
\pa_t - \sum_{i=1}^N \frac{\la_i}{t-z_i} + \sum_{j=1}^m
\frac{\al_{i_j}}{t-w_j},
$$
where $w_1,\ldots,w_m$ satisfy the Bethe Ansatz equations \eqref{bethe
dual} and the condition \eqref{imp cond}. For example, in the case of
$\sw_n$, the $PGL_n$--oper is nothing but an $n$th order differential
operator \eqref{sln-oper}, and the Miura transformation has the form
\eqref{miura trans for sln}. So we need to write
$$
\pa_t - \sum_{i=1}^N \frac{\la_i}{t-z_i} + \sum_{j=1}^m
\frac{\al_{i_j}}{t-w_j} = \pa_t + \sum_{k=1}^n u_k(t) \epsilon_k,
$$
where we identify the dual Cartan subalgebra of $\sw_n$ with the
hyperplane $\sum_{k=1}^n \epsilon_k = 0$ of the vector space
$\on{span} \{ \epsilon_k \}_{k=1,\ldots,n}$. Then the corresponding
$PGL_n$--oper is given by formula \eqref{miura trans for sln}. One
obtains similarly the opers for other simple Lie algebras of classical
types.

Let us assume from now on that all $^L G$--opers in $\on{Op}_{^L
G}(\pone)_{(z_i),\infty;\laa,\la_\infty}$ are non-degenerate.\footnote{It
follows from the results of Mukhin and Varchenko in \cite{MV:new}
that for some $\la_1,\ldots,\la_N,\la_\infty$ this may not be the case
even for generic values of $z_1,\ldots,z_N$.} Then to each $^L
G$--oper $\tau$ in $\on{Op}_{^L
G}(\pone)_{(z_i),\infty;\laa,\la_\infty}$ corresponds a special
solution of the Bethe Ansatz equations and hence a Bethe vector
$v_\tau$. If all Bethe vectors $v_\tau$ are non-zero, then we obtain
an inverse map to the map of \corref{inj map}, which assigns to $\tau
\in \on{Op}_{^L G}(\pone)_{(z_i),\infty;\laa,\la_\infty}$ the point in
the spectrum corresponding to the eigenvector $v_\tau$ (note that a
priori it could happen that there are other eigenvectors with the same
eigenvalues $\tau$). This leads us to the following result.

\begin{prop}    \label{completeness}
Suppose that all $^L G$--opers in $\on{Op}_{^L
G}(\pone)_{(z_i),\infty;\laa,\la_\infty}$ are non-dege\-ne\-rate and
that all Bethe vectors obtained from solutions of the Bethe Ansatz
equations \eqref{bethe dual} satisfying the condition \eqref{imp cond}
are non-zero. Then there is a bijection between the spectrum of the
generalized Gaudin hamiltonians on $V^G_{\laa,\la_\infty}$ (not
counting multiplicities) and the set $\on{Op}_{^L
G}(\pone)_{(z_i),\infty;\laa,\la_\infty}$ of $^L G$--opers on $\pone$
with regular singularities at $z_1,\ldots,z_N,\infty$ which have
trivial monodromy.

Moreover, if in addition the Gaudin hamiltonians are diagonalizable
and have simple spectrum on $V^G_{\laa,\la_\infty}$, then the Bethe
vectors constitute an eigenbasis of $V^G_{\laa,\la_\infty}$.
\end{prop}

The last statement of this proposition that the Bethe vectors
constitute an eigenbasis of $V^G_{\laa,\la_\infty}$ is referred to as
the completeness of the Bethe Ansatz (sometimes completeness is taken
to mean that the Bethe vectors span $V^G_{\laa,\la_\infty}$, but we
use this term to mean that they form an eigenbasis).

For $\g=\sw_2$ and generic values of the $z_i$'s it was proved by
Scherbak and Varchenko in \cite{SV} (see also \cite{RV}) that the
Bethe vectors are all non-zero. It also follows from \cite{SV} that in
the case of $\sw_2$ all opers are non-degenerate when $z_1,\ldots,z_N$
are in generic position. Hence we obtain a bijection between the
spectrum of the Gaudin hamiltonians on $V^{SL_2}_{\laa,\la_\infty}$
and the set $\on{Op}_{PGL_2}(\pone)_{(z_i),\infty;\laa,\la_\infty}$
for generic $z_i$'s. Moreover, Scherbak \cite{S} has shown that the
eigenvalues of the Gaudin hamiltonians have no multiplicities on the
Bethe vectors for generic $z_i$'s, so we obtain the completeness of
the Bethe Ansatz as well.

We note that the completeness of the Bethe Ansatz has been previously
proved for $\g=\sw_2$ and generic values of $z_1,\ldots,z_N$ by
Varchenko and Scherbak \cite{SV} by other methods. In addition, it
follows from the results of Mukhin and Varchenko \cite{MV} and
Scherbak \cite{S2} that for $\g=\sw_n$ the number of points of
$\on{Op}_{PGL_n}(\pone)_{(z_i),\infty;\laa,\la_\infty}$ is less than
or equal to the dimension of $V_{\laa,\la_\infty}^{SL_n}$.

In the general case we have the following conjecture.

\begin{conj}    \label{conj completeness}
For generic values of $z_1,\ldots,z_N$ the generalized Gaudin
hamiltonians are diagonalizable on $V_{\laa,\la_\infty}^G$ and have
simple spectrum, and the Bethe vectors corresponding to the solutions
of the Bethe Ansatz equations \eqref{bethe dual} are all non-zero.
\end{conj}

If the statement of \conjref{conj completeness} is true and all $^L
G$--opers in $\on{Op}_{^L G}(\pone)_{(z_i),\infty;\laa,\la_\infty}$
are non-degenerate, then we obtain from \propref{completeness} the
completeness of the Bethe Ansatz and a bijection between the spectrum
of the Gaudin hamiltonians, counted with multiplicity, and the set
$\on{Op}_{^L G}(\pone)_{(z_i),\infty;\laa,\la_\infty}$.

Suppose now that there are degenerate opers in $\on{Op}_{^L
G}(\pone)_{(z_i),\infty;\laa,\la_\infty}$. Consider again the unique
Miura oper structure on one of the degenerate opers for which the
horizontal Borel reduction coincides with the oper reduction at the
point $\infty$. According to \thmref{strongest}, this Miura oper
corresponds to a connection in
$\on{Conn}(\pone)_{(z_i),\infty;(\la_i),\la_\infty}$ which has the
form \eqref{conn RS}:
\begin{equation}    \label{new conn RS}
\pa_t - \sum_{i=1}^N \frac{y_i(\la_i+\rho)-\rho}{t-z_i}
- \sum_{j=1}^m \frac{y'_j(\rho)-\rho}{t-w_j}
\end{equation}
for some $y_i, y'_j \in W$ satisfying the relation \eqref{relation}
with $y_\infty = w_0$:
$$
\sum_{i=1}^N (y_i(\la_i+\rho)-\rho) + \sum_{j=1}^m
(y'_j(\rho)-\rho) = -w_0(\la_\infty).
$$
The fact that $\tau$ is not generic means that either some of the
elements $y_i$ are not equal to $1$ or some of the elements $y'_j$
have lengths greater than $1$ (i.e., are not simple reflections). We
expect that if $z_1,\ldots,z_N$ are generic, then for every $\tau \in
\on{Op}_{^L G}(\pone)_{(z_i),\infty;\laa,\la_\infty}$ we have $y_i =
1$ for all $i=1,\ldots,N$ in formula \eqref{new conn RS}. In other
words, we expect that for generic $z_1,\ldots,z_N$ this Miura oper
still satisfies condition (1) from \secref{miura opers on pone}, but
may not satisfy condition (2), that is at least one of the $y'_j$'s is
not a simple reflection.

Then we can still attach to the connection \eqref{new conn RS} an
eigenvector of the generalized Gaudin hamiltonians in
$V_{\laa,\la_\infty}^G$ by generalizing the procedure of
\cite{FFR}. We expect that for generic $z_1,\ldots,z_N$ all of these
vectors are non-zero and that they provide an eigenbasis for the
generalized Gaudin hamiltonians in $V_{\laa,\la_\infty}^G$. For more
on this, see Sect. 5.5 of \cite{F:faro}.

\section{Opers and Bethe Ansatz equations for arbitrary Kac-Moody
  algebras}    \label{fourth}

In this section we generalize some of the results of the previous
sections to the situation where $\g$ is an arbitrary Kac-Moody
algebra. One can easily write down the Bethe Ansatz equations in this
general setting and try to describe the set of solutions of these
equations. We show that, just as in the case of a simple
finite-dimensional Lie algebra, this set is an open subset of the
(ind-)flag variety of $\g$. For that we introduce the notions of opers
and Miura opers for an arbitrary Kac-Moody algebra and show that the
set of solutions of the Bethe Ansatz equations is an open and dense
subset in the set of Miura opers on the projective line with
prescribed residues at marked points (as in the finite-dimensional
case).

\subsection{Opers and Miura opers for general Kac-Moody algebras}

Let $\g$ be the Kac-Moody algebra associated to a Cartan matrix $A$ of
size $\ell \times \ell$ (not necessarily symmetrizable) and $\h$ be
its (extended) Cartan subalgebra of dimension $\ell + d$, where
$\ell-d$ is the rank of $A$. We use the same notation as before for
coroots and roots of $\g$, which are vectors in $\h$ and $\h^*$,
respectively. We have the Cartan decomposition $\g = \n_+ \oplus \h
\oplus \n_-$, and the generators $\{ e_i \}_{i=1,\ldots,\ell}$ and $\{
f_i \}_{i=1,\ldots,\ell}$ of $\n_+$ and $\n_-$, respectively. The Lie
subalgebra $\n_-$ has a natural descending filtration by Lie ideals of
finite codimension. We consider its completion with respect to this
filtration and the corresponding completion of $\g$. From now on we
will use the symbols $\n_-$ and $\g$ to denote these completions.

For example, in the case of untwisted affine algebras, the completed
Lie algebra $\g$ has the form $\ol{\g}((t^{-1})) \oplus \C K \oplus \C
d$, where $\ol{\g}$ is a finite-dimensional simple Lie algebra, $K$
is the central element and $d$ is the vector field $t \pa_t$.

Let $\wt{G}$ be the algebraic group associated to $\g$ in
\cite{Kash}. If $\g$ is infinite-dimensional, then $\wt{G}$ is not a
group scheme, but a group ind-scheme. We denote by $G$ the quotient of
$\wt{G}$ by its center (which belongs to the Cartan subgroup of $G$
corresponding to the Lie subalgebra $\h$). It comes with the lower
unipotent and Borel subgroups $N_-$ and $B_-$ (which are proalgebraic
groups) corresponding to $\n_-$ and $\bb_- = \h/{\mathfrak c} \oplus
\n_-$, respectively, and the upper unipotent and Borel subgroups $N_+$
and $B_+$ (which are group ind-schemes) corresponding to $\n_+$ and
$\bb_+ = \h/{\mathfrak c} \oplus \n_+$, respectively (here ${\mathfrak
c}$ is spanned by those elements $x$ of $\h$ that satisfy $\langle
\al_i,x \rangle = 0, i=1,\ldots,\ell$). We denote by $H$ the intersection
$B_+ \cap B_-$. It is isomorphic to $B_+/N_+$ and to $B_-/N_-$.

We wish to define the spaces of $G$--opers
and Miura $G$--opers on $X$ (which is again a smooth curve or a disc
or a punctured disc).

First we need to introduce the notion of a $G$--bundle on $X$ and a
connection on such a bundle. A $G$--bundle on $X$ is an ind-scheme
$\F$ over $X$ equipped with fiberwise simply transitive action of $G$,
which is locally trivial in the Zariski topology. This means that $X$
may be covered by Zariski open subsets $U_i$ such that the restriction
of $\F$ to each $U_i$ is isomorphic to the trivial bundle $U_i \times
G$. Two such trivializations differ by a morphism $U_i \to G$ called
the change of trivializations. If $U = \on{Spec} R$, then the changes
of trivializations on $U$ form the group $G(R)$.

To define a connection on a $G$--bundle it suffices to define the
notion of a connection on the trivial $G$--bundle on an affine curve
$X = \on{Spec} R$ and explain how to act on these connections by the
changes of trivialization. Without loss of generality we may assume
that we are given an \'etale coordinate $t: X \to {\mathbb A}^1$ on
$X$. Let $\pa_t$ be the vector field on $X$ induced by a fixed
translation vector field on ${\mathbb A}^1$. Then a connection on the
trivial bundle is by definition an operator $\nabla = \pa_t + A(t)$,
where $A(t) \in \g(R)$. If $g \in G(R)$ is a change of trivialization,
then it acts on $\nabla$ by the usual formula
$$
\nabla \mapsto \pa_t + g A(t) g^{-1} - (\pa_t g) g^{-1}.
$$
Under a change of coordinates $t = \varphi(s)$ the operator $\nabla$
transforms in the usual way:
$$
\nabla \mapsto \pa_s + \varphi'(s) A(\varphi(s)).
$$

It is easy to render this definition into the setting of analytic
topology.

We will say that a connection $\nabla$ on $\F$ gives rise to a
trivialization of $\F$ over $X$ if there is a trivialization of $\F$
over $X$ with respect to which $\nabla = \pa_t$.

Note that for ind-groups, such as $G$, a connection does not
necessarily give rise to a trivialization of $\F$, even locally
analytically. The usual correspondence between connections on
$G$--bundles and local trivializations of the $G$--bundles does not
exist in this case, because we do not have the exponential map from
the Lie algebra $\g$ to the group $G$ (though this correspondence
exists if $X$ is a formal disc). However, in what follows we will
consider Miura opers which carry a reduction to the {\em proalgebraic}
group $B_-$ preserved by the connection. In this case a connection
does give rise to local analytic trivializations of the underlying
$B_-$--bundle, and hence the $G$--bundle as well.

The definition of $G$--opers is similar to the definition given in
\secref{opers} in the finite-dimensional case.

A $G$--oper on $X$ is a triple $(\F,\nabla,\F_{B_+})$, where $\F$ is a
principal $G$--bundle $\F$ on $X$, $\nabla$ is a connection on $\F$
and $\F_{B_+}$ is a $B_+$--reduction of $\F$, such that locally, with
a choice of a coordinate $t$ and a trivialization of $\F_{B_+}$, the
connection operator has the form
\begin{equation}    \label{oper new}
\nabla = \pa_t + \sum_{i=1}^\ell \psi_i(t) f_i + {\mb v}(t),
\end{equation}
where each $\psi_i(t)$ is a nowhere vanishing function, and ${\mb
v}(t)$ is a $\bb_+$--valued function. We denote the set of
$G$--opers on $X$ by $\on{Op}_G(X)$.

The changes of trivialization amount in this case to the gauge action
by $B_+$, so when $X = \on{Spec} R$ is affine, a $G$--oper is a gauge
equivalence class of operators \eqref{oper new}, where ${\mb v}(t)
\in \bb_+(R)$, with respect to the group of gauge transformations by
$B_+(R)$. This is the same as an $N_+(R)$ gauge equivalence class of
operators of the form
\begin{equation}    \label{oper new1}
\nabla = \pa_t + p_{-1} + {\mb v}(t), \qquad {\mb v}(t) \in \bb_+(R),
\end{equation}
where, as before, $p_{-1} = \sum_{i=1}^\ell f_i$.

Next, we give the definition of Miura $G$--opers. A Miura $G$--oper on
$X$ is a quadruple $(\F,\nabla,\F_{B_+},\F_{B_-})$, where
$(\F,\nabla,\F_{B_+})$ is a $G$--oper on $X$ and $\F_{B_-}$ is a
$B_-$--reduction of $\F$ which is preserved by $\nabla$. We denote the
set of Miura $G$--opers on $X$ by $\on{MOp}_G(X)$.

In the case when $\g$ is finite-dimensional, this definition is
equivalent to our old definition from \secref{miura opers}. Indeed,
$B_-$ is conjugate to $B_+$. Therefore a $B_-$--reduction $\F_{B_-}$
gives rise to a $B_+$--reduction $\F_{B_-} w_0$, which is preserved by
the connection. We may then take this $B_+$--reduction as the
reduction $\F'_{B_+}$ of our old definition. But in the
infinite-dimensional case the groups $B_+$ and $B_-$ are not conjugate
to each other (in fact, one of them is not even a group scheme but a
group ind-scheme) and there is an essential difference between asking
for a horizontal $B_+$--reduction or a horizontal
$B_-$--reduction.

In fact, for the purposes of the present paper it is essential that
the horizontal reduction be to a proalgebraic subgroup $B_-$ and
the oper reduction be to an ind-subgroup $B_+$. Indeed, we wish to
relate our Miura opers to Cartan connections (see \propref{map beta1}
below). We will do this by intersecting the two reductions inside
$\F$, so they need to be ``opposite'' to each other. Next, since the
connection operators of the Miura opers preserve a $B_-$--bundle
$\F_{B_-}$, and $B_-$ is a proalgebraic group (not an ind-group), it
makes sense to talk about parallel transport and horizontal sections
on $\F_{B_-}$ (and hence on the induced $G$--bundle $\F$) over an
arbitrary curve. Hence we can trivialize locally a $G$--bundle
equipped with a connection and a horizontal $B_-$--reduction. Then a
reduction to $B_+$ gives rise to locally defined maps to $G/B_+$ which
is a scheme of infinite type.

If we were to switch $B_+$ and $B_-$, our horizontal reduction would
be to an ind-group $B_+$, and the notion of parallel transport would
only make sense over a formal disc. But in what follows we need to use
this notion for arbitrary curves (particularly, for $\pone$), and this
forces us to define opers and Miura opers in this fashion.

\begin{remark}    \label{different}
In \cite{FB1} Ben-Zvi and the author have already defined ``affine
opers'' and ``affine Miura opers''. However, these objects are
different from the $G$--opers and Miura $G$--opers for an untwisted
affine Kac-Moody algebra $\g$ that we consider here, because in
\cite{FB1} we considered the completion of $\n_+$ rather than $\n_-$,
so that $\g = \ol{\g}((t)) \oplus \C K \oplus \C d$ (in the case of
opers, we had also chosen in addition to the above data a reduction to
the subgroup $\ol{G}[t^{-1}]$ of $G$). In other words, in \cite{FB1}
the roles of $B_+$ and $B_-$ were switched in the sense that in
\cite{FB1} the group $B_+$ was a proalgebraic group and $B_-$ was an
ind-group.\qed
\end{remark}

Next, we define $G$--opers on the disc $D_x$ with regular singularity
at $x$ following \secref{reg sing}: these are the
$N_+((t))$--equivalence classes of operators of the form
\begin{equation}    \label{oper with RS2}
\nabla = \pa_t + \frac{1}{t} \left( p_{-1} + {\mb v}(t) \right),
\qquad {\mb v}(t) \in \bb_+[[t]].
\end{equation}
Denote by $\on{Op}_G^{\on{RS}}(D_x)$ the space of opers on $D_x$ with
regular singularity. By definition, it is a subspace of
$\on{Op}_G(D_x^\times)$.

Finally, we define, for any dominant integral coweight $\cla \in
\h/\cc$, the notion of a $G$--oper of coweight $\cla$ on $D_x$ as an
$N_+(\K_x)$--gauge equivalence class of operators of the form
\begin{equation}    \label{psi la new1}
\nabla = \pa_t + \sum_{i=1}^\ell t^{\langle \al_i,\cla \rangle} f_i +
{\mb v}(t),
\end{equation}
where ${\mb v}(t) \in \bb_+[[t]]$. Denote the set of $G$--opers of
coweight $\cla$ on $D_x$ by $\on{Op}_G(D_x)_{\cla} \subset
\on{Op}_G^{\on{RS}}(D_x)$.

\subsection{Miura opers and Cartan connections}

We generalize the results of \secref{miura opers} to the case of an
arbitrary Kac-Moody algebra.

Consider the flag variety $G/B_-$. This is an ind-scheme with the
ind-scheme structure defined as follows. As a set, $G/B_-$ decomposes
into a disjoint union of $B_-$--orbits parameterized by the Weyl group
$W$ of $G$. We denote the orbit corresponding to $w$ by $S^w$. These
orbits are finite-dimensional, and the closure of $S^w$ is the union
of the orbits $S^y$ corresponding to the elements $y \in W$ which are
less than or equal to $w$ with respect to the Bruhat order on $W$ (see
\cite{flags}, Ch. VII, for more details). Each of these closures,
$\ol{S}^w$, has the structure of a finite-dimensional (in general,
singular) algebraic variety. We have a collection of closed embeddings
$\ol{S}^y \hookrightarrow \ol{S}^w$ of these varieties into each other
corresponding to the Bruhat order. This collection defines the
structure of a (strict) ind-scheme on $G/B_-$.

We will say that $U \subset G/B_-$ is an open (resp., dense) subset if
for sufficiently large $w \in W$, with respect to the Bruhat order,
the intersection $U \cap S^w$ is open (resp., dense) in $U \cap S^w$.

The $B_+$--orbits in $G/B_-$ are also parameterized by the Weyl
group. We denote the $B_+$--orbit $B_+ w^{-1} B_- \subset G/B_-$ by
$S_w$, so that $S_1$ is the open dense orbit.

Let $(\F,\nabla,\F_{B_+},\F_{B_-})$ be a Miura $G$--oper on a curve
$X$. Then we have the following analogue of \lemref{isom with flags}.

\begin{lem}    \label{isom with flags1}
Suppose that the oper connection $\nabla$ gives rise to a global
trivialization of the bundle $\F_{B_-}$ (and hence of $\F$) on
$X$. Then for each $x \in X$ the set of horizontal $B_-$--reductions
of $\F$ over $X$ is canonically identified with the $\F_{x}$--twist of
the flag variety $G/B_-$, which coincides with its
$\F_{B_+,x}$--twists,
\begin{equation}    \label{twist of flag1}
(G/B_-)_{\F_x} = \F_x \us{G}\times G/B_- = \F_{B_+,x} \us{B_+}\times
G/B_- = (G/B_-)_{\F_{B_+,x}}.
\end{equation}
\end{lem}

We obtain from the second description of $(G/B_-)_{\F_x}$ given in
formula \eqref{twist of flag1} that $(G/B_-)_{\F_x}$ decomposes into a
union of the $\F_{B_+,x}$--twists of the $B_+$--orbits $S_w$ which we
denote by $S_{w,\F_{B_+,x}}$. We will say that $\F_{B_+,x}$ and
$\F_{B_-,x}$ are in {\em relative position} $w$ if $\F_{B_-,x}$,
considered as point of $(G/B_-)_{\F_x}$, belongs to $S_{w,\F_{B_+,x}}$
(this agrees with the definition given in \secref{miura opers} in the
finite-dimensional case). In particular, if it belongs to the open
orbit $S_{1,\F_{B_+,x}}$, we will say that $\F_{B_+,x}$ and
$\F_{B_-,x}$ are in generic position.

A Miura $G$--oper is called {\em generic} on $U \subset X$ if the
reductions $\F_{B_+,x}$ and $\F_{B_-,x}$ of $\F_x$ are in generic
position for all $x \in U$. We denote the set of generic Miura opers
on $U$ by $\on{MOp}_G(U)_{\on{gen}}$.

Consider the $H$--bundles $\F_H = \F_{B_+}/N_+$ and $\F'_H =
\F'_{B_-}/N_-$ corresponding to a generic Miura oper
$(\F,\nabla,\F_B,\F'_B)$ on $X$. Then we have the following result
(compare with \lemref{H bundles isom} in the finite-dimensional case):

\begin{lem}    \label{H bundles isom1}
For a generic Miura oper $(\F,\nabla,\F_{B_+},\F_{B_-})$ the
$H$--bundles $\F_H$ and $\F'_H$ are isomorphic.
\end{lem}

\begin{proof}
Since $\F_{B_+}$ and $\F_{B_-}$ are in generic position, their
intersection $\F_{B_+} \cap \F_{B_-}$ inside $\F$ is isomorphic to
both $\F_H$ and $\F'_H$. Hence we obtain that $\F_H \simeq \F'_H$.
\end{proof}

Since the $B_-$--bundle $\F_{B_-}$ is preserved by the oper connection
$\nabla$, we obtain a connection $\ol{\nabla}$ on $\F'_H$ and hence on
$\F_H$. We prove, in exactly the same way as in the proof of
\lemref{FH}, that $\F_H \simeq \Omega^{\crho}$, where $\crho$ is the
unique cocharacter $\C^\times \to H$ such that $\langle \al_i,\crho
\rangle = 1, i=1,\ldots,\ell$. Therefore we obtain a map ${\mb a}$
from the set of $\on{MOp}_G(U)_{\on{gen}}$ of generic Miura opers on
$U$ to the set of connections $\on{Conn}_U$ on the $H$--bundle
$\Omega^{\crho}$ on $U$.

Connections on $\Omega^{\crho}$ are described in the same way as in
the finite-dimensional case. If we choose a local coordinate $t$ on
$U$, then we trivialize $\Omega^{\crho}$ and represent the connection
as an operator $\pa_t + {\mb u}(t)$, where ${\mb u}(t)$ is an
$\h/\cc$--valued function on $U$. If $s$ is another coordinate such
that $t=\varphi(s)$, then this connection will be represented by the
operator
\begin{equation}    \label{trans for conn1}
\pa_s + \varphi'(s) {\mb u}(\varphi(s)) - \crho \cdot
\frac{\varphi''(s)}{\varphi'(s)}.
\end{equation}

\begin{prop}    \label{map beta1}
The map ${\mb a}: \on{MOp}_G(U)_{\on{gen}} \to \on{Conn}_U$ is an
isomorphism.
\end{prop}

\begin{proof}
We define a map ${\mb b}$ in the opposite direction, similarly to
the finite-dimen\-sional case. Suppose we are given a connection
$\ol\nabla$ on the $H$--bundle $\Omega^{\crho}$ on $D$. We associate
to it a generic Miura oper as follows. We set $\F = \Omega^{\crho}
\underset{H}\times G, \F_{B_\pm} = \Omega^{\crho} \underset{H}\times
B_\pm$, where we consider the adjoint action of $H$ on $G$ and on
$B_\pm$.

The space of connections on $\F$ is isomorphic to the direct product
$$
\on{Conn}_U \times \bigoplus_{\al \in \De} \Gamma(U,\Omega^{\al(\crho)
  + 1}).
$$
Its subspace corresponding to negative simple roots is isomorphic to
$\left( \bigoplus_{i=1}^\ell \g_{-\al_i} \right) \otimes R$. Having
chosen a basis element $f_i$ of $\g_{-\al_i}$ for each
$i=1,\ldots,\ell$, we now construct an element $p_{-1} =
\sum_{i=1}^\ell f_i \otimes 1$ of this space. Now we set $\nabla =
\ol\nabla + p_{-1}$. By construction, $\nabla$ has the correct
relative position with the $B_+$--reduction $\F_{B_+}$ and preserves
the $B_-$--reduction $\F_{B_-}$. Therefore the quadruple
$(\F,\nabla,\F_{B_+},\F_{B_-})$ is a generic Miura oper on $U$. We
define the map ${\mb b}$ by setting ${\mb b}(\ol{\nabla}) =
(\F,\nabla,\F_{B_+},\F_{B_-})$. It is clear that this map is
independent of the choice of the generators $f_i,
i=1,\ldots,\ell$, and that ${\mb a}$ and ${\mb b}$ are mutually
inverse maps.
\end{proof}

A {\em Miura $G$--oper of coweight} $\cla$ on $D_x$ is defined as a
quadruple $(\F,\nabla,\F_{B_+},\F_{B_-})$, where
$(\F,\nabla,\F_{B_+})$ is a $G$--oper on $D_x^\times$ which belongs to
$\on{Op}_G(D_x)_{\cla}$ and $\F_{B_-}$ is a $B_-$--reduction of $\F$
which is preserved by $\nabla$. We denote the set of Miura $G$--opers
of coweight $\cla$ on $D_x$ by $\on{MOp}_G(D_x)_{\cla}$. In
particular, if $\cla=0$, then we obtain the old definition of Miura
opers. All of the above definitions and results can be easily carried
over to the case of an arbitrary integral $\cla$.

\subsection{Bethe Ansatz equations and Miura opers}

Now we establish a connection between the Bethe Ansatz equations and
Miura opers on $\pone$, following \secref{miura opers on pone}.

Let us fix, as in the finite-dimensional case, a set of distinct
complex numbers $z_1,\ldots,z_N$ (which we will view as points of
$\pone \bs \infty$) and a set of dominant integral coweights
$\cla_1,\ldots,\cla_N \in \h/\cc$ (a dominant integral coweight is by
definition an element $\cla$ of $\h/\cc$ such that $\langle \al_i,\cla
\rangle \in \Z_+$ for all $i=1,\ldots,\ell$).

In what follows we will consider the Miura opers on $\pone$ rather
than opers. The reason is that, as we explained above, an oper
connection does not in general allow us to identify the nearby fibers
of the oper bundle; this is because the group $G$ is an
ind-group. However, if in addition to an oper structure we are given a
horizontal $B_-$--reduction $\F_{B_-}$, i.e., if we are given the
structure of a Miura oper, then, because $B_-$ is a proalgebraic
group, we can identify nearby fibers of $\F_{B_-}$ (and hence of $\F$)
using the oper connection. This makes Miura opers much easier to
handle.

Let $\on{MOp}_G(\pone)_{(z_i);(\cla_i)}$ be the set of all Miura
$G$--opers $(\F,\nabla,\F_{B_+},\F_{B_-})$ on $\pone$ without the
points $z_1,\ldots,z_N,\infty$, whose restrictions to the punctured
discs $D_{z_i}^\times$ around the points $z_i$ belong to
$\on{MOp}_{G}(D_{z_i})_{\cla_i}$, and such that the restriction of the
underlying oper to the disc around $\infty$ has regular
singularity. Then the connection induced by $\nabla$ on the
$B_-$--bundle $\F_{B_-}$ is regular everywhere on $\af$. Since $B_-$
is a proalgebraic group, we obtain that the connection identifies the
fibers of the bundle $\F_{B_-}$ (and hence of $\F$) at all points of
the affine line. Let us trivialize the fiber of $\F_{B_-}$ at some
point $z_0$ of $\af$. Then we obtain a trivialization of the bundle
$\F$ over $\af$. The oper reduction $\F_{B_+}$ gives rise to a map
$\phi: \af \to G/B_+$. Two opers underlying Miura opers from
$\on{MOp}_G(\pone)_{(z_i);(\cla_i)}$ are isomorphic if and only if the
corresponding maps $\phi$ differ by an element $g \in G$ acting on
$G/B_+$.

On the other hand, the choice of horizontal reduction $\F_{B_-}$
corresponds to the choice of a reduction to $B_-$ in the fiber of $\F$
at $z_0$ (see \lemref{isom with flags1}). Since we have trivialized
this fiber, the latter is nothing but a point of $G/B_-$. Thus, we
obtain that the space $\on{MOp}^\phi_G(\pone)_{(z_i);(\cla_i)}$ of all
Miura opers with the underlying oper map $\phi$ is isomorphic to
$G/B_-$.

Recall that any pair of points $yB_+ \in G/B_+$ and $pB_- \in G/B_-$
have a well-defined relative position. Namely, we will say that they
have relative position $w \in W$ if $p B_- \in (yB_+y^{-1}) w^{-1}
B_-$.  We will say that they are in generic position if $w=1$.

Consider the subset $(G/B_-)_\phi$ of $G/B_-$ whose points $pB_-$
satisfy the following conditions (as in \secref{miura opers on pone}):
\begin{itemize}
\item[(1)] $\phi(z_i)$ is in generic position with
$pB_-$ for all $i=1,\ldots,N$;

\item[(2)] the relative position of $\phi(x)$ and $pB_-$ is either
generic or corresponds to a simple reflection $s_i \in W$ for all $x
\in \af \bs \{ z_1,\ldots,z_N \}$.
\end{itemize}

As in the finite-dimensional case, it is clear that $(G/B_-)_\phi$ is
an open and dense subset of $G/B_-$. We claim that there is a
bijection between this set and the set of solutions of the Bethe
Ansatz equations \eqref{bethe1}
\begin{equation}    \label{bethe1}
\sum_{i=1}^N \frac{\langle \al_{i_j},\cla_i \rangle}{w_j-z_i} -
\sum_{s \neq j} \frac{\langle \al_{i_j},\chal_{i_s} \rangle}{w_j-w_s}
= 0, \qquad j=1,\ldots,m.
\end{equation}
As in the finite-dimensional case, we have an obvious action of a
product of symmetric groups permuting the points $w_j$ corresponding
to simple roots of the same kind. As before, by a {\em solution} of
the Bethe Ansatz equations we will understand a solution defined up to
these permutations. We will also adjoin to the set of all solutions
associated to all possible collections $\{ \al_{i_j} \}$ of simple
roots of $\g$, the unique ``empty'' solution, corresponding to the
empty set of simple roots.

We start with the following

\begin{lem}    \label{lem km}
Suppose that we are given a regular Miura oper on the disc $D_x$ such
that $\F_{B_+,x}$ and $\F_{B_-,x}$ are in relative position
$s_i$. Then the oper connection on $D^\times_x$ may be brought to the
form
$$
\pa_t + p_{-1} + \frac{\chal_i}{t} + {\mb u}(t), \qquad {\mb u}(t) \in
\h[[t]]
$$
(with respect to a coordinate $t$ at $x$), where $\langle \al_i,{\mb
u}(0) \rangle = 0$.
\end{lem}

\begin{proof}
First, we observe that the two reductions $\F_{B_+}$ and $\F_{B_-}$
are in generic relative position on the punctured disc
$D^\times_x$. Indeed, let $V^l_{-\omega_i}$ be the lowest weight
integrable module over $\g$ with lowest weight $-\omega_i$. Consider
its two-dimensional submodule over the $\sw_2$ subalgebra
corresponding to the $i$th simple root, generated by a lowest weight
vector $v_{-\omega_i}$. Let $v_{-\omega_i+\al_i} = e_i
v_{-\omega_i}$. By our assumption, the oper connection has the form
$$
\nabla = \pa_t + \sum_{i=1}^\ell f_i + {\mb v}(t), \qquad {\mb v}(t)
\in \bb_+[[t]].
$$
Let $\Phi(t)$ be the unique solution of the equation $\nabla \Phi(t) =
0$ such that $\Phi(0) = 1$. Such a solution exists by our assumption that
our oper carries a horizontal $B_-$--reduction. The above statement is
equivalent to the assertion that $\Phi(t) \cdot v_{-\omega_i+\al_i}
\in V^l_{-\omega_i}[[t]]$ is a linear combination of weight vectors 
which contains the lowest weight vector $v_{-\omega_i}$ with a non-zero
coefficient. But this follows immediately from the observation that
$f_i v_{-\omega_i+\al_i} = - v_{-\omega_i}$.

This further implies, in the same way as in the proof of \propref{isom
w} that by gauge transformation with an element of $N_+((t))$ we can
bring the oper connection $\nabla$ to the form
$$
\pa_t + p_{-1} + \frac{\chal_i}{t} + {\mb u}(t), \qquad {\mb u}(t) \in
\h[[t]].
$$
We associate to this Miura $G$--oper a Miura $SL_2$--oper in the same
way as in \lemref{si}. Following the argument used in the proof of
\lemref{si}, we find that the monodromy of this oper is non-trivial
unless $\langle \al_i,{\mb u}(0) \rangle = 0$. This completes the
proof.
\end{proof}

Now consider the Miura oper in
$\on{MOp}^\phi_G(\pone)_{(z_i);(\cla_i)}$ corresponding to a point of
the subset $(G/B_-)_\phi$. Then by \lemref{lem km}, the corresponding
connection operator may be brought to the form $\pa_t + p_{-1} + {\mb
u}(t)$, where ${\mb u}(t)$ is a rational function on $\af$ which is
regular at all points other than $z_1,\ldots,z_N$, and whose expansion
at $z_i$ has the form
$$
\pa_t + p_{-1} - \frac{\cla_i}{t-z_i} + \on{reg}.,
$$
and the expansion at $w_j$ has the form
$$
\pa_t + p_{-1} + \frac{\chal_{i_j}}{t-w_j} + {\mb u}_j(t-w_j), \qquad
   {\mb u}_j(t-w_j) \in \h[[t-w_j]],
$$
and $\langle \al_{i_j},{\mb u}_j(0) \rangle = 0$. By our assumption, our
oper has regular singularity at $\infty$, which implies that
\begin{equation}    \label{new miura}
\nabla = \pa_t + p_{-1} - \sum_{i=1}^N \frac{\cla_i}{t-z_i} +
\sum_{j=1}^m \frac{\chal_{i_j}}{t-w_j}.
\end{equation}
The condition $\langle \al_{i_j},{\mb u}_j(0) \rangle = 0$ from
\lemref{lem km} is precisely the $j$th Bethe Ansatz equation
\eqref{bethe1}. Thus, we obtain a map from $(G/B_-)_\phi$ to the set
of solutions of equations \eqref{bethe1}.

Let us construct the inverse map. Given a solution of the Bethe Ansatz
equations, we define a Miura $G$--oper in
$\on{MOp}_G(\pone)_{(z_i);(\cla_i)}$. We set $\F = \Omega^{\crho}
\underset{H}\times G, \F_{B_\pm} = \Omega^{\crho} \underset{H}\times
B_\pm$ and define the connection operator by formula \eqref{new
miura}. Clearly, the two reductions $\F_{B_\pm}$ satisfy the
conditions of a Miura oper. It remains to show that its restriction to
the punctured disc at $z_i$ (resp., $w_j$) belongs to
$\on{MOp}_{G}(D_{z_i})_{\cla_i}$ (resp., $\on{MOp}_{G}(D_{w_j})$), and
that it has regular singularity at $\infty$.

The expansion at $z_i$ of the connection \eqref{new miura} reads
$$
\pa_t + p_{-1} - \frac{\cla_i}{t-z_i} + \on{reg}.,
$$
which after conjugation by $\cla_i(t-z_i)^{-1}$ becomes
$$
\pa_t + \sum_{k=1}^\ell (t-z_i)^{\langle \al_k,\cla_i \rangle} f_k +
\on{reg}.
$$
Therefore the restriction to the punctured disc at $z_i$ belongs to
$\on{MOp}_{G}(D_{z_i})_{\cla_i}$ as desired.

Next, consider the expansion at the point $w_j$. We find that it has
the form
$$
\pa_t + p_{-1} + \frac{\chal_{i_j}}{t-w_j} + {\mb u}_j(t-w_j), \qquad
   {\mb u}(t-w_j) \in \h[[t-w_j]].
$$
Moreover, we find that $\langle \al_{i_j},{\mb u}_j(0) \rangle$ is
given by the expression appearing in the $j$th Bethe Ansatz
equation. The Bethe Ansatz equation means that $\langle \al_{i_j},{\mb
u}(0) \rangle = 0$, which ensures that this connection becomes
regular after conjugation with $\exp(-e_{i_j}/(t-w_j))$ (compare with
\lemref{si} in the finite-dimensional case). Therefore the restriction
to the punctured disc at $w_j$ belongs to $\on{MOp}_{G}(D_{w_j})$ as
desired.

Finally, using the transformation formula for the oper connection that
is identical to the one obtained in the finite-dimensional case (see
formula \eqref{change of var}), we find the restriction of the oper
\eqref{new miura} to the punctured disc at $\infty$:
\begin{equation}    \label{exp infty}
\pa_u + p_{-1} + u^{-1} \left( \sum_{i=1}^N \cla_i - \sum_{j=1}^m
\chal_{i_j} + 2 \crho \right) + \on{reg}.,
\end{equation}
where $u = t^{-1}$. Thus, it has regular singularity at $\infty$

Thus, we obtain a bijection between the set of solutions of the Bethe
Ansatz equations and the union of the sets of points of open dense
subsets $(G/B_-)_\phi$ of the flag variety $G/B_-$, just as in the
finite-dimensional case. Now we show that the residues of the
connection at $\infty$ correspond to the $B_-$--orbits in $G/B_-$.

Consider the action of the group $N_+$ on $G/B_-$. It translates into
a rational action of $N_+$ on the set of solutions of the Bethe Ansatz
equations. Let $SL_2^{(i)}$ be the $SL_2$ subgroup of $G$
corresponding to the $i$th simple root and set $N_{\al_i} = N_+ \cap
SL_2^{(i)}$, $B_{\al_i} = B_- \cap SL_2^{(i)}$. Note that $N_{\al_i}$
is the one-parameter additive subgroup $\{ \exp(ae_i) \}_{a \in \C}
\subset N_+$.

Observe that the $SL_2^{(i)}$--orbits in $G/B_-$ give us a partition
of $G/B_-$ into a disjoint union of $\pone \simeq
G/B_{\al_i}$. Furthermore, if $p B_-$ is a point in the Schubert cell
$B_- y B_- \subset G/B_-$, then there are two possibilities. The first
case is that the intersection of the $SL_2^{(i)}$--orbit passing
through this point and $B_- y B_-$ is an affine line. Then $l(s_i y) <
l(y)$ and the remaining point of this $SL_2^{(i)}$--orbit belongs to
the smaller Schubert cell $B_- s_i y B_-$ which is in the closure of
$B_- y B_-$; this point is then stable under $B_{\al_i}$. The second
case is that this intersection is the point $p B_-$, which is
therefore stable under $B_{\al_i}$. Then $l(s_i y) > l(y)$ and the
remaining part of the $SL_2^{(i)}$--orbit belongs to the larger
Schubert cell $B_- s_i y B_-$ which contains $B_- y B_-$ in its closure.

On the other hand, the action of $N_{\al_i}$ on solutions of the Bethe
Ansatz equations may be computed explicitly as in \secref{action of
N}. Namely, we find that the element $\exp(ae_i)$ acts on the
connection $\pa_t + p_{-1} + {\mb u}(t)$ by sending ${\mb u}(t)$ to
$\wt{\mb u}(t) = {\mb u}(t) + f(t) \chal_i$, where $f(t)$ is the
solution of the equation
\begin{equation}    \label{diff eq1}
f'(t) + f(t)^2 + f(t) u_i(t) = 0, \qquad u_i(t) = \langle \al_i,{\mb
u}(t) \rangle,
\end{equation}
with the initial condition $f(0) = a$. Now observe that
\begin{equation}    \label{ui}
u_i(t) = - \sum_{k=1}^N \frac{\langle \al_i,\cla_k \rangle}{t-z_i} +
\sum_{j \not\in S_i} \frac{\langle \al_i,\chal_j \rangle}{t-w_j} +
\sum_{j \in S_i} \frac{2}{t-w_j},
\end{equation}
where $S_i \subset \{ 1,\ldots,m \}$ is the set of those $j$'s for
which $i_j = i$. Hence \eqref{diff eq1} looks exactly like the
corresponding equation for the action of $\exp(a e)$ on the solution
of the Bethe Ansatz equation in the case of $\sw_2$ corresponding to
the connection $\pa_t + p_{-1} + u_i(t)$, where $u_i(t)$ is given by
\eqref{ui}.

This solution corresponds to the situation where we have dominant
coweights $\langle \al_i,\cla_k \rangle$ of $\sw^{(i)}_2$ attached to
the point $z_k$ for $k=1,\ldots,N$, dominant coweights $\langle
\al_i,\chal_{i_j} \rangle$ attached to the points $w_j$ with $j \in \{
1,\ldots,m \} \bs S_i$, and the variables of the Bethe Ansatz
equations are $w_j, j \in S_i$. Clearly, the Bethe Ansatz equations
\eqref{bethe1} with $j \in S_i$ imply that these $w_j$'s, $j \in S_i$
indeed solve the Bethe Ansatz equations for $\sw^{(i)}_2$ in the above
situation. Hence we find that the action of $\exp(a e_i)$ on our
solution can be read off of the action of $\exp(a e)$ on the
corresponding solution of the Bethe Ansatz equation for $\sw_2$. But
we know from \corref{fixed oper} and the discussion of \secref{action
of N} that the latter corresponds to the action of the unipotent
subgroup of $SL_2$ on the flag manifold $SL_2/B_- = \pone$. The
closure of the orbit of any solution of the $SL_2$ Bethe Ansatz
equation under the action of $N$ coincides with this $\pone$. Further,
it has two $B_-$--orbits: a point and an affine line. Suppose that the
one point orbit belongs to the open subset of $SL_2/B_-$ of points
corresponding to the solutions of the Bethe Ansatz equations (this
will be so for a generic collection of points $z_1,\ldots,z_N$). Using
\corref{fixed oper}, we find that then this point corresponds to the
unique solution for which
\begin{equation}    \label{domi}
\left\langle \al_i,\sum_{k=1}^N \cla_k - \sum_{j=1}^m \chal_{i_j}
\right\rangle
\end{equation}
is a non-negative integer; we denote this integer by $n_i$. The
points of the other, one-dimensional, cell correspond to a
one-parameter family of solutions for which the number \eqref{domi} is
a negative integer equal to $-n_i-2$. Moreover, the number \eqref{domi}
is always an integer and is never equal to $-1$.

For any solution of the Bethe Ansatz equation call the expression
\begin{equation}    \label{at infty}
\sum_{k=1}^N \cla_k - \sum_{j=1}^m \chal_{i_j}
\end{equation}
the {\em residue of the solution at $\infty$}. Note that it can be
obtained from the expansion \eqref{exp infty} of the connection
\eqref{new miura} around $\infty$.

The above analysis leads us to the following conclusion.

\begin{lem}    \label{acting with ei}
Let $p B_-$ be a point of $(G/B_-)_\phi$ which belongs to the Schubert
cell $B_- y B_-$. Denote the residue of the corresponding solution of
the Bethe ansatz equations by $\cmu_\infty$. Then $\langle
\al_i,\cmu_\infty \rangle$ is an integer not equal to $-1$. It is
non-negative if and only if $l(s_i y) > l(y)$ and negative if and only
if $l(s_i y) > l(y)$.
\end{lem}

This lemma implies that we may set up our bijection between solutions
of the Bethe Ansatz equations corresponding to a fixed oper and the
set of points of an open dense subset of $G/B_-$ in such a way that
the residue of the solution at $\infty$ is always equal to
\begin{equation}    \label{constraint}
\sum_{k=1}^N \cla_k - \sum_{j=1}^m \chal_{i_j} =
y(\cla_\infty+\crho)-\crho
\end{equation}
for some dominant integral coweight $\cla_\infty$ and $y \in W$, and
in this case the corresponding point of $G/B_-$ belongs to the
$B_-$--orbit $B_- y B_-$ of $G/B_-$.

Consider first the simplest case when $\cla_\infty = \sum_{k=1}^N
\cla_k$. Then to $y=1$ corresponds the ``empty'' solution of the Bethe
Ansatz equations, when the set of the $w_j$'s is empty. This solution
has residue $\sum_{k=1}^N \cla_k$ and hence indeed corresponds to the
one-point $B_-$--orbit $B_- \in G/B_-$. Now we apply induction on the
length $l(y)$ of $y$. Suppose we have proved the result for all $y$
whose length is less than or equal to $N$. Let us prove it for those
elements whose length is $N+1$. Those may be written in the form $s_i
y$, where $s_i$ runs over the list of all simple reflections which
satisfy $l(s_iy) > l(y)$. This is equivalent to the following
property: for any dominant integral coweight $\cla_\infty$ we have
$\langle \al_i,y(\cla_\infty+\crho)-\crho \rangle = n_i \in
\Z_+$. Then $\langle \al_i,s_iy(\cla_\infty+\crho)-\crho \rangle =
-n_i-2$. In this case the union of the $B_-$--orbits $B_- y B_-$ and
$B_- s_i y B_-$ is the union of the closures of the orbits $N_{\al_i}
\cdot g B_-, g B_- \in B_- y B_-$. Our inductive assumption and the
above computation then shows that the solutions corresponding to the
points of $B_- s_i y B$ have residue $s_iy(\cla_\infty+\crho)-\crho$.

Consider now the open subset $(G/B_-)_\phi$ corresponding to a general
oper $\phi$. Suppose that the one-point $B_-$--orbit $B_- \in G/B_-$
belongs to this subset. Then we claim that the residue $\cmu_\infty$
of the corresponding solution of the Bethe Ansatz equations is a
dominant coweight, i.e., $y=1$ in formula \eqref{constraint}. Indeed,
for each $i=1,\ldots,\ell$ the point $B_- \in G/B_-$ is the one point
$B_{\al_i}$--orbit in the $SL_2^{(i)}$--orbit passing through $B_-$
for all $i=1,\ldots,\ell$. Hence it follows from \lemref{acting with
ei} that $\langle \al_i,\cmu_\infty \rangle \geq 0$ for all
$i=1,\ldots,\ell$. Next, we consider the solutions whose residue
belongs to the orbit of $\cmu_\infty$ under the action of the Weyl
group. Using induction on the length in the same way as above, we
obtain that the points that belong to $B_- y B_- \cap (G/B_-)_\phi$
correspond to solutions satisfying \eqref{constraint}.

Finally, suppose that the one-point $B_-$--orbit $B_- \in G/B_-$ does
not belong to $(G/B_-)_\phi$. Then we pick a point of $(G/B_-)_\phi$
that belongs to the Schubert cell of the smallest possible
dimension. Consider the closures of the orbits of this point under the
action of the subgroups $SL_2^{(i)}$. Then we consider the
$SL_2^{(i)}$--orbits of the points obtained this way, and so on. As
the result, we can reach any point of $G/B_-$ in finitely many
steps. According to \eqref{acting with ei}, each time we cross from a
smaller Schubert cell $B_- y B_-$ to a larger one $B_- s_iy B_-$ via
the $SL_2^{(i)}$--orbit, the residue of the corresponding solutions of
the Bethe Ansatz equations changes from being dominant with respect to
the $i$th simple root, say, $n_i \in \Z_+$, to being anti-dominant
$-2-n_i$, while the pairing with the roots $\al_j, j \neq i$, remains
unchaged. Hence consistency with \lemref{acting with ei} requires that
the solutions corresponding to the points of $B_- y B$ have
residue $y(\cla_\infty+\crho)-\crho$.

Therefore we obtain the following result.

\begin{thm}    \label{final km}
There is a bijection between the set of solutions of the Bethe Ansatz
equations corresponding to the same underlying $G$--oper $\phi$ and an
open dense subset of the ind-flag variety $G/B_-$ such that the set
of solutions which satisfy \eqref{constraint} is in bijection with an
open subset of the $B_-$--orbit $B_- y B_- \subset G/B_-$.
\end{thm}

We also have an analogue of \thmref{strongest} establishing a
bijection between the set of all points of $G/B_-$ and a certain set
of connections on the $H$--bundle $\Omega^{\crho}$ over $\af$ of the
form
$$
\pa_t - \sum_{i=1}^N \frac{y_i(\cla_i+\crho)-\crho}{t-z_i}
- \sum_{j=1}^m \frac{y'_j(\crho)-\crho}{t-w_j}.
$$

\end{document}